\magnification=\magstep1
\input amstex
\documentstyle{amsppt}
\loadbold

\define\defeq{\overset{\text{def}}\to=}
\define\Sigmadagbyp{\Sigma^{\dag}\text{-by-}\{p\}}

\define\tor{\operatorname{tor}}
\define\Spec{\operatorname{Spec}}
\define\Gal{\operatorname{Gal}}
\define\ab{\operatorname{ab}}
\define\Aut{\operatorname{Aut}}
\define\cl{\operatorname{cl}}

\define\pr{\operatorname{pr}}

\define\sep{\operatorname{sep}}
\define\Sub{\operatorname{Sub}}

\define\ur{\operatorname{ur}}

\def \Div{\operatorname {Div}}

\def \ab{\operatorname {ab}}
\def \isom {\buildrel \sim \over \rightarrow}
\def \Ker{\operatorname {Ker}}
\def \Hom{\operatorname {Hom}}
\def \char{\operatorname {char}}
\def \id{\operatorname {id}}

\def \mod{\operatorname {mod}}

\define\Primes{\frak{Primes}}
\def \t{\operatorname {t}}
\def \ur{\operatorname {ur}}

\def \Im{\operatorname {Im}}
\def \PDiv{\operatorname {PDiv}}
\def \deg{\operatorname {deg}}
\def \gon{\operatorname {gon}}

\def \Dec{\operatorname {Dec}}
\def \cd{\operatorname {cd}}

\NoRunningHeads
\NoBlackBoxes

\document

\head
A refined version of Grothendieck's birational anabelian conjecture 
for 
curves over finite fields
\endhead
\bigskip
\centerline {MOHAMED SA\"IDI and AKIO TAMAGAWA}

\bigskip
\definition{Abstract}
In this paper we prove a 
refined version of Uchida's theorem 
on isomorphisms 
between absolute Galois groups of global fields in positive characteristics,
where one ``ignores'' the information provided by a ``small'' set of primes.
\enddefinition

\bigskip

\noindent
\S 0. Introduction

\noindent
Part I.

\noindent
\S 1.  Generalities on Galois Groups of Function Fields of Curves

\noindent
\S 2. Isomorphisms between Geometrically Pro-$\Sigma$ Galois Groups

\noindent
Part II.

\noindent
\S 3. Small and Large Sets of Primes

\noindent
\S 4. The Main Theorem



\subhead
\S 0. Introduction
\endsubhead
Let $k$ be a finite field of characteristic $p>0$ and $X$ a proper, smooth, and geometrically connected 
algebraic curve over $k$. Let 
$
K$ be the function field of $X$, with separable closure 
$
K^{\sep}$, 
and let $\bar k$ be the algebraic closure of $k$ in $K^{\sep}$.
We have the following exact sequence of profinite groups:
$$1\to {\overline G}_K\to G_{K}  @>\text {pr}>> G_k\to 1.$$
Here, $G_k$ is the absolute Galois group $\Gal (\bar k/k)$ of $k$, 
$G_{K}$ is 
the absolute Galois group $\Gal (K^{\sep}/K)$ 
of 
$
K$, and ${\overline G}_K$ is the absolute Galois group
$\Gal (K^{\sep}/K\bar k)$ of the function field $K\bar k$ of 
$\overline X\defeq X\times _k\bar k$. The following result is fundamental
in the birational anabelian geometry of curves over finite fields.

\proclaim{Theorem A (Uchida)} Let $X$, $Y$ be proper, smooth, and geometrically 
connected curves over finite fields $k$, $l$, respectively. 
Let $K$, $L$ be the function fields of $X$, $Y$, respectively. 
Let $G_{K}=\Gal (K^{\sep}/K)$, $G_{L}=\Gal (L^{\sep}/L)$ be the absolute Galois groups of $K$, $L$, 
respectively. Let 
$$\sigma : G_{K}\isom G_{L}$$
be an isomorphism of profinite groups. Then $\sigma$ arises from a uniquely 
determined commutative diagram of field extensions:
$$
\CD
L^{\sep} @>{\sim}>> K^{\sep} \\
@AAA              @AAA \\
L @>{\sim}>> K \\
\endCD
$$
in which the horizontal arrows are isomorphisms, and the vertical arrows are 
the natural field extensions. 
\endproclaim

This theorem was proved by Uchida [Uchida]. A stronger result involving fundamental groups
of hyperbolic curves over finite fields was proved by Tamagawa [Tamagawa] and Mochizuki [Mochizuki]
(see also [Sa\"\i di-Tamagawa2] for a survey of recent results in the 
anabelian geometry of hyperbolic curves over finite fields).
Uchida's theorem implies in particular that one can embed a suitable category of curves over finite fields 
into the category of profinite groups
via the absolute Galois group functor. 
It is essential in the anabelian philosophy of Grothendieck, as was 
formulated in [Grothendieck], to be able to determine the image of this 
functor. Recall that the full structure of the absolute Galois group $G_{K}$ is unknown, 
though one knows the structure of the closed  subgroup  ${\overline G}_K$ of 
$G_{K}$ by a result of Pop and Harbater. Namely ${\overline G}_K$
is a free profinite group on 
countably infinitely many generators 
(cf. [Pop1], [Harbater]),
though one has no precise description of a free set of generators of
${\overline G}_K$.
Thus, the problem of determining the image of the above functor 
seems to be quite difficult, at least for the moment. 
It is quite natural to address the following question:

\definition{Question 1} 
Is it possible to prove  any result analogous to Theorem A where $G_{K}$ is replaced 
by some (continuous) quotient of $G_{K}$ whose structure is better understood?
\enddefinition

The first quotients that come into mind are the following. 
Let $\Primes$ denote the set of all prime numbers.
Let $\Sigma \subset \Primes$ be a set of prime numbers not containing 
the characteristic $p$. Let $\Cal C$ be the full class of finite groups whose cardinality 
is divisible only by primes in $\Sigma$. Let ${\overline G}_K^{\Sigma}$ be the 
maximal pro-$\Cal C$ quotient of ${\overline G}_K$. 
Here, the structure of ${\overline G}_K^{\Sigma}$ is well understood: 
${\overline G}_K^{\Sigma}$ is isomorphic to the projective limit of 
the maximal pro-$\Sigma$ quotients $\pi_1(\overline U)^{\Sigma} $ of the fundamental groups 
$\pi_1(\overline U)$, where $\overline U$ runs over all non-empty open subschemes of $\overline X$, and
$\pi_1(\overline U)^{\Sigma} $ is isomorphic to the pro-$\Sigma$ completion
of a certain well-known finitely generated discrete group (i.e., either a free 
group or a surface group). 

Let $G_{K}^{(\Sigma)}\overset \text {def}\to= G_{K}/\Ker({\overline G}_K\twoheadrightarrow 
{\overline G}_K^{\Sigma})$. (Note that $\Ker({\overline G}_K\twoheadrightarrow 
{\overline G}_K^{\Sigma})$ is a normal subgroup of $G_K$ since it is a characteristic subgroup of
${\overline G}_K$.) We shall refer to 
$G_{K}^{(\Sigma)}$ as the maximal 
geometrically pro-$\Sigma$ quotient of  the absolute Galois group $G_{K}$ 
(or, in short, the geometrically pro-$\Sigma$ Galois group of $K$). 

\definition{Question 2} 
Is it possible to prove any result analogous to Theorem A
where $G_{K}$ is replaced by $G_{K}^{(\Sigma)}$, for
a given set of prime numbers $\Sigma\subset \Primes$ (not containing the characteristic $p$)?
\enddefinition

The first set $\Sigma$ to consider is the set $\Sigma\defeq \Primes 
\setminus \{\char(k)\}$. In this case we shall refer to 
$G_{K}^{(\prime)}\defeq G_{K}^{(\Sigma)}$ as the maximal 
geometrically prime-to-characteristic quotient
of the absolute Galois group $G_{K}$. We have the following result which was proved by Sa\"\i di and Tamagawa 
(cf. [Sa\"\i di-Tamagawa1], Corollary 3.10).

\proclaim 
{Theorem B (Prime-to-$p$ Version of Uchida's Theorem)} 
Notations as in Theorem A, 
let $G_{K}^{(\prime)}$, $G_{L}^{(\prime)}$ be the maximal geometrically 
prime-to-characteristic quotients of $G_K$, $G_L$, respectively. Let 
$$\sigma : G_{K}^{(\prime)}\isom G_{L}^{(\prime)}$$
be an isomorphism of profinite groups. Then $\sigma$ arises from a uniquely 
determined commutative diagram of field extensions:
$$
\CD
L^{(\prime)} @>{\sim}>> K^{(\prime)} \\
@AAA              @AAA \\
L @>{\sim}>> K \\
\endCD
$$
in which the horizontal arrows are isomorphisms, and the vertical arrows are the 
extensions corresponding to the Galois groups $G_{L}^{(\prime)}$,  $G_{K}^{(\prime)}$, respectively.
Thus, $L^{(\prime)}/L$ (resp. $K^{(\prime)}/K$) is the subextension of $L^{\sep}/L$ (resp. $K^{\sep}/K$)
with Galois group $G_{L}^{(\prime)}$ (resp. $G_{K}^{(\prime)}$).
\endproclaim

Let $\Sigma \subset \Primes$ be a set of primes, 
and set $\Sigma'\defeq\Primes\setminus\Sigma$. 
We say that $\Sigma$ is {\it $k$-large} if the following condition is satisfied: the $\Sigma'$-cyclotomic character 
$\chi _{\Sigma'}:G_{k}\to 
\prod _{\ell\in \Sigma'\setminus\{p\}}\Bbb Z_{\ell}^{\times}$ is not injective ($p=\char(k)$). We say that  
$\Sigma$ satisfies the condition $(\epsilon_X)$ if the following holds:

\medskip
$(\epsilon_X)$: 
For each finite extension $k'$ of $k$ in $\bar k$, 
there exists an (infinite or finite) 
extension $k''$ of 
$k'$ in $\bar k$, such that $2\sharp (J_X(k'')\{\Sigma'\})<\sharp(k'')
\leq\infty
$ 
(hence, in 
particular, 
$\sharp (J_X(k'')\{\Sigma'\})<\infty$), 

\medskip\noindent
where $J_X$ denotes the jacobian 
variety of $X$ over $k$ and $J_X(k'')\{\Sigma'\}$ denotes the 
$\Sigma'$-primary part of the torsion group 
$J_X(k'')$.

Sets of prime numbers $\Sigma\subset \Primes$ which are $k$-large and satisfy the condition $(\epsilon _X)$
(in the above sense) include those such that 
$\Primes \setminus \Sigma$ is a finite set. There exist sets of primes $\Sigma$ which are $k$-large and satisfy
$(\epsilon _X)$ such that $\Primes \setminus \Sigma$ is an infinite set. However, a finite set of prime numbers 
is never $k$-large.

Our main result in this paper is the following refined version of the above Theorems A and B.

\proclaim{Theorem C (A Refined Version of Uchida's Theorem)} 
Notations as in Theorem A, let 
$\Sigma_X, \Sigma_Y\subset \Primes$ be sets of primes. 
Assume that $\Sigma_X$ is $k$-large and satisfies the condition $(\epsilon_X)$. 
Let  $G_{K}^{(\Sigma_X)}$ (respectively, $G_{L}^{(\Sigma_Y)}$) be the maximal 
geometrically pro-{$\Sigma_X$} quotient
of  $G_{K}$ (respectively, 
the maximal geometrically pro-{$\Sigma_Y$} quotient
of  $G_{L}$). Let 
$$\sigma : G_{K}^{(\Sigma_X)}\isom G_{L}^{(\Sigma_Y)}$$
be an isomorphism of profinite groups. Then $\sigma$ arises from a uniquely 
determined commutative diagram of field extensions:
$$
\CD
L\sptilde @>{\sim}>> K\sptilde \\
@AAA              @AAA \\
L @>{\sim}>> K \\
\endCD
$$
in which the horizontal arrows are isomorphisms and the vertical arrows are 
the field extensions corresponding to the Galois groups $G_{L}^{(\Sigma_Y)}$,  
$G_{K}^{(\Sigma_X)}$, 
respectively. Thus, $L\sptilde/L$ (resp. $K\sptilde/K$) is the subextension of $L^{\sep}/L$ (resp. $K^{\sep}/K$)
with Galois group $G_{L}^{(\Sigma_Y)}$ (resp. $G_{K}^{(\Sigma_X)}$).
\endproclaim


\definition{Note} When the authors announced the result of the present paper 
in [Sa\"\i di-Tamagawa2] (cf. loc. cit. Theorem 1.5), they overlooked the 
necessity to assume condition $(\epsilon_X)$. For the time being, they do not know if one could remove 
this extra assumption in general. (It is not difficult to see that we can remove it at least when 
the genus of $X$ is $\leq 1$.) 
\enddefinition

{\it Strategy of Proof.}\ 
In what follows  
we explain the steps/ideas of the proof. 

{\it Step 1.}\ Starting from an isomorphism
$$\sigma : G_{K}^{(\Sigma_X)}\isom G_{L}^{(\Sigma_Y)}$$
between profinite groups, one can first, using well-known results on the group-theoretic 
characterization
of decomposition groups in Galois groups (the so-called local theory), 
establish a set-theoretic 
bijection
$$\phi: X^{\cl}\isom Y^{\cl}$$
between the sets of closed points of $X$, $Y$, respectively, such that $\sigma (D_x)=D_{\phi (x)}$
where $D_x$, $D_{\phi (x)}$ denote the decomposition groups of $x$, $\phi (x)$ in 
$G_{K}^{(\Sigma_X)}$, $G_{L}^{(\Sigma_Y)}$, respectively (which are 
only defined up to conjugation).

{\it Step 2.}\ It is not difficult to prove 
that $p\defeq \char(k)=\char(l)$, 
that $\Sigma\defeq\Sigma_X=\Sigma_Y$, and 
that $\Sigma$ is both $k$-large and $l$-large and satisfies both $(\epsilon_X)$ and $(\epsilon_Y)$.

{\it Step 3.}\ Using global class field theory
(one could also use Kummer theory in this step) one can reconstruct, 
naturally from $\sigma$, 
finite index subgroups $\overline H_K$, $\overline H_L$ 
of the groups of principal divisors $K^{\times}/k^{\times}$, $L^{\times}/l^{\times}$, 
respectively, 
finite index subgroups $H'_K$, $H'_L$ 
of the multiplicative groups 
$(K^{\times})^{(\Sigma)}\defeq K^{\times}/(k^{\times}\{\Sigma'\})$, 
$(L^{\times})^{(\Sigma)}\defeq L^{\times}/(l^{\times}\{\Sigma'\})$, respectively, 
and a commutative diagram: 
$$
\CD
H'_K   @>{\rho}>> H'_L \\
@VVV              @VVV \\
\overline H_K  @>{\bar \rho}>> \overline H_L\\
\endCD
$$
where the vertical arrows are the natural surjective homomorphisms and the 
horizontal arrows are natural isomorphisms induced by $\sigma$. 
Here $\Sigma'\defeq \Primes \setminus \Sigma$, and
$k^{\times}\{\Sigma'\}$ (resp. $l^{\times}\{\Sigma'\}$)
is the $\Sigma'$-primary part of the multiplicative group
$k^{\times}$ (resp.  $l^{\times}$).
Using, among other facts, that the set $\Sigma$ 
is $k$-large and satisfies $(\epsilon_X)$, we show that the equalities $\overline H_K=K^{\times}/k^{\times}$,
$\overline H_L=L^{\times}/l^{\times}$, $H'_K=(K^{\times})^{(\Sigma)}  $, and
$H'_L=(L^{\times})^{(\Sigma)}$ hold.
Thus, one deduces naturally from the isomorphism 
$\sigma : G_{K}^{(\Sigma)}\isom G_{L}^{(\Sigma)}$ a commutative diagram:

$$
\CD
(K^{\times})^{(\Sigma)}   @>{\rho}>> (L^{\times})^{(\Sigma)} \\
@VVV              @VVV \\
K^{\times}/k^{\times}   @>{\bar \rho}>>  L^{\times}/l^{\times}\\
\endCD
$$
where the vertical arrows are the natural surjective homomorphisms and the 
horizontal arrows are the isomorphisms induced by $\sigma$.

{\it Step 4.}\ One shows that the isomorphism $\bar \rho:K^{\times}/k^{\times}
\isom L^{\times}/l^{\times}$ between principal divisors has the property that it preserves 
the valuation 
of the functions, or equivalently divisors, with respect to the  set-theoretic bijection
$\phi: X^{\cl}\isom Y^{\cl}$
established (in Step 1) between the sets of closed points of $X$, $Y$, respectively.

{\it Step 5.}\ We think of the elements of 
$(K^{\times})^{(\Sigma)}=K^{\times}/(k^{\times}\{\Sigma'\})$ and 
$(L^{\times})^{(\Sigma)}=L^{\times}/(l^{\times}\{\Sigma'\})$ 
as  ``{\it pseudo-functions}'', i.e., 
classes of rational functions modulo $\Sigma'$-primary constants. In particular, given a
pseudo-function $f'\in (K^{\times})^{(\Sigma)}$ 
(resp. $g'\in (L^{\times})^{(\Sigma)}$), and a closed point $x\in X^{\cl}$ 
(resp. $y\in Y^{\cl}$) it makes sense to consider the $\Sigma$-value $f'(x)$ (resp. $g'(y)$)
of $f'$ (resp. $g'$) (cf. discussion before Lemma 4.6). Then the isomorphism 
$\rho:(K^{\times})^{(\Sigma)}
\isom (L^{\times})^{(\Sigma)}$ has the property that it preserves the 
$\Sigma$-values of the pseudo-functions with respect to the  set-theoretic bijection
$\phi: X^{\cl}\isom Y^{\cl}$
established between the sets of closed points of $X$, $Y$, respectively.

{\it Step 6.}\ We think of the elements of $K^{\times}/k^{\times}$ (respectively, $L^{\times}/l^{\times}$)
as points of the infinite-dimensional projective spaces associated to the $k$ (resp. $l$)-vector 
spaces $K$ (resp. $L$). Using again, in an essential way, the fact that the set 
$\Sigma$ is $k$-large, and satisfies $(\epsilon_X)$, as well as the above property of
the isomorphism $\rho:(K^{\times})^{(\Sigma)}
\isom (L^{\times})^{(\Sigma)}$, we show that the isomorphism 
$\bar \rho:K^{\times}/k^{\times}\isom L^{\times}/l^{\times}$ in the above diagram, viewed as 
one between points of projective spaces, preserves colineations. Thus, by the fundamental theorem of 
projective geometry (cf. [Artin]), it arises from 
a uniquely determined $\psi _0$-isomorphism 
$$\psi: (K,+) \isom  (L,+),\  \psi(1)=1,$$
where  $\psi _0:k\isom l$ is a uniquely determined field isomorphism and 
$\psi$ is an isomorphism of abelian groups which is compatible with $\psi _0$. 
Finally, we show that the isomorphism  $\psi: (K,+) \isom  (L,+)$ preserves multiplication
so that it is a field isomorphism. By passing to open subgroups of 
$G_{K}^{(\Sigma)}$ and $G_{L}^{(\Sigma)}$ which correspond to each other via $\sigma$, 
one constructs a field isomorphism $K\sptilde \isom  L\sptilde$ 
which is compatible with $\psi$, and the inverse $L\sptilde \isom  K\sptilde$ 
of this isomorphism is the desired isomorphism. 

Note that the above idea
to resort to the fundamental theorem 
of projective geometry is not new in anabelian geometry 
(see, e.g., [Bogomolov], [Pop2]), 
while the above idea
to consider ``{\it pseudo-functions}'' and ``$k$-{\it largeness}'' 
is (to the best of our knowledge) new in anabelian geometry.

This paper is divided in two main parts. 
Part I is mostly of local nature. 
In Part I, \S 1, we review some basic facts on the Galois theory of function fields of algebraic curves 
and the main (well-known) results of the so-called local theory on the characterization
of the decomposition subgroups in Galois groups. In Part I, \S2, 
we reconstruct, using the local theory in \S1, various information encoded 
in the geometrically pro-$\Sigma$ absolute Galois group of a function field of a curve over a finite field. 
Part II is of global nature.
In Part II, \S 3, we define and give various characterizations of the notions of small and large sets of 
primes, and we also prove the main Proposition 3.11 which plays a crucial role in the proof of our main result.
Finally, in 
Part II, \S 4, we state and prove our main result Theorem 4.1.


\head
Part I
\endhead

In this first part we describe the local information encoded in the geometrically pro-$\Sigma$ 
absolute Galois group of the function field of a
curve over a finite field, 
and how much of this information is preserved under isomorphisms
between geometrically pro-$\Sigma$ absolute Galois groups. 

\subhead
\S 1. Generalities on Galois Groups of Function Fields of Curves
\endsubhead
In this section we fix some notations that we will use in this paper, and review some basic facts
on Galois groups of function fields of algebraic curves.
Let $k$ be a finite field of characteristic $p>0$. Let $X$ be a proper, 
smooth, and geometrically connected curve over $k$. 
Let $K$ be the function field of $X$. Let $\eta=\Spec K$ be the generic 
point of $X$ and $\bar \eta=\Spec \varOmega$ a geometric point of $X$ 
above $\eta$.  Write $K^{\sep}$ (resp. $\bar k=k^{\sep}$) for the separable 
closure of $K$ (resp. $k$) in $\varOmega$.  
Write $G=G_K\defeq \Gal (K^{\sep}/K)$ and 
$G_k\defeq \Gal (\bar k/k)$ for the absolute Galois groups of $K$ and $k$, 
respectively. We have the following exact sequence of profinite groups:
$$1\to \overline G\to G@>\text {pr}>> G_k\to 1,\tag{1.1}$$
where $\overline G$ is the absolute Galois group $\Gal (K^{\sep}/K\bar k)$ 
of 
$K\bar k$, and $\text{pr}:G\twoheadrightarrow G_k$ is the canonical projection.
It is well-known that  the kernel $\overline G$ of the projection 
$\text{pr}:G\twoheadrightarrow G_k$ is a free profinite group of 
countably infinite rank 
(cf. [Pop1] 
and [Harbater]). However, the structure of the extension (1.1) 
is not known.

We shall consider a variant of (1.1) above. Let $\Cal C$ be a full 
class of finite groups, i.e., $\Cal C$ is closed under taking subgroups, 
quotients, finite products, and extensions. For a profinite group $H$, denote
by $H^{\Cal C}$ the maximal pro-$\Cal C$ quotient of $H$. 
Given a profinite group $H$ and 
a closed normal subgroup $\overline H$ of $H$, 
we set $H^{(\Cal C)}\defeq 
H/\Ker (\overline H\twoheadrightarrow \overline H^{\Cal C})$. 
(Observe that $\Ker (\overline H\twoheadrightarrow \overline H^{\Cal C})$
is a normal subgroup of $H$ since it is a characteristic subgroup of $\overline H$.) 
Note that $H^{(\Cal C)}$ coincides with $H^{\Cal C}$ if
and only if the quotient $A\defeq H/\overline H$ is a pro-$\Cal C$ group. 
Let $\Primes$ denote the set of all prime numbers. 
When $\Cal C$ is the class of finite $\Sigma$-groups, where $\Sigma\subset \Primes$ is a 
set of prime numbers, write $H^{\Sigma}$ and $H^{(\Sigma)}$, instead of $H^{\Cal C}$ 
and $H^{(\Cal C)}$, respectively. 
(In later sections, the notation $H^{(\Sigma)}$ is used for a slightly more general 
setting where $H$ is 
a (not necessarily profinite) topological group and 
$\overline H$ is a closed normal subgroup of $H$ which is profinite.) 
By 
definition, we have the following commutative diagram:
$$
\CD
1@>>>    \overline H    @>>>     H   @>>>    A     @>>> 1\\
  @.        @VVV            @VVV              @V\id VV \\
1     @>>> \overline H^{\Cal C}  @>>>   H^{(\Cal C)}    @>>>A@>>>
1
\endCD
$$
where the rows are exact and the columns are surjective. 

\proclaim{Lemma 1.1}
Let $\ell$ be a prime number and $i\geq 0$. Assume either $\cd_{\ell}(N)\leq 1$ 
or $\Bbb F_{\ell}\not\in\Cal C$, where 
$N\defeq\Ker(\overline{H}\twoheadrightarrow\overline{H}^{\Cal C})$. 

\noindent
{\rm (i)} Let $M$ be a finite discrete $\ell$-primary $\overline{H}^{\Cal C}$-module. 
Then
$$
H^i(\overline{H}^{\Cal C}, M)=
\cases
H^i(\overline{H}, M), &\text{if $\Bbb F_{\ell}\in\Cal C$ or $i=0$},\\
0, &\text{if $\Bbb F_{\ell}\not\in\Cal C$ and $i>0$}. 
\endcases
$$

In particular, 
$$\cd_{\ell}(\overline{H}^{\Cal C})
\cases 
\leq\cd_{\ell}(\overline{H}), &\Bbb F_{\ell}\in\Cal C,\\
=0, &\Bbb F_l\not\in\Cal C. 
\endcases
$$

\noindent
{\rm (ii)} Let $M$ be a finite discrete $\ell$-primary $H^{(\Cal C)}$-module. 
Then
$$
H^i(H^{(\Cal C)}, M)=
\cases
H^i(H, M), &\text{if $\Bbb F_{\ell}\in\Cal C$},\\
H^i(A, M^{\overline{H}}), &\text{if $\Bbb F_{\ell}\not\in\Cal C$}. 
\endcases
$$

In particular, 
$$\cd_{\ell}(H^{(\Cal C)})
\cases
\leq\cd_{\ell}(H), &\Bbb F_{\ell}\in\Cal C,\\
=\cd_{\ell}(A), &\Bbb F_{\ell}\not\in\Cal C. 
\endcases$$
\endproclaim

\demo{Proof}
(i) If $\Bbb F_{\ell}\not\in\Cal C$, then $\overline{H}^{\Cal C}$ is of order prime to $\ell$, 
hence 
$$H^i(\overline{H}^{\Cal C},M)=
\cases
0, &i>0,\\
M^{\overline{H}^{\Cal C}}=M^{\overline{H}}=H^0(\overline{H},M), &i=0,
\endcases
$$
as desired. If $\Bbb F_{\ell}\in\Cal C$, then $\cd_{\ell}(N)\leq 1$ by assumption, hence 
$H^j(N,M)=0$ for $j>1$. Further, $H^1(N,M)=\Hom(N,M)=0$, as $\overline{H}^{\Cal C}$ 
is the maximal pro-$\Cal C$ quotient of $\overline{H}$. From these, we have 
$$H^i(\overline{H},M)=H^i(\overline{H}^{\Cal C}, H^0(N,M))=H^i(\overline{H}^{\Cal C}, M),$$
as desired (cf. [Neukirch-Schmidt-Winberg], (1.6.6)Proposition). 
The second assertion follows from the first. 

\noindent
(ii) If $\Bbb F_{\ell}\not\in\Cal C$, then, by (i), 
$$H^i(H^{(\Cal C)}, M)=H^i(A, H^0(\overline{H}^{\Cal C}, M))=H^i(A, M^{\overline{H}}),$$
as desired. If $\Bbb F_{\ell}\in\Cal C$, then, similarly to the proof of (i),  
$$H^i(H,M)=H^i(H^{(\Cal C)}, H^0(N,M))=H^i(H^{(\Cal C)}, M),$$
as desired. 
The second assertion follows from the first. (For the inequality 
$\cd_{\ell}(H^{(\Cal C)})\geq\cd_{\ell}(A)$, note that any finite discrete $\ell$-primary 
$A$-module $M$ can be regarded naturally as a finite discrete $\ell$-primary 
$H^{(\Cal C)}$-module with $M^{\overline{H}}=M$.) 
\qed\enddemo

Applying the above construction to 
(1.1) we obtain the exact sequence:
$$1\to \overline G^{\Cal C}\to G^{(\Cal C)}@>\pr>> G_k\to 1.
$$
We shall refer to the quotient
$G^{(\Cal C)}$ of $G$ as the maximal geometrically pro-$\Cal C$ quotient of $G$. 
Let $U$ be an open subgroup of $G^{(\Cal C)}$ and $H$ the inverse 
image of $U$ in $G$ via the canonical map $G\twoheadrightarrow G^{(\Cal C)}$. 
Then $H$ is an open subgroup of $G$ corresponding to a finite extension $K'/K$ of $K$. 
Further, $H$ (resp. $U$) is naturally identified with the absolute Galois group 
$\Gal (K^{\sep}/K')$ of $K'$ (resp. with $H^{(\Cal C)}=\Gal (K^{\sep}/K')^{(\Cal C)}$).

Let 
$$\Sigma \subset \Primes$$ 
be a set of prime numbers. 
Set
$$\Sigma ^{\dag}\defeq \Sigma \setminus \{p\}$$ 
and 
$$\Sigma'\defeq\Primes\setminus\Sigma.$$
Write 
$$\hat \Bbb Z^{\Sigma}\defeq \prod_{\ell\in \Sigma} \Bbb Z_{\ell}.$$ 

For a field $\kappa$ of characteristic $p>0$, with a 
separable closure $\kappa ^{\sep}$, we shall write
$$M_{\kappa^{\sep}}^{\Sigma}
\defeq \Hom 
(\Bbb Q/\Bbb Z,(\kappa^{\sep})^{\times})\otimes_{\hat \Bbb Z}
{\hat \Bbb Z^{\Sigma^{\dag}}}.
$$
Thus, $M_{\kappa ^{\sep}}^{\Sigma}$ is a free ${\hat \Bbb Z^{\Sigma^{\dag}}}$-module of rank 
one. Further, $M_{\kappa ^{\sep}}^{\Sigma}$ has a natural structure of $G_{\kappa}
\defeq \Gal (\kappa ^{\sep}/\kappa)$-module, which is isomorphic to the Tate 
twist ${\hat \Bbb Z^{\Sigma^{\dag}}}(1)$, i.e., $G_{\kappa }$ acts on $M_{{\kappa}
^{\sep}}^{\Sigma}$ via the $\Sigma$-part of the cyclotomic character 
$\chi _{\Sigma}:G_{\kappa}\to (\hat \Bbb Z^{\Sigma^{\dag}})^{\times}$.
In particular, we write
$$M_{\bar k}^{\Sigma}= 
\Hom (\Bbb Q/\Bbb Z,
\bar k^{\times})\otimes_{\hat \Bbb Z}
{\hat \Bbb Z^{\Sigma^{\dag}}}.
$$
Similarly, we shall write
$$M_X^{\Sigma}\defeq M_{K^{\sep}}^{\Sigma}=
\Hom (\Bbb Q/\Bbb Z,(K^{\sep})^{\times})\otimes_{\hat \Bbb Z}
{\hat \Bbb Z^{\Sigma^{\dag}}}.
$$
Note that $M_X^{\Sigma}$ has a natural structure of $G$-module, which is  naturally
identified with the $G$-module $M_{\bar k}^{\Sigma}$ (the natural $G$- and 
$G_k$-module structures of $M_{\bar k}^{\Sigma}$ are compatible with respect to the 
natural projection $\pr:G\twoheadrightarrow G_k$).

For a scheme $T$ denote by $T^{\cl}$ the set of closed points of $T$. Let $K\sptilde /K$
be the subextension of $K^{\sep}/K$ corresponding to the subgroup 
$\Ker (G\twoheadrightarrow G^{(\Sigma)})$ of $G$,
and let $\tilde X$ be the normalization of $X$ in $K\sptilde $.
The Galois group $G^{(\Sigma)}$ acts naturally on the set ${\tilde X}^{\cl}$ and the quotient  
of ${\tilde X}^{\cl}$ by this action is naturally identified with $X^{\cl}$. 
For a point $\tilde x\in {\tilde X}^{\cl}$, with residue field 
$k(\tilde x)$ (which is an algebraic closure of the residue field $k(x)$ of $x$), we define its decomposition 
group $D_{\tilde x}$ and inertia group  $I_{\tilde x}$ by 
$$D_{\tilde x}\defeq \{\gamma \in G^{(\Sigma)}\  |\  \gamma (\tilde x)=\tilde x\}$$
and
$$I_{\tilde x}\defeq \{\gamma \in D_{\tilde x}\  |\  \gamma \ \text {acts trivially on}\  
k (\tilde x)\},$$
respectively. We have a canonical exact sequence:
$$1\to I_{\tilde x}\to D_{\tilde x}\to G_{k(x)}\defeq \Gal (k (\tilde x)/
k(x))\to 1.$$

For a profinite group $H$ we write $\Sub (H)$ for the set of closed subgroups of $H$.
Let $\Sigma \subset \Primes$ be a set of prime numbers. 
The following are well-known facts concerning the decomposition 
and inertia subgroups of the geometrically pro-$\Sigma$ Galois group $G^{(\Sigma)}$.

\proclaim 
{Proposition 1.2}\ (Properties of Decomposition and Inertia Subgroups) 
Let $\tilde x\in \tilde X^{\cl}$, and $x$ the image in 
$X^{\cl}$ of $\tilde x\in \tilde X^{\cl}$. 

\noindent 
{\rm (i)} Let $X'$ be the normalization of $X$ in $K^{\sep}$ and 
$x'$ a point of $(X')^{\cl}$ above $\tilde x$. Let 
$I_{x'}\subset D_{x'}\subset G$ be the inertia and the decomposition subgroups of $G$ at $x'$. 
(Thus, $D_{x'}\defeq \{\gamma \in G\  |\  \gamma (x')=x'\}$, and 
$I_{x'}\defeq \{\gamma \in D_{x'}\  |\  \gamma \ \text {acts trivially on}\  
k (x')\}$.) Then we have 
$D_{\tilde x}=D_{x'}^{(\Sigma)}$. More precisely, we have the 
following commutative diagram: 
$$
\matrix
1&\to &I_{x'}&\to &D_{x'}&\to &G_{k(x)}&\to &1\\
&&&&&&&&\\
&& \downarrow && \downarrow && \Vert &&\\
&&&&&&&&\\
1&\to &I_{x'}^{\Sigma}&\to &D_{x'}^{(\Sigma)}&\to &G_{k(x)}&\to &1\\
&&&&&&&&\\
&& \Vert && \Vert && \Vert &&\\
&&&&&&&&\\
1&\to &I_{\tilde x}&\to &D_{\tilde x}&\to &G_{k(x)}&\to &1
\endmatrix
$$
where the horizontal rows are exact and the vertical arrows are surjective. 

\noindent
{\rm (ii)} 
The inertia subgroup $I_{\tilde x}$ possesses a unique $p$-Sylow 
subgroup $I_{\tilde x}^{\text w}$. The quotient
$I_{\tilde x}^{\t}\defeq I_{\tilde x}/I_{\tilde x}^{\text w}$ is isomorphic to $\hat \Bbb 
Z^{\Sigma^{\dag}}$, and is naturally identified with the Galois group 
$\Gal (K^{\t}_x/K^{\ur}_x)$, where $K^{\ur}_x$ (resp. $K^{\t}_x$) is
the maximal unramified (resp. tamely ramified) extension of the $x$-adic completion $K_x$ of $K$. 
We have a natural exact sequence:
$$1\to I_{\tilde x}^{\t}\to D_{\tilde x}^{\t}\to G_{k(x)}\to 1,$$
where $D_{\tilde x}^{\t}\defeq \Gal (K_x^{\t}/K_x)$.

In particular, 
$I_{\tilde x}^{\t}$
has a natural structure of $G_{k(x)}$-module. Further, there exists a natural 
identification $I_{\tilde x}^{\t}\isom M_{k(\tilde x)}^{\Sigma}$ of $G_{k(x)}$-modules.
\endproclaim

\demo
{Proof} 
(i) The only nontrivial point in the assertion is that the natural 
homomorphism $I_{x'}^{\Sigma}\to
\overline{G}^{\Sigma}$ (whose image coincides with $I_{\tilde x}$) 
is injective. A proof of this fact is as follows. 

{\it Step 1.} First, if $p\not\in\Sigma$, this follows from 
[Sa\"\i di-Tamagawa3], Lemma 1.3. 

{\it Step 2.} Next, consider the special case $\Sigma=\{p\}$. Then we have 
$$\cd_p(I_{\tilde x})\leq\cd_p(\overline{G}^{\{p\}})\leq \cd_p(\overline{G})\leq 1.$$
Indeed, the first inequality follows from the fact that 
$I_{\tilde x}\subset \overline{G}^{\{p\}}$ ([Serre1]), 
the second inequality follows from Lemma 1.1 (i) (and the fact 
that $\Ker(\overline{G}\twoheadrightarrow \overline{G}^{\{p\}})$ is the 
absolute Galois group of a field of characteristic $p$ ([Serre1])), 
and the third inequality 
follows from the fact that $\overline{G}$ is the absolute Galois group of 
a field of characteristic $p$ ([Serre1]). 
In particular, the surjective 
homomorphism 
$I_{x'}^{\{p\}}\twoheadrightarrow I_{\tilde x}$ 
of pro-$p$ groups 
admits a section $s: I_{\tilde x} \to I_{x'}^{\{p\}}$. 
Now, the 
homomorphism 
$I_{x'}^{\{p\}}\rightarrow \overline{G}^{\{p\}}$ 
induces a homomorphism 
$(I_{x'}^{\{p\}})^{\ab}/p\to 
(\overline{G}^{\{p\}})^{\ab}/p$ of pro-$p$ abelian groups killed by $p$. 
By Artin-Schreier theory, the (Pontryagin) dual of this 
last homomorphism is identified with 
$K\bar k/\wp(K\bar k)\to (K\bar k)_{\bar x}/\wp((K\bar k)_{\bar x})$. 
Here, $\bar x$ denotes the image in $(X\times_k\bar k)^{\cl}$ 
of $\tilde x\in \tilde X^{\cl}$, 
$(K\bar k)_{\bar x}$ denotes the $\bar x$-adic completion of 
$K\bar k$, and $\wp:\alpha\mapsto\alpha^p-\alpha$. 
Observe that $\wp((K\bar k)_{\bar x})\supset \wp (\Cal O_{\bar x})= \Cal O_{\bar x}$ 
by Hensel's lemma, where $\Cal O_{\bar x}$ denotes the ring of integers of $(K\bar k)_{\bar x}$. 
Now, since the natural homomorphism $K\bar k\to (K\bar k)_{\bar x}/\Cal O_{\bar x}$ 
is surjective (as follows from the Riemann-Roch theorem), the homomorphism 
$K\bar k/\wp(K\bar k)\to (K\bar k)_{\bar x}/\wp((K\bar k)_{\bar x})$ is also surjective. 
Thus, we have $(I_{x'}^{\{p\}})^{\ab}/p \hookrightarrow (\overline{G}^{\{p\}})^{\ab}/p$, 
hence, a fortiori, 
$(I_{x'}^{\{p\}})^{\ab}/p \hookrightarrow (I_{\tilde x})^{\ab}/p$. 
Thus, $(I_{x'}^{\{p\}})^{\ab}/p \isom (I_{\tilde x})^{\ab}/p$. 
In particular, $s(I_{\tilde x})$ must surject onto $(I_{x'}^{\{p\}})^{\ab}/p$. 
By the Frattini property, this implies that $s(I_{\tilde x})=I_{x'}^{\{p\}}$, 
hence $I_{x'}^{\{p\}}\isom I_{\tilde x}$, as desired. 

{\it Step 3.} 
Finally, consider a general $\Sigma$. By Step 1, we may assume that $p \in \Sigma$. For a profinite 
group $H$, denote by $H^{\Sigmadagbyp}$ the (unique) maximal quotient of $H$ 
which is an extension of a pro-$\Sigma^{\dag}$ group by a (normal) pro-$p$ group. 
Let 
$H(p)$ denote the kernel of $H^{\Sigmadagbyp}\twoheadrightarrow H^{\Sigma^{\dag}}$. 
In general, the natural homomorphism $H^{\Sigma}\twoheadrightarrow H^{\Sigmadagbyp}$
is not an isomorphism. But we have $I_{x'}^{\Sigma}\isom 
I_{x'}^{\Sigmadagbyp}$, since the wild inertia subgroup (i.e. the $p$-Sylow 
subgroup of the inertia group) is normal. 
Thus, we have the following commutative diagram
$$
\matrix
1&\to& I_{x'}(p)&\to& I_{x'}^{\Sigma}&\to&I_{x'}^{\Sigma^{\dag}}&\to&1\\
&&&&&&&&\\
&&\downarrow&&\downarrow&&\downarrow&&\\
&&&&&&&&\\
1&\to& \overline{G}(p)&\to& \overline{G}^{\Sigmadagbyp}&\to&\overline{G}^{\Sigma^{\dag}}&\to&1
\endmatrix
$$
where the horizontal rows are exact. Here, the right vertical arrow is injective by Step 1, 
as $p\not\in\Sigma^{\dag}$. On the other hand, the left vertical arrow can be obtained as the projective 
limit of homomorphisms $(I_{x'}\cap \overline{H})^{\{p\}}\to 
\overline{H}^{\{p\}}$, where $\overline{H}$ runs over all open subgroups of $\overline{G}$ 
that contain the kernel of $\overline{G}\twoheadrightarrow \overline{G}^{\Sigma^{\dag}}$. 
Thus, the left vertical arrow is injective by Step 2, hence the middle vertical 
arrow $I_{x'}^{\Sigma}\to \overline{G}^{\Sigmadagbyp}$ 
is also injective. Now, the homomorphism 
$I_{x'}^{\Sigma}\to \overline{G}^{\Sigma}$ is, a fortiori, injective, as desired. 

\noindent 
(ii) This follows from (i), together with 
well-known facts on ramification theory (cf. [Serre2], Chapitre IV).
\qed
\enddemo

In fact, the decomposition subgroups of $G^{(\Sigma)}$ are completely determined 
by their group-theoretic structure. More precisely, we have the 
following fundamental result.

\proclaim
{Proposition 1.3} (Galois Characterization of Decomposition Subgroups) 
Consider the natural map 
$D=D_{\tilde X}^{(\Sigma)}: {\tilde X}^{\cl}\to \Sub (G^{(\Sigma)}),\ \tilde x \mapsto D_{\tilde x}$. 

\noindent
{\rm (i)} The map $D$ is Galois-equivariant. More precisely, 
for $g\in G^{(\Sigma)}$ and $\tilde x\in \tilde X^{\cl}$, 
we have $D_{g\tilde x}=gD_{\tilde x}g^{-1}$. 

\noindent
{\rm (ii)} Assume $\Sigma\neq\emptyset$. Then $D$ is injective. More precisely, 
let $\tilde x\neq \tilde x_1$ be two elements of ${\tilde X}^{\cl}$, 
then $D_{\tilde x}\cap D_{\tilde x_1}$ is pro-$\Sigma'$ and is of infinite index 
both in $D_{\tilde x}$ and in $D_{\tilde x_1}$. 

\noindent
{\rm (iii)} Assume $\Sigma^{\dag}\neq\emptyset$ and $\ell\in \Sigma^{\dag}$. 
Let $\Dec_{\ell}(G^{(\Sigma)})\subset\Sub(G^{(\Sigma)})$ be the set of closed subgroups 
$\frak D$ of $G^{(\Sigma)}$ satisfying the following property: 
There exists an open subgroup $\frak D_0$ of $\frak D$ such that for 
any open subgroup $\frak D'\subset\frak D_0$, 
$\dim_{\Bbb F_{\ell}}H^2(\frak D',\Bbb F_{\ell})=1$. Define $\Dec^{\max}_{\ell}(G^{(\Sigma)})\subset
\Dec_{\ell}(G^{(\Sigma)})$ 
to be the set of maximal elements of $\Dec_{\ell}(G^{(\Sigma)})$ with respect to 
the inclusion relation. 
Then the image of $D$ coincides with $\Dec^{\max}_{\ell}(G^{(\Sigma)})$. 
(In particular, $\Dec^{\max}_{\ell}(G^{(\Sigma)})$ does not depend on the choice of 
$\ell\in\Sigma^{\dag}$.) 

Thus, $D: \tilde X^{\cl}\to \Sub(G^{(\Sigma)})$ induces a natural, Galois-equivariant 
bijection 
$\tilde X^{\cl}\isom \Dec^{\max}_{\ell}(G^{(\Sigma)})$. 
\endproclaim

\demo
{Proof} 
(i) This follows from the definition of decomposition group. 

\noindent
(ii) Let $\ell\in\Sigma$. If $\ell\neq p$, then 
$D_{\tilde x}\cap D_{\tilde x_1}$ is of order prime to $\ell$ 
by [Sa\"\i di-Tamagawa3], Proposition 1.5 (i). The case $\ell=p$ can be treated 
along the same lines, by resorting to Artin-Schreier theory instead of Kummer theory. 
More precisely, 
let $D_p$ be a $p$-Sylow subgroup of $D_{\tilde x}\cap D_{\tilde x_1}$, 
and suppose that $D_p\neq 1$. 
We have 
$$\cd_p(D_p)\leq \cd_p(G^{(\Sigma)})\leq \cd_p(G)\leq 1<\infty.$$
Indeed, the first inequality follows from the fact that 
$D_p\subset G^{(\Sigma)}$ ([Serre1]), 
the second inequality follows from Lemma 1.1 (ii) (and the fact 
that $\Ker(\overline{G}\twoheadrightarrow \overline{G}^{\Sigma})=\Ker (G\twoheadrightarrow G^{(\Sigma)})$ is the 
absolute Galois group of a field of characteristic $p$ ([Serre1])), 
and the third inequality 
follows from the fact that $G$ is the absolute Galois group of 
a field of characteristic $p$ ([Serre1]). In particular, 
$D_p$ is torsion-free, hence is infinite. 
Thus, one may replace 
$G^{(\Sigma)}$ by any open subgroup, and assume that 
the images $x, x_1$ in $X^{\cl}$ of $\tilde x,\tilde x_1\in\tilde X^{\cl}$ 
are distinct, and that 
the image of $D_p$ in $(G^{(\Sigma)})^{\ab}/
p$ is nontrivial. 
In particular, this implies that the natural map 
$$D_{\tilde x}^{\ab}/
p \times 
D_{\tilde x_1}^{\ab}/
p \to 
(G^{(\Sigma)})^{\ab}/
p$$ 
induced by the group operation of $(G^{(\Sigma)})^{\ab}/p$ 
is not injective. By Artin-Schreier theory and Proposition 1.2 (i), 
this last condition is equivalent 
to saying that the natural map 
$$K/\wp(K)\to K_x/\wp(K_x)\times K_{x_1}/\wp(K_{x_1})$$ 
is not surjective, where $\wp: \alpha\mapsto \alpha^p-\alpha$. 
(Observe that by Hensel's lemma 
$\wp(K_x)$ (resp. $\wp(K_{x_1})$) 
contains the maximal ideal of the ring of integers of $K_x$ 
(resp. $K_{x_1}$), hence is 
open in $K_x$ (resp. $K_{x_1}$).) 
This contradicts the approximation theorem (cf. 
[Neukirch], Lemma 8). 

Thus, 
$D_{\tilde x}\cap D_{\tilde x_1}$ is pro-$\Sigma'$, which also implies that 
it is of infinite index both in $D_{\tilde x}$ and in $D_{\tilde x_1}$. 
(Note that for $\ell\in\Sigma$ the pro-$\ell$-Sylow subgroups of $D_{\tilde x}$ 
and $D_{\tilde x_1}$ are infinite.) In particular, $D$ must be injective. 

\noindent
(iii) This is a special case of [Sa\"\i di-Tamagawa3], Proposition 1.5 (ii). 
(This result goes back to [Uchida], where the case 
$\Sigma=\Primes$ is treated.) 
\qed
\enddemo

\definition{Remark 1.4} 
For other characterizations of decomposition groups, 
see [Sa\"\i di-Tamagawa3], Remark 1.6. 
\enddefinition

\subhead 
\S 2. Isomorphisms between Geometrically Pro-$\Sigma$ Galois Groups
\endsubhead
In this section we follow the notations in $\S 1$. 
Let $k$, $l$ be finite fields of characteristic 
$p_k$, $p_l$, respectively, and of cardinality $q_k$, $q_l$, respectively. 

Let $X$, $Y$ be smooth, proper, and geometrically connected 
curves over $k$, $l$, respectively. 
Let $K$, $L$ be the function fields of $X$, $Y$, respectively. 
We will write $G_{K}\defeq \Gal (K^{\sep}/K)$, 
$G_{L}\defeq \Gal (L^{\sep}/L)$ 
for the absolute Galois groups of $K$, $L$, respectively. 

Let $\Sigma_X, \Sigma_Y\subset \Primes$
be sets of prime numbers. 
We assume that 
the $\Sigma_X$-cyclotomic character $\chi _{\Sigma_X}:G_{k}\to 
\prod _{\ell\in \Sigma_X\setminus\{p_k\}}\Bbb Z_{\ell}^{\times}$ is injective. 
(In the terminology of \S 3 (cf. Definition/Proposition 3.1), this is equivalent to 
saying that $\Sigma_X$ is not $k$-small.) 
Write $G_K^{(\Sigma_X)}$ (resp. $G_L^{(\Sigma_Y)}$)
for the  maximal geometrically pro-$\Sigma_X$ (resp. $\Sigma_Y$) quotient of $G_K$ (resp. $G_L$). 
Thus, we have exact sequences:
$$1\to {\overline G}_K^{\Sigma _X}\to G_K^{(\Sigma _X)}@>\pr>> G_{k}\to 1,$$
and 
$$1\to {\overline G}_L^{\Sigma _Y}\to G_L^{(\Sigma _Y)}@>\pr>> G_{l}\to 1,$$
where 
$G_{k}\defeq \Gal (\bar k/k)$ 
(resp. $G_{l}\defeq \Gal (\bar l/l)$) 
is the absolute Galois group of $k$ (resp. $l$), and
${\overline G}_K^{\Sigma_X}$ 
(resp. ${\overline G}_L^{\Sigma_Y}$) 
is the maximal pro-$\Sigma _X$ 
(resp. pro-$\Sigma_Y$) quotient
of the absolute Galois group ${\overline G}_K\defeq\Gal (K^{\sep}/K\bar k)$ 
(resp. ${\overline G}_L\defeq\Gal (L^{\sep}/L\bar l)$) 
of 
$K\bar k$ (resp. 
$L\bar l$). 

For the rest of this section we will consider an isomorphism of profinite groups
$$\sigma : G_K^{(\Sigma _X)}\isom G_L^{(\Sigma _Y)}$$ 
between the maximal geometrically pro-$\Sigma _X$ (resp. pro-$\Sigma_Y$) quotient of the absolute Galois group 
$G_K$ (resp. $G_L$). We write $\tilde X$ (resp. $\tilde Y$) 
for the normalization of $X$ (resp. $Y$) in $K\sptilde$ (resp. $L\sptilde$). 
Here, 
$K\sptilde/K$ (resp. $L\sptilde/L$)
is the subextension of $K^{\sep}/K$ (resp. $L^{\sep}/L$) 
with Galois group  $G_K^{(\Sigma _X)}$ (resp. $G_L^{(\Sigma _Y)}$).

Recall $\Sigma_X^{\dag}=\Sigma _X\setminus \{p_k\}$, and  $\Sigma_Y^{\dag}=\Sigma _Y\setminus \{p_l\}$.

\proclaim
{Lemma 2.1} (Invariance of Sets of Primes)

\noindent
{\rm (i)} We have $\Sigma _X^{\dag}=\Sigma_Y^{\dag}$,
$\Sigma_X\cap\{p_k\}=\Sigma_Y\cap\{p_l\}$, and 
$\Sigma_X=\Sigma_Y$. 
Set $\Sigma^{\dag}\defeq\Sigma_X^{\dag}=\Sigma_Y^{\dag}$, 
$\Sigma\defeq\Sigma_X=\Sigma_Y$, and 
$\Sigma'\defeq \Primes \setminus \Sigma$. 

\noindent
{\rm (ii)} $\Sigma$ is infinite. 

\endproclaim

\demo
{Proof} 
(i) It follows from global class field theory for $K$ that 
(for a prime number $\ell$) the maximal pro-$\ell$ quotient $(G_K^{(\Sigma_X)})^{\ab,\ell}$ of 
the maximal abelian quotient  $(G_K^{(\Sigma_X)})^{\ab}$ of $G_K^{(\Sigma_X)}$is described as follows: 
$$
(G_K^{(\Sigma_X)})^{\ab,\ell}
\simeq 
\cases
\Bbb Z_{\ell}, &\ell\not\in\Sigma_X, \\
\Bbb Z_{\ell}\times \overline{(\text{infinite torsion group})}, &\ell\in\Sigma_X^{\dag}, \\
\Bbb Z_{\ell}^{\aleph_0}, &\ell\in \Sigma_X\cap\{p_k\}. 
\endcases
$$
Here $\overline{(\text{infinite torsion group})}$ denotes the closure of the torsion 
subgroup (which is infinite) 
of 
$(G_K^{(\Sigma_X)})^{\ab,\ell}$.
A similar description holds for  $(G_L^{(\Sigma_Y)})^{\ab,\ell}$.
This implies $\Sigma_X^{\dag}=\Sigma_Y^{\dag}$, 
$\Sigma_X\cap\{p_k\}=\Sigma_Y\cap\{p_l\}$, and 
$\Sigma_X=\Sigma_Y$. 

\noindent
(ii) This follows immediately from the assumption that 
the $\Sigma_X$-cyclotomic character $\hat\Bbb Z\simeq G_{k}\to 
\prod _{\ell\in \Sigma_X\setminus\{p_k\}}\Bbb Z_{\ell}^{\times}$ is injective. 
\qed\enddemo

\proclaim
{Lemma 2.2} (Set-Theoretic  Correspondence between Points) The isomorphism $\sigma$ induces naturally a 
bijection:
$$\tilde \phi :{\tilde X}^{\cl}\isom {\tilde Y}^{\cl},\ \tilde x\mapsto \tilde y,$$
such that 
$$\sigma (D_{\tilde x})=D_{\tilde y},\ \forall \tilde x\in {\tilde X}^{\cl},$$ 
where $D_{\tilde x}$ (resp. $D_{\tilde y}$) is the decomposition subgroup of  
$G_K^{(\Sigma)}$ (resp. $G_L^{(\Sigma)}$) at $\tilde x$ (resp. $\tilde y$),  and 
$\tilde\phi$ induces a bijection 
$$\phi :X^{\cl}\isom Y^{\cl},\ x\mapsto y,$$
where $x$ (resp. $y$) is the image of $\tilde x$ (resp. $\tilde y$) in $X^{\cl}$ (resp. $Y^{\cl}$). Thus, in particular, 
$\sigma$ induces naturally a 
bijection:
$$\phi : \Div_{X}\isom \Div_{Y}$$
between the groups of divisors of $X$ and $Y$, respectively.
\endproclaim

\demo
{Proof} This follows from Proposition 1.3 and Lemma 2.1 (to ensure 
$\Sigma^{\dag}=\Sigma_X^{\dag}=\Sigma_Y^{\dag}\neq\emptyset$). 
\qed
\enddemo

Let $x\in X^{\cl}$, and $y\defeq \phi (x)\in Y^{\cl}$. 
Write $K_x$ (resp. $L_y$) for the completion of $K$ (resp. $L$) at $x$ (resp. $y$). Denote the ring of integers  
of $K_x$ (resp. $L_y$) by $\Cal O_x$ (resp. $\Cal O_y$). Write
$D_x\defeq D_{\tilde x}$ (resp. $D_y\defeq D_{\tilde y}$) and $I_x\defeq I_{\tilde x}$
(resp.  $I_y\defeq I_{\tilde y}$)
for the decomposition and the inertia subgroups of $G_K^{(\Sigma)}$ (resp. $G_L^{(\Sigma)}$)
at $\tilde x$ (resp. $\tilde y$), where $\tilde x\in \tilde X^{\cl}$ (resp. $\tilde y\in \tilde Y^{\cl}$)
is a point above $x$ (resp. $y$).
Thus, $D_x$ (resp. $D_y$) is defined only up to conjugation.  
By Proposition 1.2 (i) and local class field theory (cf., e.g., [Serre3]), 
we have natural isomorphisms 
$$(K_{x}^{\times})^{\wedge,(\Sigma)}\isom D_x^{\ab},$$
and
$$(L_{y}^{\times})^{\wedge,(\Sigma)}\isom D_y^{\ab},$$
where $(K_{x}^{\times})^{\wedge}$ (resp. $(L_{y}^{\times})^{\wedge}$)
is the profinite completion of the 
topological group $K_{x}^{\times}$ (resp. $L_{y}^{\times}$),  
and we set 
$$(K_{x}^\times)^{\wedge,(\Sigma)}\defeq 
(K_{x}^\times)^{\wedge}/\Ker (\Cal O_{x}^{\times}\twoheadrightarrow (\Cal O_{x}^{\times})^{\Sigma})$$
and 
$$(L_{y}^\times)^{\wedge,(\Sigma)}\defeq 
(L_{y}^\times)^{\wedge}/\Ker (\Cal O_{y}^{\times}\twoheadrightarrow (\Cal O_{y}^{\times})^{\Sigma}),$$ 
where $(\Cal O_{x}^{\times})^{\Sigma}$, $(\Cal O_{y}^{\times})^{\Sigma}$ stand for the 
maximal pro-$\Sigma$ quotients of the profinite groups 
$\Cal O_{x}^{\times}$, $\Cal O_{y}^{\times}$, respectively. 

More concretely, we have 
$$(K_{x}^{\times})^{\wedge,(\Sigma)}=(K_{x}^{\times})^{\wedge}/N_x,\  
(\Cal O_{x}^{\times})^{\Sigma}=\Cal O_{x}^{\times}/N_x$$
with 
$$N_x\defeq\Ker (\Cal O_{x}^{\times}\twoheadrightarrow (\Cal O_{x}^{\times})^{\Sigma})=
\cases
U_{x}^1 (\Cal O_{x}^{\times,\tor}\{\Sigma'\}), &\Sigma=\Sigma^{\dag},\\
\Cal O_{x}^{\times,\tor}\{\Sigma'\}, &\Sigma\neq\Sigma^{\dag}, 
\endcases
$$
and we have a similar description for $(L_{y}^{\times})^{\wedge,(\Sigma)}$.
Here, 
$U_{x}^1$ is the group of principal units in $\Cal O_{x}^{\times}$, 
and $\Cal O_{x}^{\times,\tor}\{\Sigma'\}$ is 
the group of $\Sigma'$-primary torsion of $\Cal O_{x}^{\times}$. 
(Observe that $\Cal O_{x}^{\times,\tor}\{\Sigma'\}\isom k(x)^{\times}\{\Sigma'\}.$)

We have the following commutative diagram 
$$
\matrix
0&\to& \Cal O_{x}^{\times} &\to& K_{x}^{\times} &\to& \Bbb Z &\to&0\\
&&&&&&&&\\
&& \downarrow && \downarrow && \cap &&\\
&&&&&&&&\\
0&\to& (\Cal O_{x}^{\times})^{\Sigma} &\to& (K_{x}^\times)^{\wedge,(\Sigma)} &\to& \hat\Bbb Z &\to&0\\
&&&&&&&&\\
&& \phantom{\wr}\downarrow\wr && \phantom{\wr}\downarrow\wr && \phantom{\wr}\downarrow\wr &&\\
&&&&&&&&\\
0&\to& \Im(I_x) &\to& D_x^{\ab} &\to& G_{k} &\to&0, 
\endmatrix
$$
where the horizontal rows are exact. Here, the map $K_{x}^{\times} \to \Bbb Z$ 
is the $x$-adic valuation, $\Im (I_x)$ is the image of $I_x$ in $D_x^{\ab}$, and the map $\hat\Bbb Z\isom G_{k}$ sends $1\in\hat\Bbb Z$ 
to the $q_k$-th power Frobenius element in $G_{k}$. 

Further, the natural filtration 
$$(k(x)^{\times})^{\Sigma}, (U_{x}^1)^{\Sigma}
\subset(\Cal O_{x}^{\times})^{\Sigma}
\subset (K_{x}^{\times})^{(\Sigma)}
\subset (K_{x}^{\times})^{\wedge,(\Sigma)},$$
where $(U_{x}^1)^{\Sigma}$ denotes the maximal pro-$\Sigma$ quotient of
$U_{x}^1$ and
$(K_{x}^{\times})^{(\Sigma)}$ is the image of $K_x^{\times}$ in $(K_{x}^{\times})^{\wedge,(\Sigma)}$,
induces, 
via the above isomorphism $(K_{x}^{\times})^{\wedge,(\Sigma)}\isom D_x^{\ab}$, 
a filtration 
$$\Im((k(x)^{\times})^{\Sigma}), \Im((U_{x}^1)^{\Sigma})
\subset\Im((\Cal O_{x}^{\times})^{\Sigma})
\subset \Im((K_{x}^{\times})^{(\Sigma)})
\subset\Im((K_{x}^{\times})^{\wedge,(\Sigma)})=D_x^{\ab}.$$
Here, $\Im((\Cal O_{x}^{\times})^{\Sigma})$ coincides with the image 
$\Im(I_x)$ in $D_x^{\ab}$ of $I_x$. Similar statements and filtrations 
hold for  $(L_{y}^{\times})^{\wedge,(\Sigma)}$ and $D_y^{\ab}$.

Let 
$$\sigma _{x,y}:D_x\isom D_y$$ 
be the isomorphism of profinite groups induced by $\sigma$ (which is only defined up to conjugation)
(cf. Lemma 2.2). Write 
$$\sigma _{x,y}^{\ab}:D_x^{\ab}\isom D_y^{\ab}$$ 
for the induced isomorphism between the 
maximal abelian quotients of $D_x$ and $D_y$, respectively.

\proclaim
{Lemma 2.3} (Invariants of Isomorphisms between Geometrically Pro-$\Sigma$
Decomposition Groups) 

\noindent
{\rm (i)} The isomorphism $\sigma _{x,y}^{\ab}:D_x^{\ab}\to D_y^{\ab}$ 
preserves the images $\Im ((k(x)^{\times})^{\Sigma})$ and $\Im ((k(y)^{\times})^{\Sigma})$ , hence 
it induces naturally an isomorphism 
$$\tau _{x,y}:(k(x)^{\times})^{\Sigma} \isom (k(y)^{\times})^{\Sigma}$$ 
between the maximal pro-$\Sigma$ quotients of the multiplicative 
groups of the residue fields at $x$ and $y$, respectively. 

\noindent
{\rm (ii)}  The isomorphism $\sigma_{x,y}$ induces naturally an isomorphism $M_{\overline{k(x)}}^{\Sigma}\isom
M_{\overline{k(y)}}^{\Sigma}$, which is Galois-equivariant with respect to $\sigma_{x,y}$. 
In particular, $\sigma_{x,y}$ commutes with the $\Sigma$-parts of the cyclotomic 
characters $\chi _x:D_x\to ({\hat \Bbb Z}^{\Sigma^{\dag}})^{\times}$ (resp. 
$\chi _y:D_y\to ({\hat \Bbb Z}^{\Sigma^{\dag}})^{\times}$) 
of $D_x$ (resp. $D_y$),  
i.e., we have a commutative diagram:

$$
\CD
({\hat \Bbb Z}^{\Sigma^{\dag}})^{\times}  @= ({\hat \Bbb Z}^{\Sigma^{\dag}})^{\times}\\
@A{\chi _x}AA    @A{\chi _y}AA   \\
D_x  @>\sigma_{x,y}>> D_y \\
\endCD 
$$

\noindent
{\rm (iii)}
The isomorphism $\sigma_{x,y}$ preserves $I_x$ and $I_y$. 
\endproclaim

\demo
{Proof} The proofs of (i)(ii)(iii) are similar to those 
of [Sa\"\i di-Tamagawa3], Proposition 2.1 (iii)(iv)(v), 
respectively. More precisely: 

\noindent
(i) By Proposition 1.2 (i) and local class field theory, 
$\Im((k(x)^{\times})^{\Sigma})\subset D_x^{\ab}$ coincides with 
the torsion subgroup $D_{x}^{\ab,\tor}$ of $D_{x}^{\ab}$, and 
a similar statement holds for 
$\Im((k(y)^{\times})^{\Sigma})\subset D_y^{\ab}$. From this, 
the assertion follows. 

\noindent
(ii) By applying (i) to open subgroups of $D_x$, $D_y$, (which 
correspond to each other via $\sigma_{x,y}$), and 
passing to the projective limit, 
we obtain a natural isomorphism
$M_{\overline{k(x)}}^{\Sigma}\isom M_{\overline{k(y)}}^{\Sigma}$ 
between the modules of roots of unity. More precisely, let 
$E$ be a finite extension of $K_{x}$ corresponding to 
an open subgroup $H$ of $D_x$ and $h$ the residue field of $E$. 
Then the following diagram commutes:
$$
\matrix
(h^{\times})^{\Sigma}&\subset& (E^{\times})^{\wedge,(\Sigma)} &\isom&  H^{\ab} \\
\downarrow &&\downarrow && \downarrow\\
(k(x)^{\times})^{\Sigma}&\subset&(K_{x}^{\times})^{\wedge,(\Sigma)} &\isom&  D_x^{\ab}, 
\endmatrix
$$
where the map $H^{\ab}\to D_x^{\ab}$ is induced
by the natural inclusion $H\subset D_x$ 
and the map $(E^{\times})^{\wedge,(\Sigma)}\to (K_{x}^{\times})^{\wedge,(\Sigma)}$ 
is induced by the norm map $E^{\times} \to K_{x}^{\times}$. 
The map $(h^{\times})^{\Sigma}\to (k(x)^{\times})^{\Sigma}$ is induced by 
the (norm) map $(E^{\times})^{\wedge,(\Sigma)}\to (K_{x}^{\times})^{\wedge,(\Sigma)}$, 
hence coincides with the $e$-th power of the map 
$(h^{\times})^{\Sigma}\to (k(x)^{\times})^{\Sigma}$ induced by 
the norm map $h^{\times}\to k(x)^{\times}$, where $e$ denotes 
the ramification index of $E/K_{x}$. 
Thus, if we consider 
the projective subsystem formed by the open subgroups of $D_x$ that are 
obtained as the inverse image of an open subgroup of $D_x^{\ab}/(\text{torsion})$, 
we get a projective system $((h^{\times})^{\Sigma})$ with surjective 
transition homomorphisms (as all the ramification indices are powers of $p$) 
whose limit is identified with $M_{\overline{k(x)}}^{\Sigma}$. Indeed, 
(again as all the ramification indices are powers of $p$) 
the limit is unchanged if the projective system is replaced with the subsystem 
indexed by the open subgroups $H\subseteq D_x$ 
obtained as inverse image of open subgroups of
$D_x/I_x=G_{k(x)}$. 
Then the above norm map $h^{\times}\to k(x)^{\times}$ is just 
the $a$-th power map, where $a=\sum_{i=0}^{[h:k(x)]-1}|k(x)|^i=|h^{\times}|/|k(x)^{\times}|$. 
[This sort of precise argument involving suitable projective subsystems should have 
been inserted also in the proof of [Sa\"\i di-Tamagawa3], Proposition 2.1 (iv).] 
Further, this identification is (by construction) Galois-compatible with respect to the 
isomorphism $\sigma_{x,y}$, as desired. 
The second assertion follows from this Galois-compatibility. 

\noindent
(iii) The character 
$\chi _x:D_x\to ({\hat \Bbb Z}^{\Sigma^{\dag}})^{\times}$ 
(resp. 
$\chi _y:D_y\to ({\hat \Bbb Z}^{\Sigma^{\dag}})^{\times}$) 
factors as 
$D_x\twoheadrightarrow D_x/I_x=G_{k(x)}\overset{\chi_{k(x)}}\to{\to}({\hat \Bbb Z}^{\Sigma^{\dag}})^{\times}$ 
(resp. 
$D_y\twoheadrightarrow D_y/I_y=G_{k(y)}\overset{\chi_{k(y)}}\to{\to}({\hat \Bbb Z}^{\Sigma^{\dag}})^{\times}$), 
where $\chi_{k(x)}$ (resp. $\chi_{k(y)}$) is the $\Sigma$-cyclotomic character of 
$G_{k(x)}$ (resp. $G_{k(y)}$). Further, since the $\Sigma$-cyclotomic character 
of $G_k$ is assumed to be injective, $\chi_{k(x)}$ is also injective. 
Thus, $I_x$ coincides with the kernel of $\chi_x$ and 
$I_y$ is included in the kernel of $\chi_y$. 

Now, it follows from (ii) that $\sigma_{x,y}(I_x)\supset I_y$, hence 
$$\hat\Bbb Z\simeq G_{k(x)}=D_x/I_x\overset{\sigma_{x,y}}\to{\isom}D_y/\sigma_{x,y}(I_x)\twoheadleftarrow D_y/I_y
=G_{k(y)}\simeq\hat\Bbb Z.$$
As any surjective homomorphism $\hat\Bbb Z\to\hat\Bbb Z$ is an isomorphism, this shows 
$\sigma_{x,y}(I_x)=I_y$, as desired. 
\qed
\enddemo

\noindent
{\bf 2.4.} {\it Invariants of Isomorphisms between Geometrically Pro-$\Sigma$ Galois Groups.}

\proclaim
{Lemma 2.4.1}The following diagram is commutative:

$$
\CD 
(\hat \Bbb Z^{\Sigma^{\dag}})^{\times}    @=   (\hat \Bbb Z^{\Sigma^{\dag}})^{\times} \\
@A\chi_{k}AA                                  @A\chi_{l}AA \\
G_{k}  @.   G_{l} \\
@A\pr_KAA                    @A\pr_LAA   \\
G_K^{(\Sigma)}   @>\sigma>>  G_L^{(\Sigma)} \\
\endCD 
\tag{2.1}
$$
where $\chi_{k}$ (resp. $\chi _l$) is the $\Sigma$-part of the cyclotomic character of $G_{k}$ (resp. $G_l$). 
\endproclaim

\demo{Proof} For each $\tilde x\in 
\tilde X^{\cl}$, with $\tilde y\defeq \tilde \phi (\tilde x)\in \tilde Y^{\cl}$, we 
have the following diagram:
$$
\CD 
(\hat \Bbb Z^{\Sigma^{\dag}})^{\times}    @=   (\hat \Bbb Z^{\Sigma^{\dag}})^{\times} \\
@A\chi_KAA                    @A\chi_LAA   \\
G_K^{(\Sigma)}   @>\sigma>>  G_L^{(\Sigma)}  \\
@AAA                    @AAA   \\
D_{\tilde x} @>{\sigma _{\tilde x,\tilde y}}>>   D_{\tilde y} \\
\endCD\tag{2.2}
$$
where the maps $D_{\tilde x}\to G_K^{(\Sigma)}$ 
and $D_{\tilde y}\to G_L^{(\Sigma)}$ 
are the natural inclusions,
the lower square is commutative, and $\chi_K$ (resp. $\chi_L$) 
is the $\Sigma$-part of the cyclotomic 
character of $G_K^{(\Sigma)}$ (resp. $G_L^{(\Sigma)}$).
Since the restriction of  $\chi_K$ (resp. $\chi_L$) to $D_{\tilde x}$ 
(resp. $D_{\tilde y}$) 
coincides with the $\Sigma$-part of the cyclotomic character 
of $D_{\tilde x}$ (resp. $D_{\tilde y}$), 
the exterior square of (2.2) is commutative 
by Lemma 2.3  (ii). 
Hence the upper square of (2.2) is also commutative, 
since  $G_K^{(\Sigma)}$ is (topologically) generated by the 
decomposition subgroups $D_{\tilde x}$, $\forall \tilde x\in {\tilde X}^{\cl}$, 
as follows from Chebotarev's density theorem. 
The commutativity of the diagram (2.1) follows from this, since 
$\chi_{k}\circ \pr_K=\chi_K$ and 
$\chi_{l}\circ \pr_L=\chi_L$. 
\qed
\enddemo

\proclaim {Lemma 2.4.2}
The isomorphism $\sigma$ commutes with the canonical 
surjections 
$G_K^{(\Sigma)}\twoheadrightarrow \pi_1(X)^{(\Sigma)}$
(resp. 
$G_L^{(\Sigma)}\twoheadrightarrow \pi_1(Y)^{(\Sigma)}$),
where $\pi_1(X)^{(\Sigma)}\defeq\pi_1(X)/
\Ker(\pi_1(\overline{X})\twoheadrightarrow\pi_1(\overline{X})^{\Sigma})$
(resp. $\pi_1(Y)^{(\Sigma)}\defeq\pi_1(Y)/
\Ker(\pi_1(\overline{Y})\twoheadrightarrow\pi_1(\overline{Y})^{\Sigma})$)
is the maximal geometrically pro-$\Sigma$ quotient of the fundamental group 
$\pi_1(X)$ (resp. $\pi_1(Y)$) of $X$ (resp. $Y$) (with respect to the geometric point $\Spec(K^{\sep})\to X$ (resp.  $\Spec(L^{\sep})\to Y$)). 
More precisely, we have a commutative 
diagram: 
$$
\CD
G_K ^{(\Sigma)}   @>>>    \pi_1(X)^{(\Sigma)}    \\
@V\sigma VV              @VVV \\
G_L ^{(\Sigma)} @>>>   \pi_1(Y)^{(\Sigma)}
\endCD 
$$
where the vertical arrows are isomorphisms. 
\endproclaim

\demo{Proof} Let $\Cal I_X$ (resp. $\Cal I_Y$) denote the closed normal subgroup of $G_K^{(\Sigma)}$
(resp. $G_L^{(\Sigma)}$) (topologically) generated by the inertia subgroups. Then the isomorphism $\sigma $ 
maps $\Cal I_X$ onto $\Cal I_Y$ by Lemma 2.2 and Lemma 2.3 (iii). Since 
$\pi_1(X)^{(\Sigma)}=G_K^{(\Sigma)}/\Cal I_X$ and $\pi_1(Y)^{(\Sigma)}=G_L^{(\Sigma)}/\Cal I_Y$, 
the assertion follows.
\qed
\enddemo

\proclaim {Lemma 2.4.3}
The isomorphism $\sigma$ commutes with the canonical 
projections $\pr_K:G_K^{(\Sigma)}\twoheadrightarrow G_{k}$ and
$\pr_L:G_L^{(\Sigma)}\twoheadrightarrow G_{l}$,
i.e., we have a commutative 
diagram: 
$$
\CD
G_K ^{(\Sigma)}   @>\pr_K>>    G_{k}    \\
@V\sigma VV              @VVV \\
G_L ^{(\Sigma)} @>\pr _L>>   G_{l} 
\endCD 
\tag{2.3}$$
where the vertical arrows are isomorphisms. 
\endproclaim

\demo{Proof}
This follows from Lemma 2.4.2, since we have $G_k=\pi_1(X)^{(\Sigma),\ab}/(\text {torsion})$
as a quotient of $G_K^{(\Sigma)}$ (cf. [Tamagawa], Proposition 3.3 (ii)). 
Similarly, $G_l=\pi_1(Y)^{(\Sigma),\ab}/(\text {torsion})$
as a quotient of $G_L^{(\Sigma)}$. Here, $\pi_1(X)^{(\Sigma),\ab}$ (resp. $\pi_1(Y)^{(\Sigma),\ab}$)
is the maximal abelian quotient of $\pi_1(X)^{(\Sigma)}$ (resp. $\pi_1(Y)^{(\Sigma)}$).
(Alternatively, this follows from Lemma 2.4.1. Indeed, since $\chi_k\circ \pr_K=\chi_K$, 
$\chi_l\circ \pr_L=\chi_L$, and $\pr_K$, $\pr_L$ are surjective, $\Im (\chi_K)\subset (\hat\Bbb Z^{\Sigma^{\dag}})^{\times}$
coincides with $\Im(\chi_k)$ (similarly, $\Im (\chi_L)\subset (\hat\Bbb Z^{\Sigma^{\dag}})^{\times}$
coincides with $\Im(\chi_l)$).
By assumption 
$\chi_{k}$ is injective. Since this injectivity condition is equivalent to 
requiring $\Im(\chi_{k})\simeq \hat\Bbb Z$ (as abstract profinite groups), 
we see that both $\chi_{k}$ and $\chi_l$ are injective. 
In summary, we have 
$$G_{k}\isom\Im(\chi_{k})=\Im(\chi_K)\subset(\hat\Bbb Z^{\Sigma^{\dag}})^{\times}$$ 
and 
$$G_{l}\isom\Im(\chi_{l})=\Im(\chi_L)\subset(\hat\Bbb Z^{\Sigma^{\dag}})^{\times}.$$ 
Now, the assertion follows from the commutativity of the diagram (2.1).)  
\qed
\enddemo

\proclaim {Lemma 2.4.4}
For each subset $T\subset \Sigma$, 
the isomorphism $\sigma$ commutes with the canonical 
surjections $G_K^{(\Sigma)}\twoheadrightarrow G_K^{(T)}$,
and  $G_L^{(\Sigma)}\twoheadrightarrow G_L^{(T)}$,
i.e., we have a commutative 
diagram: 
$$
\CD
G_K ^{(\Sigma)}   @>>>    G_K^{(T)}    \\
@V\sigma VV              @VVV \\
G_L ^{(\Sigma)} @>>>   G_L^{(T)} 
\endCD
$$
where the vertical arrows are isomorphisms. 
\endproclaim

\demo{Proof} This follows from Lemma 2.4.3, since the quotient $G_K^{(\Sigma)}\twoheadrightarrow G_K^{(T)}$ 
can be characterized as 
$$G_K^{(T)}=G_K^{(\Sigma)}/\Ker(\Ker(\pr_K)\twoheadrightarrow(\Ker(\pr_K))^T),$$ 
and 
a similar statement holds for $G_L^{(\Sigma)}\twoheadrightarrow G_L^{(T)}$. 
\qed
\enddemo

\proclaim {Lemma 2.4.5}
The bijection $\phi:X^{\cl}\isom Y^{\cl}$ commutes 
with the degree functions $\deg_{X}: X^{\cl}\to \Bbb Z_{>0}$, 
$x\mapsto [k(x):k]$, and $\deg_{Y}: Y^{\cl}\to \Bbb Z_{>0}$, 
$y\mapsto [k(y):l]$. 
\endproclaim

\demo{Proof}This follows from Lemmas 2.2 and 2.4.3. 
Indeed, for each $x\in X^{\cl}$, take $\tilde x\in \tilde X^{\cl}$ 
above $x$ and set $y=\phi(x)$ and $\tilde y=\tilde\phi(\tilde x)$ 
(which is above $y$). 
Then we have 
$$\deg_{X}(x)=(G_{k}:\pr_K(D_{\tilde x}))=(G_{l}:\pr_L(D_{\tilde y}))=
\deg_Y(y).\ \qed$$ 
\enddemo

\proclaim {Lemma 2.4.6} For each integer $n>0$, let $k\subset k_n\subset \bar k$ 
(resp. $l\subset l_n\subset \bar l$) 
denote the unique extension with $[k_n:k]=[l_n:l]=n$. Then we have 
$\sharp(X(k_n))=\sharp(Y(l_n))$ for all $n>0$. 
\endproclaim

\demo{Proof}This follows from Lemma 2.4.5, since 
$$\sharp(X(k_{n}))=\sum_{0< d\mid n}d\cdot\sharp(\deg_{X}^{-1}(d)),$$
and
$$\sharp(Y(l_{n}))=\sum_{0< d\mid n}d\cdot\sharp(\deg_{Y}^{-1}(d)).\  \qed$$
\enddemo

\proclaim {Lemma 2.4.7}
{\rm (i)} We have $q_k=q_l$. In particular, $p_k=p_l$. 

\noindent
{\rm (ii)} Notations as in Lemma 2.4.6, we have 
$\sharp(k_n)=\sharp(l_n)$ for all $n>0$. 
\endproclaim

\demo{Proof} (i) 
This follows from Lemma 2.4.6 
(cf. [Pop2], Lemma 2.3). 
More precisely, by the Weil estimate, we have 
$$1+q_k^{n}-2g_Xq_k^{\frac{n}{2}}\leq \sharp(X(k_{n}))\leq 1+q_k^{n}+2g_Xq_k^{\frac{n}{2}},$$
where $g_X$ denotes the genus of $X$, hence 
$$\frac{\sharp(X(k_{n}))}{q_k^n}\to 1\ (n\to\infty),$$
and similarly $\frac{\sharp(Y(l_{n}))}{q_l^n}\to 1\ (n\to\infty)$.
Now, by Lemma 2.4.6, we obtain 
$$\left(\frac{q_k}{q_l}\right)^n=\frac{q_k^n}{q_l^n}\to 1\ (n\to \infty),$$
which implies $q_k=q_l$, as desired. 

\noindent
(ii) This follows immediately from (i), as 
$\sharp(k_n)=q_k^n$ and $\sharp(l_n)=q_l^n$. 
\qed
\enddemo

Set $p\defeq p_k=p_l$ and $q\defeq q_k=q_l$.

\proclaim {Lemma 2.4.8}
The bijection $\phi:X^{\cl}\isom Y^{\cl}$ commutes 
with the norm functions $N_X: X^{\cl}\to \Bbb Z_{>0}$, 
$x\mapsto \sharp(k(x))$, and $N_Y: Y^{\cl}\to \Bbb Z_{>0}$, 
$y\mapsto \sharp(k(y))$. 
\endproclaim

\demo{Proof}This follows from Lemmas 2.4.5 and 2.4.7 (i). 
\qed
\enddemo

\proclaim {Lemma 2.4.9}
The isomorphism
$G_{k}\to G_{l}$ induced naturally by $\sigma$ 
(cf. Lemma 2.4.3) 
maps the $q$-th 
power Frobenius element $\varphi _{k}$ of $G_{k}$ to the $q$-th 
power Frobenius element $\varphi _{l}$ of $G_{l}$.
\endproclaim

\demo{Proof} As shown in the (alternative) proof of Lemma 2.4.3, 
$\chi_{k}:G_{k}\to(\hat\Bbb Z^{\Sigma^{\dag}})^{\times}$ 
and 
$\chi_{l}:G_{l}\to(\hat\Bbb Z^{\Sigma^{\dag}})^{\times}$ are 
injective. Thus, $\varphi_{k}\in G_{k}$ can be characterized 
by the property $\chi_{k}(\varphi_{k})=q\in \hat \Bbb Z^{\Sigma^{\dag}}$, 
and similarly $\chi_l(\varphi_{l})=q\in \hat \Bbb Z^{\Sigma^{\dag}}$. 
Now, the assertion follows from the commutativity of the diagram
(2.1) in Lemma 2.4.1, and the diagram (2.3) in Lemma 2.4.3. 
\qed
\enddemo

\proclaim{Lemma 2.4.10}
In the notation of Lemma 2.4.6, we have 
$\sharp(J_{X}(k_{n}))=\sharp(J_{Y}(l_{n}))$ 
for all $n>0$, where $J_{X}$ (resp. $J_Y$) denotes the jacobian variety of $X$ (resp. $Y$). 
\endproclaim

\demo{Proof} By Lemmas 2.4.2 and 2.4.3, 
the isomorphism $\pi_1(X)^{(\Sigma)}\isom\pi_1(Y)^{(\Sigma)}$ induced by 
$\sigma$ preserves $\pi_1(\overline{X})^{\Sigma}=\Ker(\pi_1(X)^{(\Sigma)}
\twoheadrightarrow G_{k})$ and $\pi_1(\overline{Y})^{\Sigma}=\Ker(\pi_1(Y)^{(\Sigma)}
\twoheadrightarrow G_{l})$. It follows from this that $\sigma$ induces 
an isomorphism
$$T(J_{X})^{\Sigma}=(\pi_1(\overline{X})^{\Sigma})^{\ab}\isom
(\pi_1(\overline{Y})^{\Sigma})^{\ab}=T(J_{Y})^{\Sigma}$$ 
(where $T(J_{X})\defeq \prod_{\ell\in\Primes}T_{\ell}(J_{X})$ 
(resp. $T(J_{Y})\defeq \prod_{\ell\in\Primes}T_{\ell}(J_{Y})$)
is the full Tate module of $J_{X}$ (resp. $J_Y$) and $T(J_{X})^{\Sigma}$ (resp. $T(J_{Y})^{\Sigma}$)
is its maximal pro-$\Sigma$ quotient), which is Galois-equivariant with respect to 
the isomorphism $G_{k}\isom G_{l}$ in Lemma 2.4.3. Thus, it follows from Lemma 2.4.9 that 
$P_{X,n}=P_{Y,n}$ for all $n>0$, where 
$P_{X,n}$ (resp. $P_{Y,n}$) denotes 
the characteristic polynomial for the action of $\varphi_{k}^n$ (resp. $\varphi_{l}^n$)
on the free $\hat\Bbb Z^{\Sigma^{\dag}}$-module 
$T(J_{X})^{\Sigma^{\dag}}$ (resp. $T(J_{Y})^{\Sigma^{\dag}}$). 
Now, we have 
$$\sharp(J_{X}(k_{n}))=P_{X,n}(1)=P_{Y,n}(1)=\sharp(J_{Y}(l_{n})),$$
as desired. 
\qed
\enddemo

\definition {Remark 2.4.11} 
Let $\ell$ be a prime $\neq p=p_k$. When $\Sigma=\{\ell\}$, most of the results presented in 
Lemmas 2.4.1-2.4.9 are proved in [Pop2], Part I, 2, without resorting to Lemma 2.3 which relies 
heavily on local class field theory. (In fact, in [Pop2], function fields with arbitrary transcendence 
degree are also treated.)

Further, when $\Sigma=\Sigma^{\dag}$ (i.e., $\Sigma\not\ni p$), the quotient 
$G_K^{(\Sigma)}\twoheadrightarrow G_k$ can be identified with 
$G_K^{(\Sigma)}\twoheadrightarrow 
(G_K^{(\Sigma)})^{\ab}/
\overline{
(G_K^{(\Sigma)})^{\ab,\tor}
}$, 
where $(G_K^{(\Sigma)})^{\ab,\tor}$ is the torsion subgroup of 
$(G_K^{(\Sigma)})^{\ab}$ and $\overline{(G_K^{(\Sigma)})^{\ab,\tor}}$ 
is its closure in $(G_K^{(\Sigma)})^{\ab}$
(cf. the proof of Lemma 2.1 (i)). It follows from this that for each $\ell\in\Sigma$, 
the quotient $G_K^{(\Sigma)}\twoheadrightarrow G_K^{(\{\ell\})}$ can be recovered 
group-theoretically from $G_K^{(\Sigma)}$. Thus, 
most of the results presented in Lemmas 2.4.1-2.4.9 for this case could be reduced to the case 
$\Sigma=\{\ell\}$ basically. 

However, in the general case where $\Sigma$ may contain $p$, 
the authors do not know any quick way (without establishing Lemma 2.3 first) 
of reconstructing the quotient $G_K^{(\Sigma)}\twoheadrightarrow G_k$ 
and reducing to the case $\Sigma=\{\ell\}$. 
\enddefinition

\proclaim{Lemma 2.5}
(Invariance of Filtrations of Geometrically Pro-$\Sigma$ Decomposition Groups) 
Let the notations be as in Lemma 2.3 and the discussion before Lemma 2.3. 
Then the isomorphism $\sigma _{x,y}^{\ab}:D_x^{\ab}\to D_y^{\ab}$ 
preserves the filtrations 
$$\Im((k(x)^{\times})^{\Sigma}), \Im((U_{x}^1)^{\Sigma}) 
\subset\Im((\Cal O_{x}^{\times})^{\Sigma})
\subset \Im((K_{x}^{\times})^{(\Sigma)})
\subset\Im((K_{x}^{\times})^{\wedge,(\Sigma)})=D_x^{\ab},$$ 
and
$$\Im((k(y)^{\times})^{\Sigma}), \Im((U_{y}^1)^{\Sigma}) 
\subset\Im((\Cal O_{y}^{\times})^{\Sigma})
\subset \Im((L_{y}^{\times})^{(\Sigma)})
\subset\Im((L_{y}^{\times})^{\wedge,(\Sigma)})=D_y^{\ab}.$$
\endproclaim

\demo{Proof}
First, $\sigma^{\ab}_{x,y}$ preserves $\Im((k(x)^{\times})^{\Sigma})$ 
and $\Im((k(y)^{\times})^{\Sigma})$ 
by Lemma 2.3 (i). Next, 
$\sigma^{\ab}_{x,y}$ preserves 
$\Im((\Cal O_{x}^{\times})^{\Sigma})=\Im(I_x)$ 
and 
$\Im((\Cal O_{y}^{\times})^{\Sigma})=\Im(I_y)$ 
by Lemma 2.3 (iii), 
and preserves $\Im((U_{x}^1)^{\Sigma})$ and $\Im((U_{y}^1)^{\Sigma})$
by Lemma 2.4.7 (i), since 
$\Im((U_{x}^1)^{\Sigma})$ (resp. $\Im((U_{y}^1)^{\Sigma})$) 
is the pro-$p$ Sylow group of 
$\Im((\Cal O_{x}^{\times})^{\Sigma})$ (resp. $\Im((\Cal O_{y}^{\times})^{\Sigma})$). 
Finally, $\sigma^{\ab}_{x,y}$ preserves 
$\Im((K_{x}^{\times})^{(\Sigma)})=\pr_X^{-1}(\varphi_{k}^{\Bbb Z})$
and
$\Im((L_{y}^{\times})^{(\Sigma)})=\pr_Y^{-1}(\varphi_{l}^{\Bbb Z})$
by Lemma 2.4.3 and lemma 2.4.9.  
\qed
\enddemo

\head
Part II 
\endhead

In this part we introduce the notion of ``small'' and ``large'' sets of primes, and 
we state and prove our main results. 

\subhead 
\S 3. Small and Large Sets of Primes
\endsubhead
Let $\Primes$ be the set of all prime numbers and 
$\Sigma \subset \Primes$ a 
subset. Set
$\Sigma '\defeq \Primes \setminus \Sigma$. Let $k$ be a finite field of characteristic $p>0$ 
and set $\Sigma^{\dag}=\Sigma\setminus\{p\}$. 
Write 
$$\hat \Bbb Z^{\Sigma}\defeq \prod_{\ell\in \Sigma} \Bbb Z_{\ell}.$$ 
For a prime number $\ell\in \Primes \setminus \{p\}$ let 
$$\chi_{\ell}:G_k\to\Bbb Z_{\ell} ^{\times}$$ 
be the $\ell$-adic cyclotomic character of $k$,
and define the $\Sigma$-part of the cyclotomic character of $k$ by: 
$$\chi_\Sigma\defeq (\chi_{\ell})_{\ell\in\Sigma^{\dag}}: G_k\to 
(\hat \Bbb Z^{\Sigma^{\dag}})^{\times}=\prod_{\ell\in\Sigma^{\dag}}\Bbb Z_{\ell}^{\times}.$$ 
Thus, we have 
$$\bar k^{\Ker(\chi_\Sigma)}=
k_{\Sigma}\defeq k(\zeta _{\ell^j}\ \vert \ \ell\in \Sigma^{\dag}, j\in 
\Bbb Z_{\ge 0}).$$
For a prime number $\ell\in \Primes$,
let $G_{k,\ell}\subset G_k$ be the pro-$\ell$-Sylow subgroup of $G_k$. 
(Recall that $G_k\simeq\hat \Bbb Z$ and $G_{k,\ell}\simeq\Bbb Z_{\ell}$.) 

\definition 
{Definition/Proposition 3.1} (Small Set of Primes) {\sl Let $\Sigma \subset \Primes$ 
be a set of prime numbers. We say that the set $\Sigma$ is $k$-small if the following 
equivalent conditions are satisfied:

\noindent
{\rm (i)} $k_{\Sigma} \neq \bar k$.

\noindent
{\rm (ii)} The $\Sigma$-part $\chi_{\Sigma}$ of the cyclotomic character is not injective.

\noindent
{\rm (iii)} There exists a prime number $\ell_0\in \Primes$, such that $\sharp
(\chi _{\Sigma}(G_{k,\ell_0}))< \infty$.

\noindent
{\rm (iii$'$)} There exists a prime number $\ell_0\in \Primes$, such that 
$\ell_0\not\in \Sigma^{\dag}$ and that there exists $N_0\in\Bbb Z_{\geq 0}$ 
satisfying that for any $\ell\in\Sigma^{\dag}$, 
the order of $p\mod \ell\in \Bbb F_{\ell}^{\times}$ is not divisible by $\ell_0^{N_0}$. 

\noindent
{\rm (iv)} There exists a subfield $k\subset k' \subset \bar k$
such that $(G_k:G_{k'})=\infty$ and that 
$(\chi_\Sigma(G_k): \chi_\Sigma(G_{k'}))<\infty$}.
\enddefinition

\demo
{Proof} Easy.
\qed
\enddemo

\definition 
{Definition 3.2} (Large Set of Primes) Let $\Sigma \subset \Primes$ 
be a 
set of prime numbers. 
We say that the set $\Sigma$ is {\it $k$-large} if the set $\Sigma '=\Primes \setminus \Sigma$ is $k$-small.
\enddefinition

\proclaim{Proposition 3.3} 
Let $\Sigma \subset \Primes$ be a set of prime numbers. Consider the following 
conditions: 

\noindent
{\rm (i)} $\Sigma$ is cofinite, i.e., $\Sigma'$ is finite. 

\noindent
{\rm (ii)} $\Sigma$ is $k$-large. 

\noindent
{\rm (ii$'$)} $\Sigma$ is not $k$-small. 

\noindent
{\rm (i$'$)} $\Sigma$ is infinite. 

\noindent
Then we have the following implications: 
$\text{\rm (i)}
\implies\text{\rm (ii)}
\implies\text{\rm (ii$'$)}
\implies\text{\rm (i$'$)}$. 
\endproclaim

\demo{Proof}
To prove the implication $\text{\rm (ii)}\implies\text{\rm (ii$'$)}$, 
suppose that $\Sigma$ is $k$-large and $k$-small at a time, or, equivalently, that 
both $\Sigma$ and $\Sigma'$ are $k$-small. This contradicts 
[Grunewald-Segal], Theorem A, as $\Primes=\Sigma\cup\Sigma'$. 
To prove the implication $\text{\rm (ii$'$)}\implies\text{\rm (i$'$)}$, 
suppose that $\Sigma$ is finite, then 
there is no injective homomorphism $\hat\Bbb Z\to (\hat\Bbb Z^{\Sigma^{\dag}})^{\times}$. 
In particular, $\chi_{\Sigma}:G_{k}\to (\hat\Bbb Z^{\Sigma^{\dag}})^{\times}$ is not injective, 
i.e., $\Sigma$ is $k$-small, which is a contradiction. 
The implication $\text{\rm (i)}\implies\text{\rm (ii)}$ is obtained by applying 
the implication $\text{\rm (ii$'$)}\implies\text{\rm (i$'$)}$ to $\Sigma'$. 
\qed
\enddemo

\definition {Remarks 3.4}
\enddefinition
\definition {3.4.1} Consider the following 
conditions: 

\noindent
{\rm (i)} $\Sigma$ is cofinite.

\noindent
{\rm (ii)} $\Sigma$ is $k$-large. 

\noindent
{\rm (iii)} $\Sigma$ is $\Bbb F_p$-large

\noindent
{\rm (iv)} $\Sigma$ is of (natural) density $1$. 

\noindent
{\rm (iv$'$)} $\Sigma$ is of (natural) density $\neq 0$. 

\noindent
{\rm (iii$'$)} $\Sigma$ is not $\Bbb F_p$-small. 

\noindent
{\rm (ii$'$)} $\Sigma$ is not $k$-small. 

\noindent
{\rm (i$'$)} $\Sigma$ is infinite. 

\noindent
Then we have the following implications: 
$$\matrix
\text{\rm (ii)}&\iff & \text{\rm (iii)} 
&\implies &\text{\rm (iii$'$)} &\iff &\text{\rm (ii$'$)} \\
&&&&&&\\
\Uparrow&&&&&&\Downarrow\\
&&&&&&\\
\text{\rm (i)} &\implies&\text{\rm (iv)} &\implies &\text{\rm (iv$'$)} &\implies& \phantom{,}\text{\rm (i$'$)}, 
\endmatrix$$
$\text{\rm (iii)}\implies\text{\rm (iv$'$)}$,  and $\text{\rm (iv)}\implies\text{\rm (iii$'$)}$.

Indeed, 
the implications (i)$\implies$(ii), 
(ii$'$)$\implies$(i$'$), and 
(iii)$\implies$(iii$'$) are proved in Proposition 3.3.
The implications (ii)$\iff$(iii) and 
(i)$\implies$(vi)$\implies$(vi$'$)$\implies$(i$'$)  
are immediate. To prove the implication 
(iii)$\implies$(iv$'$), suppose that $\Sigma$ is of density $0$. Then by [Grunewald-Segal], Theorem A, 
$\chi_{\Sigma'}:G_{\Bbb F_p}\to (\hat\Bbb Z^{(\Sigma')^{\dag}})^{\times}$ is injective. 
This is equivalent to saying that $\Sigma'$ is not $\Bbb F_p$-small, or 
that $\Sigma$ is not $\Bbb F_p$-large. The implications 
(iv)$\implies$(iii$'$)$\iff$(ii$'$) are obtained by applying the 
implications (ii)$\iff$(iii)$\implies$(iv$'$) to $\Sigma'$. 
\enddefinition

\definition {3.4.2} The implication (ii)$\implies$(i) in 3.4.1 does not always hold. 
To construct such an example, set $k=\Bbb F_p$ (for simplicity), 
consider a prime number $r\neq p$, $r\nmid p-1$, and define $\Sigma$ to be 
the set of prime numbers which do not divide $p^{r^m}-1$ for any $m\geq 0$. 
Then $\Sigma'$ is infinite. Indeed, we have 
$$(p-1) \mid(p^r-1)\mid\cdots\mid (p^{r^{m}}-1)\mid (p^{r^{m+1}}-1)\mid\cdots,$$
$$\frac{p^{r^{m+1}}-1}{p^{r^m}-1}
>1,$$
and 
$$\left(p^{r^m}-1,\frac{p^{r^{m+1}}-1}{p^{r^m}-1}\right)=(p^{r^m}-1,r)=(p-1,r)=1$$
by the Euclidean algorithm. 
(Here, to prove the second equality, use the fact that $p^{r^m}
\equiv p \pmod{r}$.) 
Thus, for each $m\geq 0$, there exists a prime number $\ell_{m+1}$ 
such that $\ell_{m+1}\mid (p^{r^{m+1}}-1)$ and that $\ell_{m+1}\nmid (p^{r^m}-1)$. 
This implies that $\Sigma'\ni \ell_1,\ell_2,\dots,\ell_m,\dots$ is infinite. 
On the other hand, $\Sigma'$ is $k$-small. Indeed, take any $\ell_0\in 
\Sigma\setminus\{r\}
$ 
(e.g., $\ell_0=p$). Then we claim that
$\chi_{\Sigma'}(G_{k,\ell_0})=\{1\}$. To prove this, it suffices to show that 
$\chi_s(G_{k,\ell_0})=\{1\}$ for all prime number $s\in (\Sigma')^{\dag}$. 
Further, since $\Ker((\Bbb Z_s)^{\times}\to(\Bbb F_s)^{\times})$ is pro-$s$ 
and $\ell_0\neq s$, it suffices to show that 
$\chi_s(G_{k,\ell_0})\bmod{s}=\{1\}$. 
But $\chi_s(G_{k,\ell_0})\bmod{s}$ is the $\ell_0$-Sylow subgroup of 
$\chi_s(G_k)\bmod{s}=\langle p\bmod{s}\rangle\subset (\Bbb F_s)^{\times}$. 
Since $s\mid (p^{r^m}-1)$ for some $m\geq 0$, the order of 
$\langle p\bmod{s}\rangle\subset (\Bbb F_s)^{\times}$ is a power of $r$. 
Now, as $\ell_0\neq r$, $\chi_s(G_{k,\ell_0})\bmod{s}=\{1\}$, as desired. 

Taking $\Sigma'$ in this example as $\Sigma$, we also see that 
the implication (i$'$)$\implies$(ii$'$) in 3.4.1 does not always hold. 
\enddefinition

\definition {3.4.3}
The implication (iii)$\implies$(iv) in 3.4.1 does not always hold. 
In fact, for any $\varepsilon>0$, there exists a set of prime numbers $\Sigma$ 
of density $<\varepsilon$, such that $\Sigma$ is $\Bbb F_p$-large. 
Indeed, for each $N\in\Bbb Z_{>0}$, set 
$\Sigma_i(N)\defeq\{\ell\in\Primes\mid \ell\equiv i\pmod{N}\}$, $i=1,\dots,N$. 
Then the density of $\Sigma_i(N)$ is $1/\varphi(N)$ (resp. $0$) 
if $(i,N)=1$ (resp. $(i,N)\neq 1$). Now, choose a prime number 
$N$ such that $\varphi(N)=N-1>1/\varepsilon$, and set 
$\Sigma\defeq \Sigma_1(N)\cup\Sigma_0(N)=\Sigma_1(N)\cup\{N\}$, 
whose density is $1/\varphi(N)<\varepsilon$. We claim that 
$\Sigma$ is $\Bbb F_p$-large. To see this, we have to prove that 
$\chi_{\Sigma'}: G_{\Bbb F_p}\to (\hat\Bbb Z^{(\Sigma')^{\dag}})^{\times}$ 
is not injective. But this follows from the fact that 
$G_{\Bbb F_p}$ ($\simeq \hat\Bbb Z$) has a nontrivial pro-$N$-Sylow group 
($\simeq \Bbb Z_N$), while $(\hat\Bbb Z^{(\Sigma')^{\dag}})^{\times}$ 
($=\prod_{\ell\in(\Sigma')^{\dag}}\Bbb Z_{\ell}^{\times}$) has trivial pro-$N$-Sylow group. 
(For the latter, observe that, for each $\ell\in\Sigma'$, 
$\Bbb Z_{\ell}^{\times}$ has trivial pro-$N$-Sylow group, since 
$\ell\not\equiv 1 \pmod{N}$, and $\ell\neq N$.) 

Taking $\Sigma'$ in this example as $\Sigma$, we also see that 
the implication (iv$'$)$\implies$(iii$'$) in 3.4.1 does not always hold. 
\enddefinition

\definition {3.4.4} The implication (iii$'$)$\implies$(iv$'$) in 3.4.1 does not always hold. 
In fact, there exists a set of prime numbers $\Sigma$ 
of density $0$, such that $\Sigma$ is not $\Bbb F_p$-small. 
To see this, take a sequence of positive integers $N_1\mid N_2\mid \cdots\mid N_k\mid\cdots$ 
such that $N_k \to\infty\  (k\to\infty)$, hence $\varphi(N_k) \to\infty\  (k\to\infty)$. 
Identify $\{1,\dots,N\}\isom\Bbb Z/N\Bbb Z$ by $i\mapsto i\mod N$ and 
set
$$I(N)\defeq\{i\in \Bbb Z/N\Bbb Z\mid \text{$\Sigma_i(N)$ is not $\Bbb F_p$-small}\}.$$
(See 3.4.3 for the definition of $\Sigma_i(N)$.) Then, by [Grunewald-Segal], Theorem A, 
$I(N)\neq\emptyset$. It is easy to see that $(I(N_k))_{k=1,2,\dots}$ is a projective subsystem 
of $(\Bbb Z/N_k\Bbb Z)_{k=1,2,\dots}$. As $I(N_k)\neq\emptyset$ for all $k\geq 1$, we 
see that $\varprojlim I(N_k) \neq\emptyset$. Fix any element $(i_k)_{k=1,2,\dots}$ of 
this projective limit. Then, for each $k\geq 1$, $\Sigma_{i_k}(N_k)$ is not 
$\Bbb F_p$-small, hence there exist $r_k\geq 1$, 
$\ell_{k,1},\dots,\ell_{k,r_k}\in\Sigma_{i_k}(N_k)^{\dag}$, and 
$e_{k,1},\dots,e_{k,r_k}\in\Bbb Z_{>0}$, such that the order of 
$p\mod \ell_{k,1}^{e_{k,1}}\cdots \ell_{k, r_k}^{e_{k,r_k}}$ in 
$(\Bbb Z/\ell_{k,1}^{e_{k,1}}\cdots \ell_{k, r_k}^{e_{k,r_k}}\Bbb Z)^{\times}$ is divisible 
by $k!$. Now, set 
$$\Sigma\defeq\{\ell_{k,j}\mid k\geq 1, 1\leq j\leq r_k\}.$$
Then it follows easily from the construction that $\Sigma$ is not $\Bbb F_p$-small. 
On the other hand, for each $k\geq 1$, we have 
$$\Sigma\subset \Sigma_{i_k}(N_k)\cup\{\ell_{k',j}\mid 1\leq k'<k, 1\leq j\leq r_k\}.$$
Since the density of $\Sigma_{i_k}(N_k)$ is (at most) $1/\varphi(N_k)$, we see that 
$\Sigma$ must be of density $0$, as desired. 

Taking $\Sigma'$ in this example as $\Sigma$, we also see that 
the implication (iv)$\implies$(iii) in 3.4.1 does not always hold. 

In particular, the implication (iii$'$)$\implies$(iii) in 3.4.1 does not always hold. 
(Namely, there exists a set of prime numbers which is neither $k$-large nor $k$-small.)
Indeed, if it held, then, combining it with the implication (iii)$\implies$(iv$'$), 
we would have the implication (iii$'$)$\implies$(iv$'$), which is absurd. 
\enddefinition

Next, and throughout the paper, for a subfield 
$\kappa \subset \bar k$, 
we write 
$\kappa ^{\times}\{\Sigma'\}\subset \kappa^{\times}$
for the $\Sigma '$-primary part of the (torsion) multiplicative group $\kappa ^{\times}$ 
and  
$(\kappa ^{\times})^{\Sigma}\defeq \kappa^{\times}/\kappa^{\times}\{\Sigma'\}$ 
for the maximal $\Sigma$-primary quotient of $\kappa ^{\times}$.

Let $X$ be a proper, smooth, and geometrically connected 
curve over $k$. In the following discussion $f,g$: $X\to\Bbb P^1$ will 
be non-constant $k$-morphisms. Define the open 
subschemes $U\defeq X\setminus (f^{-1}(\infty)\cup g^{-1}(\infty))$ and 
$U'\defeq U\setminus (f^{-1}(0)\cup g^{-1}(0))$ of $X$. We have the following commutative diagram:

$$\matrix
(f,g): 
&X &\to &\Bbb P^1_k \times \Bbb P^1_k \\
& \cup && \cup \\
&U &\to &\Bbb A^1_k \times \Bbb A^1_k \\
& \cup && \cup \\
&U' &\to & \ \Bbb G_{m,k} \times \Bbb G_{m,k} 
\endmatrix
$$
where $(f,g): X \to \Bbb P^1_k \times \Bbb P^1_k$ is the natural morphism determined by 
$f$ and $g$, the vertical inclusions are the natural open immersions, and the squares 
are fiber products. 

\definition 
{Definition/Proposition 3.5} {\sl We say that the pair $(f,g)$ has the property $P_\Sigma$
(respectively, $Q_\Sigma$, $Q_{0, \Sigma}$, $Q_{1,\Sigma}$ and $Q_{\infty, \Sigma}$)
if the following holds:

\medskip
\noindent
$P_\Sigma(f,g)$: 
$\exists a,b \in k^{\times}\{\Sigma'\}$, such that $f=a+bg$.

\smallskip\noindent
$Q_\Sigma(f,g)$: 
$\forall x\in U^{\cl}$, $\exists a_x,b_x\in k(x)^{\times}\{\Sigma'\}$, such that $f(x)=a_x+b_x g(x)$. 

\smallskip\noindent
$Q_{0, \Sigma}(f,g)$: 
$\forall x\in (U')^{\cl}$, $\exists a_x,b_x\in k(x)^{\times}\{\Sigma'\}$, such that $f(x)=a_x+b_x g(x)$. 

\smallskip\noindent
$Q_{1,\Sigma}(f,g)$: 
$\forall' x\in U^{\cl}$, $\exists a_x,b_x\in k(x)^{\times}\{\Sigma'\}$, such that $f(x)=a_x+b_x g(x)$. 

\smallskip\noindent
$Q_{\infty, \Sigma}(f,g)$: 
$\exists\infty\  x\in U^{\cl}$, for which 
$\exists a_x,b_x\in k(x)^{\times}\{\Sigma'\}$ such that $f(x)=a_x+b_x g(x)$. 

\medskip\noindent
Here the sign $\forall '$ means ``for all but finitely many'' and the sign 
$\exists\infty$ means ``there exist infinitely many''. 

Further, We say that the pair $(f,g)$ has the property $\overline P_\Sigma$
(respectively, $\overline Q_\Sigma$, $\overline Q_{0, \Sigma}$, $\overline Q_{1,\Sigma}$ and 
$\overline Q_{\infty, \Sigma}$) if the following holds:

\medskip
\noindent
$\overline P_\Sigma(f,g)$: 
$\exists a,b \in \bar k^{\times}\{\Sigma'\}$, such that $f=a+bg$.

\smallskip
\noindent
$\overline Q_\Sigma(f,g)$: 
$\forall x\in U^{\cl}$, $\exists a_x,b_x\in \bar k^{\times}\{\Sigma'\}$, such that $f(x)=a_x+b_x g(x)$. 

\smallskip\noindent
$\overline Q_{0, \Sigma}(f,g)$: 
$\forall x\in (U')^{\cl}$, $\exists a_x,b_x\in \bar k^{\times}\{\Sigma'\}$, such that $f(x)=a_x+b_x g(x)$. 

\smallskip\noindent
$\overline Q_{1, \Sigma}(f,g)$: 
$\forall' x\in U^{\cl}$, $\exists a_x,b_x\in \bar k^{\times}\{\Sigma'\}$, such that $f(x)=a_x+b_x g(x)$. 

\smallskip\noindent
$\overline Q_{\infty,\Sigma}(f,g)$: 
$\exists\infty\  x\in U^{\cl}$, for which 
$\exists a_x,b_x\in \bar k ^{\times}\{\Sigma'\}$, such that $f(x)=a_x+b_x g(x)$. 

\medskip
Then we have the following implications: 

$$
\matrix
P_\Sigma(f,g) &\iff & \overline P_\Sigma(f,g) \\
\Downarrow && \Downarrow \\
Q_\Sigma(f,g) &\implies &\overline Q_\Sigma(f,g) \\
\Downarrow && \Downarrow \\
Q_{0, \Sigma}(f,g) &\implies &\overline Q_{0, \Sigma}(f,g) \\
\Downarrow && \Downarrow \\
Q_{1,\Sigma}(f,g) &\implies &\overline Q_{1, \Sigma}(f,g) \\
\Downarrow && \Downarrow \\
Q_{\infty, \Sigma}(f,g) &\implies &\overline Q_{\infty,\Sigma}(f,g) 
\endmatrix
$$}

\enddefinition
\noindent

\demo
{Proof} For the proof of ``$\impliedby$'' in the first row, consider the action of 
$\Gal (K\bar k/K)\simeq G_k$, where $K$ denotes the function field of
$X$, and resort to the fact that $g$ is non-constant. 
The remaining implications are immediate.
\qed
\enddemo

\proclaim
{Proposition 3.6} Assume that $\Sigma \cup \{p\}\subsetneq \Primes$. Then the property 
$Q_{\infty, \Sigma}(f,g)$ (hence also the property $\overline Q_{\infty,\Sigma}(f,g)$) 
always holds. 
\endproclaim

\demo
{Proof} First note that the condition $\Sigma \cup \{p\}\subsetneq \Primes$ is equivalent to
saying that $\bar k^{\times}\{\Sigma'\}$ is an infinite set. For the proof it suffices to consider the 
following three cases: 1) $f-g$ is a non-constant function. 2) $\frac {f-1}{g}$ is a non-constant function,
and finally 3) both $f-g$ and $\frac {f-1}{g}$ are constant functions.
In case 1, consider the non-constant (hence dominant) 
$k$-morphism $f-g: U\to \Bbb A^1_k$. For all but 
finitely many $a\in \bar k^{\times}\{\Sigma'\}\subset \bar k=\Bbb A^1(\bar k)$, 
there exists $\bar x\in U(\bar k)$ that maps to $a$. Then, 
for the image $x$ of $\bar x$ in $U$, we have 
$f(x)-g(x)=a$, or $f(x)=a+1\cdot g(x)$. This completes 
the proof of case 1, as the equality $f(x)-g(x)=a$ also shows 
$a\in k(x)$. The proof of case 2 is similar to that of case 1: 
consider the non-constant $k$-morphism $\frac {f-1}{g}: U'\to \Bbb A^1_k$ 
and take a point in the fiber at 
$b\in \bar k^{\times}\{\Sigma'\}\subset \bar k=\Bbb A^1(\bar k)$. 
In case 3, we have $f=a_0+g=1+b_0g$ for some $a_0,b_0\in k$. 
As $g$ is non-constant, the second equality forces $a_0=b_0=1$, 
or, equivalently, $f=1+g$. Thus, $P_\Sigma(f,g)$ holds, hence, a fortiori, 
$Q_{\infty,\Sigma}(f,g)$ holds
\qed
\enddemo

\proclaim
{Proposition 3.7} 
\noindent
{\rm (i)} Assume that $\Sigma$ is $k$-small. Then the property $\overline Q_{0,\Sigma}(f,g)$ holds.

\noindent
{\rm (ii)} Assume that $\Sigma$ is finite. Then the property $Q_{1,\Sigma}(f,g)$ holds.
\endproclaim

\demo
{Proof} Fix $x\in U'$ and write $c=f(x), d=g(x) \in k(x) ^{\times} \subset\bar k ^{\times}$.

\noindent
(i) Since $\Sigma$ is $k$-small, 
we have 
$k\subset \exists k' \subset \bar k$, such that 
$(G_k:G_{k'})=\infty$ and that 
$(\chi_\Sigma(G_k): \chi_\Sigma(G_{k'}))<\infty$. 
Here, the first property says that $[k': k]=\infty$, 
while the second implies that
$N'\defeq \sharp((k')^{\times}\{\Sigma\})<\infty$. 
Replacing $k'$ by the finite extension $k'(c,d)$, 
we may assume that $k(c,d)\subset k'$. 
Consider 
the $k(c,d)$-curve 
$$Z_{N'}\defeq \{(u,v)\mid c=u^{N'}+dv^{N'}\} 
\subset \Bbb G_m\times\Bbb G_m.$$ 
This is a twist of the $N'$-th 
Fermat curve (minus cusps), hence, in particular, it is smooth and 
geometrically connected. Thus, (by means of the Weil bound) we have 
$\sharp(Z_{N'}(k'))=\infty$. Take $(u_0,v_0)\in 
Z_{N'}(k')$ and set $a\defeq u_0^{N'}, 
b\defeq v_0^{N'}$. Then we have $c=a+bd$. This completes 
the proof, since we have $a,b\in (k')^{\times}  \{\Sigma'\}(\subset \bar k^{\times} 
\{\Sigma'\})$ by the definition of $N'$. 

\noindent 
(ii) The proof of (ii) is similar to (i) 
but a little bit more subtle. First, 
for each $n\in\Bbb Z_{>0}$, we define 
$n_{\Sigma}$  
to be the greatest divisor of $n$ all of 
whose prime divisors belong to $\Sigma$.
Next, set $q\defeq q_{c,d}\defeq\sharp(k(c,d))$ and 
$N\defeq N_{c,d}\defeq \sharp(k(c,d)^{\times}\{\Sigma\})$. Thus, 
we have $N=(q-1)_{\Sigma}$. 

As in (i), consider the $k(c,d)$-curve 
$$Z_{N}\defeq \{(u,v)\mid c=u^{N}+dv^{N}\} \subset \Bbb G_m\times\Bbb G_m.$$ 
This is a twist of the $N$-th Fermat curve (minus cusps), 
hence, in particular, it is smooth and geometrically 
connected. The genus $g$ of $Z_{N}$ equals $(N-1)(N-2)/2$, and 
the cardinality $r$ of the set of geometric points which are cusps is 
$3N$. Thus, by means of the Weil bound, 
we have $\sharp(Z_{N}(k(c,d)))>0$, if 
$1+q-2g\sqrt{q} -r>0$. This last inequality can be rewritten as: 
$$q>N^2\sqrt{q}-\{(3N-2)(\sqrt{q}-1)-1\}.$$
Thus, it holds if $q\geq 4$ and $q>N^2\sqrt{q}$ hold, 
or, equivalently, if $q\geq 4$ and $q>N^4$ hold. By Lemma 3.8
below, these inequalities are satisfied (hence $\sharp(Z_{N}(k(c,d))>0$) 
for all but finitely many $q=p^m$. (Here, we resort to the fact that 
$\Sigma$ is finite.) 

Thus, for all but finitely many pairs $(c,d)$, 
$\sharp(Z_{N}(k(c,d)))>0$ holds. For such $(c,d)$, 
take $(u_0,v_0)\in Z_N(k(c,d))$ and set $a\defeq u_0^{N}, 
b\defeq v_0^{N}$. Then we have $c=a+bd$. This completes 
the proof, since we have $a,b \in k(c,d)^{\times}\{\Sigma'\} 
(\subset k(x)^{\times}\{\Sigma'\})$ 
by the definition of $N=N_{c,d}$. 
\qed
\enddemo

\proclaim{Lemma 3.8}
Let $p$ be a prime number, and $\Sigma$ a finite subset of $\Primes$. 
Then there exists a constant $C>0$ depending on $p$ and $\Sigma$, such that, 
for all $m\in\Bbb Z_{>0}$, $(p^m-1)_\Sigma\leq Cm$ holds. (For the notation 
$n_\Sigma$, see the proof of Proposition 3.7 (ii).) 
\endproclaim

\demo{Proof}
Let $f$ be the order of the image of $p$ in the multiplicative 
group $\prod_{\ell\in\Sigma^{\dag}}(\Bbb Z/\ell^{\epsilon_{\ell}}\Bbb Z) ^{\times}$, 
where 
$\Sigma^{\dag}\defeq\Sigma\setminus\{p\}$ 
and $\epsilon_{\ell}\defeq 1$ (resp. $2$) for $\ell\neq 2$ (resp. $\ell=2$). 
Then 
$$
(p^m-1)_\Sigma
\leq(p^{fm}-1)_\Sigma
= (p^f-1)_\Sigma\cdot m_{\Sigma^{\dag}}
\leq (p^f-1)_\Sigma \cdot m.
$$
Here, the first inequality follows from the fact that $p^m-1$ divides $p^{fm}-1$ 
and the equality is obtained by considering the structure of the multiplicative group 
$\Bbb Z_{\ell}  ^{\times}$ for $\ell\in\Sigma^{\dag}$. (More precisely, 
we have an isomorphism 
$1+\ell^{\epsilon_{\ell}}\Bbb Z_{\ell}\isom \ell^{\epsilon_{\ell}}\Bbb Z_{\ell}$ 
(say, the $\ell$-adic logarithm), 
which maps $1+\ell^e\Bbb Z_{\ell}$ onto $\ell^e\Bbb Z_{\ell}$ for each $e\geq\epsilon_{\ell}$. 
It follows from this that $a\in (1+\ell^e\Bbb Z_{\ell})\setminus (1+\ell^{e+1}\Bbb Z_{\ell})$ 
implies $a^m\in (1+m\ell^e\Bbb Z_{\ell})\setminus (1+m\ell^{e+1}\Bbb Z_{\ell})$, as desired.) 
Thus, $C\defeq (p^f-1)_\Sigma$ satisfies the desired property. 
\qed
\enddemo

\definition
{Remark 3.9}
The proof of Proposition 3.7 (i) can be viewed as a down-to-earth, (2-dimensional) 
torus version of the proof of [Raynaud], Proposition 2.2.1. 
\enddefinition

\definition
{Remark 3.10} (i) The proof of Proposition 3.7 (ii) shows that we may replace the 
assumption that $\Sigma$ is finite by the following: For all $m\gg 0$, 
$(p^m-1)_{\Sigma}<p^{m/4}$ holds. 

\noindent
(ii) Under the weaker assumption that $\Sigma$ is $k$-small, 
$Q_{1,\Sigma}(f,g)$ does not always hold. 
To construct a counterexample, set $k=\Bbb F_p$ and 
consider a prime number $r\neq p,\ r\nmid p-1$ and define $\Sigma$ to be 
the set of prime numbers dividing 
$p^{r^m}-1$ for some $m\geq 0$. Then, as in 3.4.2, 
$\Sigma$ is (infinite and) $k$-small. 
We define $k'$ to be the union of 
the finite fields $\Bbb F_{p^{r^m}}$ ($m\in\Bbb Z_{\geq 0}$). 
(Namely, $k'$ is the unique $\Bbb Z_r$-extension of 
the finite field $k$.) 
By definition, we have $(k')^{\times}\{\Sigma'\}=\{1\}$. 
Now, take any $X$, $f$, $g$ as above such that $f\neq 1+g$. Then 
$U_1\defeq \{x\in U\mid f(x)\neq 1+g(x)\}$ 
is a non-empty open subset of $X$. Thus, (by the Weil bound) 
we have $\sharp(U_1(k'))=\infty$. Moreover, for any 
$x$ in the image of $U_1(k')$ in $U$, there does not exist 
$a_x,b_x\in k(x)^{\times}\{\Sigma'\}$ such that $f(x)=a_x+b_x g(x)$. 
(Observe $k(x)^{\times}\{\Sigma'\}\subset (k')^{\times}\{\Sigma'\}=\{1\}$.) Thus, 
$Q_{1,\Sigma}(f,g)$ does not hold. 
\enddefinition

The following is the main result in this section, which plays a crucial role in the proof of
the main Theorem 4.1 of this paper.

\proclaim
{Proposition 3.11} Assume that $\Sigma$ is $k$-large. Then the implication  
$$\overline Q_{1,\Sigma}(f,g)\implies \overline P_{\Sigma}(f,g)$$ 
holds. 
\endproclaim

\demo
{Proof}
For each non-constant rational function $h$ 
on $\overline X\defeq X\times_k\bar k$, 
we define $\deg(h)$ to be the degree of the non-constant 
$\bar k$-morphism $h:\overline X\to\Bbb P^1$ associated with $h$ (or, equivalently, 
the degree of the pole divisor $(h)_\infty$). 
Set $d=\deg(f)+\deg(g)$. Then, for any $a,b\in\bar k$, 
either $f-(a+bg)$ is a constant function or $\deg(f-(a+bg))\leq d$. 
Assume that the property $\overline{Q}_{1,\Sigma}(f,g)$ holds. 
Then there exists a non-empty open subscheme $U_2\subset U$, such that 
$\forall x\in (U_2)^{\cl}$, $\exists a_x,b_x\in \bar k^{\times}\{\Sigma'\}$ 
such that the equality $f(x)=a_x+b_xg(x)$ holds. First, consider the case 
where $f-(a_x+b_x g)$ is constant for some $x$. 
Then, by evaluating at $x$, we see that this constant 
must be $0$. Namely, $f=a_x+b_x g$ holds, which implies that the property 
$\overline P_{\Sigma}(f,g)$ holds, as desired. So, suppose that 
$f-(a_x+b_x g)$ is non-constant for any $x$. 
Then the non-constant morphism $f-(a_x+b_x g): \overline X\to \Bbb P^1$ 
is defined over $k(a_x,b_x)\subset k_{\Sigma'}$. (For the last inclusion, 
note that $\overline k^{\times}\{\Sigma'\}\subset k_{\Sigma'}$ by definition.)  
Considering the fiber at $0$ of this non-constant morphism 
over $k(a_x,b_x)\subset k_{\Sigma'}$, we deduce: 
$$[k_{\Sigma'}(x):k_{\Sigma'}]
\leq [k(a_x,b_x)(x):k(a_x,b_x)]\leq\deg(f-(a_x+b_x g))\leq d.$$
Now, since $G_{k_{\Sigma'}}$ is (pro)cyclic as a closed subgroup of 
$G_k\simeq\hat \Bbb Z$, we conclude that there exists a finite extension 
$k'$ of $k_{\Sigma'}$, such that 
$k(x)\subset k_{\Sigma'}(x)\subset k'$ holds for any $x\in U_2^{\cl}$. 
By the assumption that $\Sigma$ is $k$-large, we have 
$k_{\Sigma'}\subsetneq \bar k$, hence $[\bar k:k_{\Sigma'}]=\infty$. 
(Observe that $G_k$ does not admit a nontrivial finite subgroup.) 
So, we also have $k'\subsetneq \bar k$. This contradicts the 
previous conclusion. Indeed, since $U_2$ is an affine curve over $k$, 
it admits a finite $k$-morphism $\phi: U_2 \to\Bbb A^1$. Take 
$a\in \bar k\setminus k'\subset \bar k=\Bbb A^1(\bar k)$ and $x\in \phi^{-1}(a)$. 
Then we have $a\in k(a)\subset k(x)\subset k'$, which is absurd. 
\qed
\enddemo

We will also use the following slight generalization of Proposition 3.11 later.

\definition
{Definition/Proposition 3.12}  {\sl For a pair $(f,g)$ as in the discussion before 
Definition/Proposition 3.5, a positive integer $m$, and a set of prime numbers 
$\Sigma \subset \Primes$, we define the following properties: 

\medskip\noindent
$P_\Sigma ^{(m)}(f,g)$: 
$\exists a,c \in k^{\times}\{\Sigma'\}$, such that $f=a(1+cg)^m$.

\smallskip\noindent
$\overline P_\Sigma ^{(m)}(f,g)$: 
$\exists a,c \in \bar k^{\times}\{\Sigma'\}$, such that $f=a(1+cg)^m$.

\smallskip\noindent
$\overline Q_{1, \Sigma}^{(m)}(f,g)$: 
$\forall' x\in U$, $\exists a_x,c_x\in \bar k^{\times}\{\Sigma'\}$, such that 
$f(x)=a_x(1+c_x g(x))^m$. 

\medskip
Then: 

\noindent
(i) The implications 
$$
P_{\Sigma}^{(m)}(f,g)
\iff \overline P_{\Sigma}^{(m)}(f,g)
\implies \overline Q_{1,\Sigma}^{(m)}(f,g)
$$ 
hold. 

\noindent 
(ii) If $\Sigma$ is $k$-large, then the implication  
$$\overline Q_{1,\Sigma}^{(m)}(f,g)\implies \overline P_{\Sigma}^{(m)}(f,g)$$ 
holds.} 
\enddefinition

\demo
{Proof} 
(i) Similar to the proof of Definition/Proposition 3.5. 

\noindent 
(ii) Similar to the proof of Proposition 3.11.
\qed
\enddemo

\subhead 
\S 4. The Main Theorem
\endsubhead
In this section we state and prove our main result. 
We follow the notations in $\S 1$, $\S2$, and $\S 3$.

Let $k$, $l$ be finite fields of characteristic 
$p_k$, $p_l$, respectively, and of cardinality $q_k$, $q_l$, respectively. 
Let $X$, $Y$ be smooth, proper, and geometrically connected 
curves over $k$, $l$, respectively. 
Let $K$, $L$ be the function fields of $X$, $Y$, respectively. 
We will write $G_K\defeq \Gal (K^{\sep}/K)$
for the absolute Galois group of $K$, and similarly $G_L=\Gal (L^{\sep}/L)$ 
for the absolute Galois group of $L$. 

Let $\Sigma_X, \Sigma_Y\subset \Primes$
be sets of prime numbers. 
Write $G_K^{(\Sigma_X)}$ (resp. $G_L^{(\Sigma_Y)}$)
for the  maximal geometrically pro-$\Sigma_X$ (resp. pro-$\Sigma_Y$) quotient of $G_K$ (resp. $G_L$). 
Thus, we have exact sequences:
$$1\to {\overline G}_K^{\Sigma _X}\to G_K^{(\Sigma _X)}@>\pr>> G_{k}\to 1,$$
resp.
$$1\to {\overline G}_L^{\Sigma _Y}\to G_L^{(\Sigma _Y)}@>\pr>> G_{l}\to 1,$$
where 
$G_{k}\defeq \Gal (\bar k/k)$ 
(resp. $G_{l}\defeq \Gal (\bar l/l)$) 
is the absolute Galois group of $k$ (resp. $l$), and
${\overline G}_K^{\Sigma_X}$ 
(resp. ${\overline G}_L^{\Sigma_Y}$) 
is the maximal pro-$\Sigma _X$ 
(resp. pro-$\Sigma_Y$) quotient
of the absolute Galois group ${\overline G}_K\defeq\Gal (K^{\sep}/K\bar k)$ 
(resp. ${\overline G}_L\defeq\Gal (L^{\sep}/L\bar l)$ 
of $K\bar k$ (resp. $L\bar l$). 
Our aim in this section is to prove the following Theorem:

\proclaim
{Theorem 4.1} Assume that $\Sigma_X$ is $k$-large (cf. Definition 3.2). 
Assume also that $\Sigma_X$ satisfies condition $(\epsilon_{X})$  
(cf. the discussion before Theorem C in \S 0). 
Let
$$\sigma:G_K^{(\Sigma_X)} \isom G_L^{(\Sigma_Y)}$$ 
be an isomorphism between profinite groups. Then $\sigma$ arises from a uniquely 
determined commutative diagram of field extensions:

$$
\CD
L\sptilde @>{\sim}>> K\sptilde \\
@AAA              @AAA \\
L@>{\sim}>> K \\
\endCD
$$
in which the horizontal arrows are isomorphisms, and the vertical arrows are the field 
extensions corresponding to the groups $G_L^{(\Sigma_Y)}$ and $G_K^{(\Sigma_X)}$, respectively. 
Thus, $L\sptilde/L$ (resp. $K\sptilde/K$) is the subextension of $L^{\sep}/L$ (resp. $K^{\sep}/K$)
with Galois group $G_{L}^{(\Sigma_Y)}$ (resp. $G_{K}^{(\Sigma_X)}$).
\endproclaim

For the rest of this section we will consider an isomorphism of profinite groups
$$\sigma : G_K^{(\Sigma _X)}\isom G_L^{(\Sigma _Y)}.$$ 
We write $\tilde X$ (resp. $\tilde Y$) for the normalization of $X$ (resp. $Y$) in 
$K\sptilde$ (resp. $L\sptilde$).

We already know the following: $\Sigma\defeq \Sigma_X=\Sigma_Y$ (cf. Lemma 2.1), 
$p\defeq p_k=p_l$, and $q\defeq q_k=q_l$ (cf. Lemma 2.4.7 (i)). 
Moreover, there exists a bijection $\tilde \phi:\tilde X^{\cl}\isom \tilde Y^{\cl},\ \tilde x\mapsto \tilde y$,
such that $\sigma(D_{\tilde x})=D_{\tilde y}$, 
which naturally induces a bijection $\phi: X^{\cl}\isom Y^{\cl}$ and 
an isomorphism $\phi:\Div_X\isom \Div _Y$ (cf. Lemma 2.2). 

\proclaim
{Lemma 4.2} (Invariance of Global Modules of Roots of Unity) 
The isomorphism $\sigma$ induces naturally an isomorphism:
$$M_{X}^{\Sigma}\isom M_{Y}^{\Sigma}$$
between the (global) modules of roots of unity which is Galois-equivariant with respect to $\sigma$. 
\endproclaim

\demo
{Proof} Let $J_{K}\defeq \prod'_{x\in X^{\cl}}K_x^{\times}$ be the id\`ele group of $K$ and 
$J_{K}^{(\Sigma)} \defeq  \prod'_{x\in X^{\cl}}(K_x^{\times})^{(\Sigma)}$ (cf. discussion before Lemma 2.3 for the definition of 
$(K_x^{\times})^{(\Sigma)}$ and the various notations below) which is a quotient of $J_K$. The  Artin reciprocity map $\psi _K:J_K\to G_K^{\ab}$ of global class field theory 
induces naturally a map $\psi_K^{(\Sigma)}: J_K^{(\Sigma)}\to G_K^{(\Sigma),\ab}$. 
(When $H$ is a profinite group, $H^{\ab}$ 
denotes the maximal abelian quotient of $H$.) 
The exact sequence 
$$1\to K^{\times}\to J_K @>\psi_K>>  G_K^{\ab}$$ 
from global class field theory induces naturally an exact sequence
$$1\to k^{\times}\to \prod _{x\in X^{\cl}}\Cal O_x^{\times} \to G_K^{\ab}\to \pi_1(X)^{\ab}\to 0,$$ 
where the map $G_K^{\ab}\to \pi_1(X)^{\ab}$ is the natural one, the map $\prod _{x\in X^{\cl}}\Cal O_x^{\times} \to G_K^{\ab}$
is the restriction of the reciprocity map $\psi_K$, and the map $k^{\times}\to {\underset {x\in X^{\cl}}\to 
\prod} \Cal O_{x}^{\times}$ is the natural diagonal embedding.  
(Here we recall that for each $x\in X^{\cl}$, the map $\psi_X:J_K
\to G_K^{\ab}$ maps the component $K_x^{\times}$ of $J_K$ into the decomposition group 
$D_x^{\ab}\subset G_K^{\ab}$ associated to $x$ via the local reciprocity map $K_x^{\times}\to D_x^{\ab}$.)
This latter sequence induces naturally an exact sequence
$$1\to (k^{\times})^{\Sigma}\to {\underset {x\in X^{\cl}}\to \prod} 
(\Cal O_{x}^{\times})^{\Sigma} \to G_K^{(\Sigma), \ab}\to \pi _1(X)^{(\Sigma),\ab}\to 0,$$ 
where the map ${\underset {x\in X^{\cl}}\to \prod} 
(\Cal O_{x}^{\times})^{\Sigma} \to G_K^{(\Sigma), \ab}$ is naturally induced by the above map
$\psi_K^{(\Sigma)}:J_K^{(\Sigma)}\to G_K^{(\Sigma),\ab}$, and the map $(k^{\times})^{\Sigma}\to {\underset {x\in X^{\cl}}\to 
\prod} (\Cal O_{x}^{\times})^{\Sigma}$ is the natural diagonal embedding.  
Further, we have the following commutative diagram:

$$
\CD
1@>>>(k^{\times})^{\Sigma}@>>> {\underset {x\in X^{\cl}}\to \prod} 
(\Cal O_{x}^{\times})^{\Sigma} @>>> G_K^{(\Sigma), \ab}@>>> \pi _1(X)^{(\Sigma),\ab}@>>> 0  \\
@.  @VVV            @VVV           @VVV           @VVV   \\
1@>>>(l^{\times})^{\Sigma}@>>> {\underset {y\in Y^{\cl}}\to \prod} 
(\Cal O_{y}^{\times})^{\Sigma} @>>> G_L^{(\Sigma),\ab}@>>> \pi _1(Y)^{(\Sigma),\ab}@>>> 0
\endCD \tag{4.1}
$$
where the map $G_K^{(\Sigma), \ab}\to G_L^{(\Sigma), \ab}$ is naturally induced by $\sigma:G_K^{(\Sigma)} \isom G_L^{(\Sigma)}$ (hence is an isomorphism),
and the map  ${\underset {x\in X^{\cl}}\to \prod} (\Cal O_{x}^{\times})^{\Sigma}\to  {\underset {y\in Y^{\cl}}\to \prod} (\Cal O_{y}^{\times})^{\Sigma}$ 
maps each component $(\Cal O_{x}^{\times})^{\Sigma}$ isomorphically onto 
$(\Cal O_{y}^{\times})^{\Sigma}$ , where
$y\defeq \phi (x)$  
(cf. Lemma 2.2 and Lemma 2.5). 
In particular, this map is an isomorphism since $\phi: X^{\cl}\isom Y^{\cl}$ 
is a set-theoretic bijection. 
Thus, the far left vertical map in the diagram (4.1) gives an isomorphism 
$(k^{\times})^{\Sigma}\isom (l^{\times}) ^{\Sigma}$. Passing to the open subgroups of  
$G_K^{(\Sigma)}$ and $G_L^{(\Sigma)}$, corresponding to extensions of the constant fields, 
to corresponding diagrams (4.1), and to the 
projective limits via the 
natural maps, we obtain the desired 
isomorphism $M_{X}^{\Sigma}\isom M_{Y}^{\Sigma}$, which is 
Galois-equivariant with respect to $\sigma$ as is easily verified by construction.
\qed
\enddemo

\proclaim
{Lemma 4.3} (Rigidity of Inertia) Let $x\in X^{\cl}$ and $y\defeq \phi (x)$. 
The following diagram is commutative:

$$
\CD
M_{X}^{\Sigma}  @>>>   M_{\overline{k (x)}}^{\Sigma}\\
@VVV             @VVV   \\
M_{Y}^{\Sigma}  @>>>   M_{\overline{k(y)}}^{\Sigma}
\endCD 
$$
where the left vertical arrow is the isomorphism in Lemma 4.2, 
the right vertical arrow is the isomorphism in Lemma 2.3 (ii), and the 
horizontal maps are the natural identifications. Further, this 
diagram is Galois-equivariant with respect to the commutative diagram:

$$
\CD
D_{\tilde x}  @>>>     G_K^{(\Sigma)}  \\
@VVV                    @V{\sigma }VV  \\
D_{\tilde y}  @>>>     G_L^{(\Sigma)} 
\endCD
$$
where $\tilde x\in \tilde X^{\cl}$ is a point above $x\in X$, 
$\tilde y\defeq \tilde\phi(\tilde x)$, and 
the horizontal maps are the natural inclusions.
\endproclaim

\demo
{Proof}
Indeed, the far left square in the diagram (4.1) induces a  
diagram:
 $$
\CD
(k^{\times})^{\Sigma}@>>> (k(x)^{\times})^{\Sigma} \\
  @VVV            @VVV     \\
(l^{\times})^{\Sigma}@>>> (k(y)^{\times})^{\Sigma}
\endCD 
$$
where the vertical arrows are the isomorphisms induced by $\sigma$ 
and the horizontal arrows are the natural inclusions (cf. Lemma 
2.3 (i) and Lemma 4.2). 
Passing to the open subgroups of  $G_K^{(\Sigma),\ab}$  
and $G_L^{(\Sigma),\ab}$, corresponding to extensions of the constant fields, and 
corresponding diagrams as above, and to the projective limits, 
we obtain the desired Galois-equivariant diagram. 
(Observe that 
the above diagram is commutative 
cofinally.) 
\qed
\enddemo

For an abelian group $A$, write $A\{\Sigma'\}$ for 
the $\Sigma'$-primary part of the torsion subgroup of $A$, and 
set $A^n\defeq \{a^n\mid a\in A\}$ for each $n\in \Bbb Z_{>0}$. 
Applying the notation $H^{(\Sigma)}$ in the beginning of \S 1 
to the (discrete) group $H=K^{\times}$ 
(resp. $H=L^{\times}$) and 
a (pro)finite subgroup $\overline{H}=k^{\times}$ 
(resp. $\overline{H}=l^{\times}$), 
we have 
$(K^{\times})^{(\Sigma)}=
K^{\times} /(k^{\times}\{\Sigma'\})$ (resp. $(L^{\times})^{(\Sigma)}=
L^{\times} /(l^{\times}\{\Sigma'\})$). 

\proclaim
{Lemma 4.4} (A Power of the Multiplicative Group 
modulo $\Sigma'$-primary Torsion).
\noindent
{\rm (i)} We have 
$m\defeq \sharp(\pi _1(X)^{\ab,\tor}\{\Sigma'\})
=\sharp(\pi _1(Y)^{\ab,\tor}\{\Sigma'\})$.

\noindent
{\rm (ii)}
The isomorphism $\sigma$ induces naturally an injective homomorphism 
$$\gamma': ((K^{\times})^{(\Sigma)})^m
\hookrightarrow (L^{\times})^{(\Sigma)}$$ 
between multiplicative groups. 

\noindent
{\rm (iii)} The homomorphism $\gamma'$ fits into the following natural commutative 
diagram: 

$$
\CD
((K^{\times})^{(\Sigma)})^m   @>{\gamma'}>> (L^{\times})^{(\Sigma)} \\
 @VVV                           @VVV \\
(K^{\times}/k^{\times})^m   @>{\bar \gamma'}>> L^{\times}/l^{\times}\\
\endCD
$$
where the vertical maps are the natural surjective homomorphisms and 
\newline
$\bar \gamma':(K^{\times}/k^{\times})^m \hookrightarrow 
L^{\times}/l^{\times}$ is an injective homomorphism naturally induced by $\gamma'$. 
\endproclaim

\demo
{Proof} 
(i) As $\pi_1(X)^{\ab,\tor}
\isom J_{X}(k)$, we have 
$$\sharp(\pi _1(X)^{\ab,\tor}\{\Sigma'\})=\sharp(J_{X}(k)\{\Sigma'\})=
\sharp(J_{X}(k))_{\Sigma'},$$
where, for each $n\in\Bbb Z_{>0}$, we define 
$n_{\Sigma'}$ 
to be the greatest divisor of $n$ all of 
whose prime divisors belong to $\Sigma'$. Similarly
$$\sharp(\pi _1(Y)^{\ab,\tor}\{\Sigma'\})=\sharp(J_{Y}(l)\{\Sigma'\})=
\sharp(J_{Y}(l))_{\Sigma'}.$$
Thus, the assertion follows from 
Lemma 2.4.10. 

\noindent
(ii) We have the following commutative diagram:
$$
\CD
1 @>>> \Ker (\psi_K^{(\Sigma)}) @>>> J_{K}^{(\Sigma)} @>{\psi _K^{(\Sigma)}}>>  G_K^{(\Sigma),\ab} \\
@.           @VVV     @VVV     @VVV \\      
1 @>>> \Ker (\psi_L^{(\Sigma)}) @>>> J_{L}^{(\Sigma)} @>{\psi _L^{(\Sigma)}}>>  G_L^{(\Sigma),\ab} \\
\endCD 
\tag{4.2}
$$
where the horizontal rows are exact. Here, 
$J_{K}^{(\Sigma)}\defeq \prod_{x\in X^{\cl}}'(K_{x}^{\times})^{(\Sigma)}$ 
(resp. $J_{L}^{(\Sigma)}\defeq \prod_{y\in Y^{\cl}}'(L_{y}^{\times})^{(\Sigma)}$) 
is a quotient of the id\`ele group $J_{K}$ 
(resp. $J_{L}$) 
of $K$ (resp. $L$), and the map $\psi_K^{(\Sigma)}:J_{K}^{(\Sigma)} \to G_K^{(\Sigma),\ab}$ 
(resp. the map $\psi_L^{(\Sigma)}:J_{L}^{(\Sigma)} \to G_L^{(\Sigma),\ab}$)
is naturally induced by Artin's reciprocity map in global class 
field theory (cf. proof of Lemma 4.2). 
The far right vertical map is naturally induced by 
$\sigma$, and  the middle vertical map $J_{K}^{(\Sigma)}\to J_{L}^{(\Sigma)}$ maps each component 
$(K_{x}^{\times })^{(\Sigma)}$ 
isomorphically onto $(L_{y}^{\times})^{(\Sigma)}$; 
where $y\defeq \phi (x)$, via the natural identification 
in Lemma 2.5, which is induced by $\sigma$. 
In particular, the map $J_{K}^{(\Sigma)}\to J_{L}^{(\Sigma)}$ is an isomorphism.
Thus, the far left vertical map in the diagram (4.2) is a natural isomorphism
$\Ker (\psi_K^{(\Sigma)})\to \Ker (\psi_L^{(\Sigma)})$ between kernels. 
We claim: 
\definition{Claim 1} There exists a canonical exact sequence:
$$1\to (K^{\times})^{(\Sigma)}\to \Ker (\psi_K^{(\Sigma)})\to \pi _1(X)^{\ab,\tor}\{\Sigma'\}\to 0$$
(resp.
$$1\to (L^{\times})^{(\Sigma)}\to \Ker (\psi_L^{(\Sigma)})\to \pi _1(Y)^{\ab,\tor}\{\Sigma'\}\to 0).$$
\enddefinition

Assuming this claim, we then have a commutative diagram:
$$
\CD
1  @>>>   (K^{\times})^{(\Sigma)}  @>>> \Ker (\psi_K^{(\Sigma)})@>>> \pi _1(X)^{\ab,\tor}\{\Sigma'\} @>>> 
0\\
@.        @.          @VVV      @.    @. \\
1  @>>>   (L^{\times})^{(\Sigma)}  @>>> \Ker (\psi_L^{(\Sigma)}) @>>> \pi _1(Y)^{\ab,\tor}\{\Sigma'\} 
@>>> 0\\
\endCD
$$
where the horizontal rows are exact, and the vertical arrow is the
above isomorphism.
This isomorphism has, a priori, no reason to map 
$(K^{\times})^{(\Sigma)}$ into 
$(L^{\times})^{(\Sigma)}$. However, since $\pi _1(Y)^{\ab,\tor}\{\Sigma'\}$ is a finite 
abelian group of exponent dividing $m$, we can conclude that the above isomorphism 
$\Ker (\psi_K^{(\Sigma)})\to \Ker (\psi_L^{(\Sigma)})$ maps $((K^{\times})^{(\Sigma)})^{m}$
injectively into $(L^{\times})^{(\Sigma)}$. 
Thus, we obtain a natural injective map $\gamma' : 
((K^{\times})^{(\Sigma)})^{m}  \to 
(L^{\times})^{(\Sigma)}$. It remains to prove Claim 1. We will only prove the assertion concerning
$\Ker (\psi_K^{(\Sigma)})$ (the assertion concerning $\Ker (\psi_L^{(\Sigma)})$ is proved in a similar way).

We have the 
following commutative diagram:
$$
\CD
@. @.   1 @. 1 \\
@. @. @AAA    @AAA  \\
1 @>>>  \Ker (\rho _K) @>>>   \Im (\psi_K)(\subset G_K^{\ab}) @>{\rho _K}>>   
\Im (\psi_K^{(\Sigma)})(\subset G_K^{(\Sigma),\ab}) @>>> 1 \\
@.   @AAA    @A{\psi _K}AA  @A{\psi _K^{(\Sigma)}}AA  \\
1@>>> \prod _{x\in X^{\cl}} 
N_x
@>>> J_{K} @>>>J_{K}^{(\Sigma)}@>>> 1 \\ 
@. @AAA  @AAA   @AAA  \\
1 @>>>   k^{\times}\{\Sigma'\} @>>>  K^{\times}   @>>>  \Ker (\psi_K^{(\Sigma)})\\
@. @AAA @AAA  @AAA \\
@. 1 @.  1 @. 1\\
\endCD 
$$
where the vertical and horizontal rows are exact.
Here, the map $\psi_K:J_{K}\to G_K^{\ab}$ is Artin's reciprocity map in 
global class field theory, and the map
$\rho_K:\Im (\psi_K) \to \Im (\psi_K^{(\Sigma)})$ is the restriction of the natural map
$G_K^{\ab}\twoheadrightarrow G_K^{(\Sigma),\ab}$. Further, $J_{K}\to J_{K}^{(\Sigma)}$ is the 
natural map which maps each component $K_{x}^{\times}$ canonically 
onto $(K_{x}^{\times})^{(\Sigma)}=K_{x}^{\times}/N_x$ (cf. the discussion before Lemma 2.3 for the definition of $N_x$). 
In particular, we deduce that the cokernel of the injective map 
$(\prod _{x\in X^{\cl}} N_x)/(k^{\times}\{\Sigma'\}) \hookrightarrow \Ker (\rho _K)$ 
is naturally isomorphic to the cokernel of the injective map  
$(K^{\times})^{(\Sigma)}\hookrightarrow \Ker (\psi_K^{(\Sigma)})$. 
Observe that $\Ker (\rho_K)$ is naturally 
identified with the kernel $\Ker (G_K^{\ab}\twoheadrightarrow G_K^{(\Sigma), \ab})$.
Further, we claim: 

\definition {Claim 2}
The cokernel of the above injective homomorphism 
$(\prod _{x\in X^{\cl}}N_x)/(k^{\times}\{\Sigma'\})  
\hookrightarrow \Ker (\rho _K)$ is naturally 
isomorphic to $\pi _1(X)^{\ab,\tor}\{\Sigma'\}$. 
\enddefinition

Indeed, we have the following 
commutative diagram:
$$
\CD
@.   1  @.      1  @.    1@.  \\
@.     @VVV      @VVV    @VVV  \\        
1 @>>>  k^{\times}\{\Sigma'\}       @>>> k^{\times}  @>>>  (k^{\times})^{\Sigma} @>>> 1 \\
@.    @VVV    @VVV    @VVV   \\
1 @>>>  \prod _{x\in X^{\cl}} N_x   @>>> \prod _{x\in 
X^{\cl}} \Cal O_{x}^{\times}@>>> \prod _{x\in X^{\cl}} 
(\Cal O_{x}^{\times})^{\Sigma}  @>>> 1   \\
@.     @VVV    @VVV    @VVV    \\
1@>>>  \Ker (\rho_K) @>>>  G_K^{\ab} @>{\rho_K}>>  G_K^{(\Sigma),\ab} @>>>  1  \\ 
@.  @VVV     @VVV   @VVV    \\
1@>>>  \Ker (\nu_X) @>>> \pi _1(X)^{\ab}  @>{\nu _X}>>  \pi _1(X)^{(\Sigma),
\ab} @>>> 1 \\
@. @VVV  @VVV    @VVV  \\
@.  1 @.   1  @. 1 \\
\endCD 
$$
where the vertical and horizontal rows are exact. Here, the maps  
$G_K^{\ab}\to \pi _1(X)^{\ab}$ and $G_K^{(\Sigma),\ab}\to \pi _1(X)^{(\Sigma),\ab}$ are 
the natural maps, the map $\prod _{x\in X^{\cl}} \Cal O_{x}^{\times}
\to  G_K^{\ab}$ is the restriction of Artin's 
reciprocity map, and the map $ k^{\times}\to \prod _{x\in X^{\cl}} 
\Cal O_{x}^{\times}$ is the natural diagonal embedding.
Further, the kernel $\Ker (\nu_X)$ of $\nu_X$ is canonically isomorphic to
$\pi _1(X)^{\ab,\tor}\{\Sigma'\}$, as follows from the structure of $\pi _1(X)^
{\ab}$. Note that the maximal pro-$\Sigma$
quotient $\pi _1(X)^{\ab,\tor,\Sigma}$ of $\pi _1(X)^{\ab,\tor}$
is naturally isomorphic to the torsion subgroup $\pi _1(X)^{(\Sigma),\ab,\tor}$ 
of $\pi _1(X)^{(\Sigma),\ab}$.
Thus, Claim 2, hence Claim 1, are proved. This completes the proof of (ii). 

\noindent
(iii) This follows from the fact that 
$(K^{\times})^{(\Sigma)}$ (resp. ($L^{\times})^{(\Sigma)}$)
modulo its torsion subgroup 
(which is naturally identified with $(k^{\times})^{\Sigma}$ (resp. $(l^{\times})^{\Sigma}$)) 
is naturally identified with $K^{\times}/k^{\times}$ (resp. $L^{\times}/l^{\times}$).
\qed\enddemo

We have a commutative diagram:
$$
\CD
1  @>>>   (K^{\times})^{(\Sigma)}  @>>> \Ker (\psi_K^{(\Sigma)})@>>> \pi _1(X)^{\ab,\tor}\{\Sigma'\} @>>> 
0\\
@.        @.          @V{\rho}VV      @.    @. \\
1  @>>>   (L^{\times})^{(\Sigma)}  @>>> \Ker (\psi_L^{(\Sigma)}) @>>> \pi _1(Y)^{\ab,\tor}\{\Sigma'\} 
@>>> 0\\
\endCD
$$
where the horizontal rows are exact, and the vertical arrow is an isomorphism naturally induced by $\sigma$
(cf. proof of Lemma 4.4). Let
$$R_K\defeq \Ker (\psi_K^{(\Sigma)})/ (\Ker (\psi_K^{(\Sigma)})^{\tor}\{\Sigma\})$$
where $\Ker (\psi_K^{(\Sigma)})^{\tor}\{\Sigma\}$ is the group of 
$\Sigma$-primary torsion of $\Ker (\psi_K^{(\Sigma)})$, which is contained in  
$(K^{\times})^{(\Sigma)}$
(since $\pi _1(X)^{\ab,\tor}\{\Sigma'\}$ is $\Sigma'$-primary),
and is naturally identified with $(k^{\times})^{\Sigma}$. 
Thus, $R_K$ naturally inserts in the following exact sequence:
$$1\to K^{\times}/k^{\times}\to R_K \to \pi _1(X)^{\ab,\tor}\{\Sigma'\}\to 0.$$
We define $R_L$ in a similar way which sits in the following exact sequence:
$$1\to L^{\times}/l^{\times}\to R_L \to \pi _1(Y)^{\ab,\tor}\{\Sigma'\}\to 0.$$

The above isomorphism 
$$\rho:\Ker (\psi_K^{(\Sigma)})\isom \Ker (\psi_L^{(\Sigma)})$$ 
induces naturally a commutative diagram:
$$
\CD
1  @>>>   K^{\times}/k^{\times}  @>>> R_K @>>> \pi _1(X)^{\ab,\tor}\{\Sigma'\} @>>> 
0\\
@.        @.          @V{\bar \rho}VV      @.    @. \\
1  @>>>   L^{\times}/l^{\times}  @>>> R_L @>>> \pi _1(Y)^{\ab,\tor}\{\Sigma'\} 
@>>> 0\\
\endCD
$$
where the horizontal rows are exact, and the vertical arrow 
is an isomorphism. 
Further, define $\overline{H}_K\subset K^{\times}/k^{\times}$ to 
be the kernel of the composite homomorphism 
$K^{\times}/k^{\times}\hookrightarrow R_K\isom R_L \twoheadrightarrow 
\pi_1(Y)^{\ab,\tor}\{\Sigma'\}$, and set 
$\overline H_L\defeq \bar \rho (\overline H_K)$. 
Then it is easy to see that $\overline H_L\subset L^{\times}/l^{\times}$ and 
that $(K^{\times}/k^{\times}: \overline H_K)=(L^{\times}/l^{\times}: \overline H_L)$ 
divides $m\defeq \sharp(\pi _1(X)^{\ab,\tor}\{\Sigma'\})
=\sharp(\pi _1(Y)^{\ab,\tor}\{\Sigma'\})$. 
In particular, we have $(K^{\times}/k^{\times})^m\subset \overline H_K$. 

In the following, we will think of the elements of  
$$\PDiv_{X}\defeq K^{\times}/k^{\times}$$
(resp. $\PDiv_{Y}\defeq L^{\times}/l^{\times})$
as principal divisors of rational functions on $X$ (resp. $Y$), and denote them by $\bar f, \bar g,\dots$, where 
$f, g, \dots$ are rational functions on $X$ (resp. $Y$). We will also denote the elements of 
$(K^{\times})^{(\Sigma)}$ (resp. $(L^{\times})^{(\Sigma)}$)
by $f', g',\dots$, and refer to them as ``{\it pseudo-functions}''$\defeq$ 
classes of rational functions on $X$ (resp. $Y$) 
modulo constants in $k^{\times}\{\Sigma'\}$ (resp. $l^{\times}\{\Sigma'\}$).
We define
$$H_K'\defeq \{f'\in (K^{\times})^{(\Sigma)}\ \vert\ \bar f\in \overline H_K\},$$
and
$$H_K^{\times}\defeq \{f\in K^{\times}\ \vert\ \bar f\in \overline H_K\}.$$
We define $H_L'\subset (L^{\times})^{(\Sigma)}$, and $H_L^{\times}\subset L^{\times}$, in a similar way.
Since $\overline H_K$ is a finite index subgroup of $K^{\times}/k^{\times}$, 
$H_K'$ (resp. $H_K^{\times}$) is a finite index subgroup of $(K^{\times})^{(\Sigma)}$ 
(resp. $K^{\times}$). 
Note that 
$(k^{\times})^{\Sigma}\subset H_K'$ and 
$k^{\times}\subset H_K^{\times}$ by definition. 
Similar statements also hold for $L$. 
Moreover, 
the isomorphism 
$$\rho:\Ker (\psi_K^{(\Sigma)})\isom \Ker (\psi_L^{(\Sigma)})$$ 
restricts to an isomorphism
$$\rho:H_K'\isom H_L'.$$ 
In summary, we have the following:

\proclaim
{Lemma 4.5} (Almost-Recovering the Group of Principal Divisors).
The isomorphism $\sigma$ naturally induces isomorphisms:
$$\rho:H_K'\isom H_L'$$
and 
$$\bar \rho: \overline H_K\isom \overline H_L$$
where $\overline H_K$ (resp. $\overline H_L$) and $H'_K$ (resp. $H'_L$)
are defined as above, which fit into the following commutative diagram:
$$
\CD
H_K'@>{\rho}>>  H_L'\\
@VVV                    @VVV   \\
\overline H_K@>{\bar \rho}>> \overline H_L\\
\endCD\tag{4.3}
$$
where the vertical maps are the natural surjective homomorphisms. 
Further, $\rho$ induces naturally an isomorphism: 
$$\tau : (k^{\times})^{\Sigma}\isom (l^{\times})^{\Sigma},$$ 
which fits into the following commutative diagram:
$$
\CD
(k(x)^{\times})^{\Sigma} @>{\tau _{x,y}}>> (k(y)^{\times})^{\Sigma}\\
@AAA                    @AAA\\
(k^{\times})^{\Sigma} @>{\tau}>>   (l^{\times})^{\Sigma}\\
\endCD\tag{4.4}
$$
where $x\in X^{\cl}$, $y\defeq \phi (x)\in Y^{\cl}$, $\tau _{x,y}$ is the isomorphism in 
Lemma 2.3 (i), and the vertical maps are the natural ones.
\endproclaim

\demo
{Proof} For the last assertion, observe that  $(k^{\times})^{\Sigma}$ (resp. $(l^{\times})^{\Sigma}$)
is naturally identified with the torsion subgroup of $H_K'$ (resp. $H_L'$).
\qed
\enddemo

Given a principal divisor $\bar f \in K^{\times}/k^{\times}$ and 
$x\in X^{\cl}$, we define $v_{x}(\bar f)\in\Bbb Z$ to be the order $v_{x}(f)$ at $x$ 
of a representative $f\in K^{\times}$ of the class 
$\bar f \in K^{\times}/k^{\times}$. Thus, $v_{x}(\bar f)$ is well-defined and does 
not depend on the choice of the representative $f$ of the class $\bar f$. 
We shall refer to $v_{x}(\bar f)$ as the valuation at $x$ of the principal divisor $\bar f$.
Similarly, we define the valuation $v_y(\bar g)$ of a principal divisor 
$\bar g$ on $Y$ at a point $y\in Y^{\cl}$. 

Given a pseudo-function $f'\in (K^{\times})^{(\Sigma)}$ and 
$x\in X^{\cl}$ with $v_{x}(\bar f)=0$ where $\bar f$ is the image of $f'$ in $K^{\times}/k^{\times}$, 
we will denote by $f'(x)$ the image of $f(x)$ in $(k(x)^{\times})^{\Sigma}$, where $f\in K^{\times}$ 
is a representative of the class of $f'\in (K^{\times})^{(\Sigma)}=K^{\times}/(k^{\times}\{\Sigma'\})$, 
via the natural surjective map $k(x)^{\times}\twoheadrightarrow 
(k(x)^{\times})^{\Sigma}$. Thus, $f'(x)$ is well-defined and does not depend on the choice of the representative
$f$ of the class $f'$. We shall refer to $f'(x)$ as the $\Sigma$-value at $x$ 
of the pseudo-function $f'$. We define the $\Sigma$-value $g'(y)\in (k(y)^{\times})^{\Sigma}$
of a pseudo-function $g'\in (L^{\times})^{(\Sigma)}$ at a point $y\in Y^{\cl}$ with $v_{y}(\bar g)=0$ 
in a similar way. 

Further, for $x\in X^{\cl}$ (resp. $y\in Y^{\cl}$) we will think of elements of $(k(x)^{\times})^{\Sigma}$ 
(resp $(k(y)^{\times})^{\Sigma}$)
as classes of elements 
of $k(x)^{\times}$ (resp. $k(y)^{\times}$) modulo elements of $k(x)^{\times}\{\Sigma'\}$ 
(resp. $k(y)^{\times}\{\Sigma'\}$) 
and denote them
by $\eta', \zeta', \dots$, where $\eta, \zeta,\dots \in k(x)^{\times}$ (resp. $\in k(y)^{\times}$) 
are elements of multiplicative groups of residue fields.

\proclaim
{Lemma 4.6} (Recovering the Valuations and the $\Sigma$-Values of Pseudo-Functions) 
Consider the commutative 
diagram (4.3) 
in Lemma 4.5. Let $x\in X^{\cl}$, and $y\defeq \phi (x)$. Then the following implications hold:

\noindent 
{\rm (i)} For $\bar f\in \overline H_K$ and $\bar g\in\overline H_L$:
$$\bar \rho (\bar f)=\bar g \Longrightarrow v_{x}(\bar f)=v_{y}(\bar g).$$ 
In particular, in terms of divisors, if: 
$$\bar f=x_1+x_2+\dots+x_n-x_1'-\dots-x_{n'}',$$ 
then: 
$$\bar g= y_1+y_2+\dots+y_n-y_1'-\dots-y_{n'}',$$ 
where $y_i\defeq \phi (x_i)$ 
(resp. $y_{i'}'\defeq \phi (x_{i'}')$) for $i\in \{1,\dots,n\}$ 
(resp. $i'\in \{1,\dots,n'\}$). In other words
the map $\bar \rho$ preserves the valuations of the classes of functions in $\overline H_K$
with respect to the bijection $\phi:X^{\cl}\isom Y^{\cl}$ between points.

\noindent
{\rm (ii)} For $f'\in H_K'$ and $g'\in H_L'$:
$$\text{$v_{x}(\bar f)=0$ and $\rho(f')=g'$}
\Longrightarrow \text{$v_{y}(\bar g)=0$ and $\tau _{x,y} (f'(x))=g'(y)$},$$
where 
$$\tau _{x,y}:(k(x)^{\times})^{\Sigma} \isom (k(y)^{\times})^{\Sigma}$$ 
is the isomorphism in Lemma 2.3 (i) and $\bar f$ (resp. $\bar g$)
is the image of $f'$ (resp. $g'$) in $K^{\times}/k^{\times}$ (resp. $L^{\times}/l^{\times}$). 
In other words
the map $\rho$ preserves the $\Sigma$-values of the pseudo-functions in $H_K'$
with respect to the bijection $\phi:X^{\cl}\isom Y^{\cl}$ between points.
\endproclaim

\demo
{Proof} As shown in (the proofs of) Lemmas 4.4 and 4.5, we have the commutative diagram 
$$
\matrix
H_K'&\subset& (K^{\times})^{(\Sigma)}& \to& J_{K}^{(\Sigma)}&
\defeq&\prod_{x\in X^{\cl}}'(K_{x}^{\times})^{(\Sigma)} \\
&&&&&&\\
\downarrow&&&&\downarrow&&\\
&&&&&&\\
H_L'&\subset& (L^{\times})^{(\Sigma)}& \to& J_{L}^{(\Sigma)}&
\defeq&\prod_{y\in Y^{\cl}}'(L_{y}^{\times})^{(\Sigma)} 
\endmatrix
$$
where the vertical arrows are the isomorphisms induced by $\sigma$. More precisely, 
$H_K'\isom H_L'$ is $\rho'$, and $J_{K}^{(\Sigma)}\isom J_{L}^{(\Sigma)}$ 
maps each component 
$(K_{x}^{\times })^{(\Sigma)}$ 
isomorphically onto $(L_{y}^{\times})^{(\Sigma)}$, 
where $y\defeq \phi (x)$. Further, the isomorphism 
$(K_{x}^{\times })^{(\Sigma)}\isom 
(L_{y}^{\times})^{(\Sigma)}$ arises from Lemma 2.5. 
It follows from this that $\bar \rho$ preserves the valuations by Lemmas 2.4.3 and 2.4.9 
and that $\rho$ preserves the $\Sigma$-values by 
Lemma 2.5. 
\qed
\enddemo

Let $U$ be an open subgroup of $G_K^{(\Sigma)}$, and let $V\defeq \sigma (U)$.
Let $K'/K$ (resp. $L'/L$) be the finite subextension of $K\sptilde /K$ (resp. $L\sptilde /L$)
corresponding 
to $U$ (resp. $V$), $k'$ (resp. $l'$) the constant field of $K'$ (resp. $L'$), and $X'$ (resp. $Y'$) the normalization of $X$ (resp. $Y$) 
in $K'$ (resp. $L'$). 
Then $\sigma$ induces, by restriction to $U$, 
an isomorphism 
$$\sigma:U(=G_{K'}^{(\Sigma)})\isom V(=G_{L'}^{(\Sigma)}),$$ 
which naturally induces by Lemma 4.5 the following commutative diagram: 
$$
\CD
H_{K'}'@>>>  H_{L'}'\\
@VVV                    @VVV   \\
\overline H_{K'}@>>> \overline H_{L'}\\
\endCD\tag{4.5}
$$
where the horizontal arrows are the isomorphisms induced by $\sigma$, and 
the vertical arrows are the natural surjective homomorphisms. 

\proclaim
{Lemma 4.7} 
The above diagram (4.5) is compatible with the diagram (4.3) in Lemma 4.5. 
More precisely, the natural injective homomorphisms 
$(K^{\times})^{(\Sigma)}\to ((K')^{\times})^{(\Sigma)}$,  
$(L^{\times})^{(\Sigma)}\to ((L')^{\times})^{(\Sigma)}$
(resp.   $K^{\times}/k^{\times}\to (K')^{\times}/(k')^{\times}$, 
$L^{\times}/l^{\times}\to (L')^{\times}/(l')^{\times}$) 
map $H_K'$ into $H_{K'}'$, $H_L'$ into $H_{L'}'$
(resp. $\overline H_{K}$ into $\overline H_{K'}$, and $\overline H_{L}$ into $\overline H_{L'}$)
and the resulting diagrams 
$$
\CD
H_{K'}'@>>> H_{L'}'\\
@AAA                    @AAA   \\
H_K' @>>> H_L'\\
\endCD\tag{4.6}
$$
and 
$$
\CD
\overline H_{K'}@>>> \overline H_{L'}\\
@AAA                    @AAA   \\
\overline H_K @>>> \overline H_L\\
\endCD\tag{4.7}
$$
are commutative. 
\endproclaim

\demo
{Proof} 
First, consider the diagram 
$$
\CD
J_{K'}^{(\Sigma)} @>>>  J_{L'}^{(\Sigma)}\\
@AAA                    @AAA   \\
J_{K}^{(\Sigma)} @>>> J_{L}^{(\Sigma)}, 
\endCD
\tag {4.8}
$$
where the horizontal arrows are the natural isomorphisms induced by $\sigma$ 
and the vertical arrows are induced by the natural inclusions 
$J_K\hookrightarrow J_{K'}$, $J_L\hookrightarrow J_{L'}$ of id\`ele groups. 
This diagram is commutative, 
since the vertical arrows arise from the (local) transfer maps. 
Now, the diagram (4.6)
commutes as a subdiagram of (4.8), and the diagram (4.7) 
commutes as a quotient diagram of (4.6). 
\qed
\enddemo

{\it From now on, we shall assume that 
$\Sigma$ satisfies condition $(\epsilon_{X})$} (cf. discussion before Theorem C in $\S0$). 
Then, by Lemmas 2.4.7 (ii) and 2.4.10, 
$\Sigma$ also satisfies condition $(\epsilon_{Y})$.  
We shall use the following lemma.

\proclaim
{Lemma 4.8} Let $k\subset \kappa \subset \bar k$ be an (infinite or finite) extension 
of $k$, and $\Cal K\defeq K\kappa$. Let $\overline{U}\subset  \Cal K^{\times}/\kappa^{\times}$ be a finite index subgroup 
and assume $\sharp(\kappa)>2(\Cal K^{\times}/\kappa^{\times}: \overline{U})$. 
Then there exists $f\in \Cal K^{\times}\setminus \kappa^{\times}$, such that 
$\bar f$, $\overline {1+f}\in  \overline{U}$, and that 
$\deg (f)=\gon (X\times _{k} \kappa)$, where $\gon (X\times _{k} \kappa)$ denotes the gonality 
of $X\times _{k} \kappa$ over $\kappa$ 
and $\deg (f)$ is the degree of the finite map $f:X\times _{k} 
\kappa\to \Bbb P^1_{\kappa}$ (equivalently, $\deg(f)$ is the degree of the pole divisor 
of $f$). 
\endproclaim

\demo
{Proof} 
Take any $g\in \Cal K$ that attains the gonality: $\deg (g)=\gon (X\times _{k} \kappa)$, 
and consider the set 
$\{g-a \mid a \in \kappa\}\subset \Cal K^{\times}$. Since 
$\sharp(\kappa)>2(\Cal K^{\times}/\kappa^{\times}: \overline{U})$ by assumption, 
there exist three distinct values $a,b,c \in \kappa$ such that the images of 
$g-a,g-b,g-c$ in the quotient group 
$(\Cal K^{\times}/\kappa^{\times})/\overline{U}$ 
are the same. Now, define $U$ to be the inverse image of $\overline{U}$ 
in $\Cal K^{\times}$ and 
set 
$$f\defeq
\frac{a-b}{b-c}\cdot\frac{g-c}{g-a}\in U.$$
Then we have 
$$1+f=\frac{a-c}{b-c}\cdot\frac{g-b}{g-a}\in U.$$
Finally, as $f$ is a linear fractional transformation of $g$, 
we have 
$\deg(f)=\deg(g)=\gon (X\times _{k} \kappa)$, as desired. 
\qed
\enddemo

\proclaim
{Lemma 4.9}  Let 
$$\tau : (k^{\times})^{\Sigma} \isom (l^{\times})^{\Sigma}$$ 
be the isomorphism in Lemma 4.5 between the maximal $\Sigma$-primary quotients 
of the multiplicative groups of the constant fields, which, by Lemma 4.7, 
extends to 
$$\tau : (\bar k^{\times})^{\Sigma} \isom (\bar l^{\times})^{\Sigma}$$ 
naturally (by passing to the open subgroups of $G_K^{(\Sigma)}$ and $G_L^{(\Sigma)}$, 
corresponding to each other via $\sigma$). 
For $\eta\in \bar k^{\times}$ and $\zeta\in \bar l^{\times}$, 
if 
$$\text{$1+\eta\neq0$ and $\tau (\eta')=\zeta'$},$$  
where $\eta'$ (resp. $\zeta'$) is the image of $\eta$ (resp. $\zeta$) in $(\bar k^{\times})^{\Sigma}$
(resp. $(\bar l^{\times})^{\Sigma}$),
then there exist $\alpha, \beta \in \bar l^{\times}\{\Sigma'\}$, such that 
$$\text{$\alpha+\beta\zeta \neq 0$ and $\tau ((1+\eta)')=(\alpha+\beta \zeta)'$}.$$ 
\endproclaim

\demo
{Proof} 
Take a finite extension $k'$ of $k$ (resp. $l'$ of $l$) such that 
$\gon(X\times_{k}\bar k)=\gon(X\times_{k}k')$ (resp.
$\gon(Y\times_{l}\bar l)=\gon(Y\times_{l}l')$). By replacing $k'$ and $l'$ with  
suitable finite extensions, we may and shall assume that 
$\sigma: G_{K}^{(\Sigma)}\isom G_{L}^{(\Sigma)}$ 
induces an isomorphism $G_{Kk'}^{(\Sigma)}\isom G_{Ll'}^{(\Sigma)}$ 
(cf. Lemma 2.4.3). 
Since $\Sigma$ satisfies condition $(\epsilon_{X})$, 
there exists an extension $k''$ of $k'$ in $\bar k$, such that 
$\sharp(k'')>2\sharp (J_{X}(k'')\{\Sigma'\})$. In particular, 
$\sharp (J_{X}(k'')\{\Sigma'\})<\infty$, hence, by replacing $k''$ 
by a suitable subfield containing $k'$ if necessary, we may and shall assume 
that $k''$ is a finite extension of $k'$. Let $l''$ be the finite extension of 
$l'$ corresponding to $k''$ via $\sigma$: 
$\sigma(G_{Kk''}^{(\Sigma)})=G_{Ll''}^{(\Sigma)}$. 
Set $\Cal K\defeq Kk''$ and $\Cal L\defeq Ll''$. 
Now, by Lemma 4.8, 
there exists $f\in \Cal K^{\times}\setminus (k'')^{\times}$, such that 
$\bar f$, $\overline {1+f}\in 
\Ker(\Cal K^{\times}/(k'')^{\times}\to \pi_1(Y\times_{l}l'')^{\ab,\tor}\{\Sigma'\})$ 
and that $\deg (f)=\gon (X\times _{k} k'')=\gon(X\times_{k}\bar k)$. 
Similarly, there exists $g_1\in \Cal L^{\times}\setminus {(l'')}^{\times}$, such that 
$\bar g_1$, $\overline {1+g_1}\in 
\Ker(\Cal L^{\times}/{(l'')}^{\times}\to \pi_1(X\times_{k}k'')^{\ab,\tor}\{\Sigma'\})$ 
and that $\deg (g_1)=\gon (Y\times _{l} l'')=\gon(Y\times_{l}\bar l)$ (use again Lemma 4.8). 
Here we used the fact that $\pi_1(Y\times_{l}l'')^{\ab,\tor}\{\Sigma'\}\isom J_Y(l'')\{\Sigma'\}$
and $\pi_1(X\times_{k}k'')^{\ab,\tor}\{\Sigma'\}\isom J_X(k'')\{\Sigma'\}$.

We may write $\rho(f')=g'$ for some $g\in 
{\Cal L}^{\times}/(l'')^{\times}
$ and 
$\rho^{-1}(g_1')=f_1'$ for some $f_1\in {\Cal K}^{\times}/(k'')^{\times}$. Thus, we have 
$$\gon(X\times_{k}\bar k) 
=\deg(f)=\deg(g)\geq\gon(Y\times_{l}\bar l)$$ 
and 
$$
\gon(Y\times_{l}\bar l) 
=\deg(g_1)=\deg(f_1)\geq\gon(X\times_{k}\bar k), 
$$
where the second equalities follow from Lemmas 2.4.5 and 4.6 (i), hence 
$$n\defeq \gon(X\times_{k}\bar k)=\gon(Y\times_{l}\bar l).$$
Further, again by replacing $k''$ and $l''$ by suitable finite extensions 
corresponding to each other, we may assume that the zeros ($\subset X^{\cl}$) of 
$1+f$ are $k''$-rational and that $\eta\in (k'')^{\times}$ and $\zeta\in (l'')^{\times}$. 

{}From now on, we may and shall assume that $k''=k$, and $l''=l$, by replacing 
$K$ and $L$ by $\Cal K$ and $\Cal L$, respectively, and $\sigma:G_K^{(\Sigma)}\isom G_L^{(\Sigma)}$ 
by the isomorphism $G_{\Cal K}^{(\Sigma)}\isom G_{\Cal L}^{(\Sigma)}$ induced by $\sigma$. 
Thus, 
$f\in H_K^{\times}\subset K^{\times}$, 
$g\in H_L^{\times}\subset L^{\times}$, 
$\rho(f')=g'$ and $\bar\rho(\bar f)=\bar g$. 

\definition {Claim} We have 
$$\rho ((1+f)')=(\alpha+\beta g)'$$ 
for some $\alpha, \beta \in \bar l^{\times}\{\Sigma'\}$. 
\enddefinition

Indeed, as $\overline {1+f}\in  \overline H_K$, we have $(1+f)'\in H_K'$, hence 
we may write $\rho ((1+f)')=h'$ for some $h\in H_L^{\times}\subset L^{\times}$. 
Write $\overline {1+f}=x_1+\cdots+x_n-(f)_{\infty}$ as a divisor, 
where $(f)_{\infty}$ denotes the pole divisor of $f$ which is equal to the pole divisor of 
$1+f$. Note that $x_1,\dots,x_n$ are $k$-rational by our choice. 
Then $\bar h=y_1+\cdots+y_n-(g)_{\infty}$, where $y_i\defeq \phi (x_i)$ 
(observe that $g$ and $h$ have the same pole divisor by Lemma  4.6 (i), since $f$ and $1+f$ do). 
Note that $y_1,\dots,y_n$ are $l$-rational by Lemma 2.4.5. 
Let $c \defeq g(y_1)\in l^{\times}$. Note that 
$-c\in l^{\times}\{\Sigma'\}$ by the preservation of the $\Sigma$-value of pseudo-functions 
(cf. Lemma 4.6 (ii)) and the fact that $f(x_1)=-1$. 
(Observe that $\tau_{x,y}(-1)=-1$.) 
Thus, $g-c$ has a zero at $y_1$ and $\overline {g-c}=y_1+E-(g)_{\infty}$ as a divisor, where 
$E$ is an effective divisor of degree $n-1$. 
Consider the function $h_1\defeq h/(g-c)$. Then $\overline {h_1}=y_2+\cdots+y_n-E$ as a divisor. 
Thus, $\deg (h_1)<n$, which implies that $h_1=\beta\in l$ is a constant (by the minimality of $n$ 
as the 
degree of a non-constant function), and $h=\beta g+\beta(-c)$. Further, 
let $w\in X^{\cl}$ be a zero of $f$ and set $z\defeq\phi(w)\in Y^{\cl}$. 
Then $z$ is a zero of $g$ by Lemma 4.6 (i), and 
$$\beta=h_1(z)=h(z)/(-c)\in \bar l^{\times}\{\Sigma'\}$$
by Lemma 4.6 (ii) and the fact that $(1+f)(w)=1$. Thus, $\alpha\defeq\beta(-c)\in 
\bar l^{\times}\{\Sigma'\}$, and 
the above claim is proved.

Let $f\in H_K^{\times}$ and $g\in H_L^{\times}$ be as above. 
In particular, 
$\rho(f')=g'$ and $\rho((1+f)')=(\alpha+\beta g)'$ for some 
$\alpha, \beta \in \bar l^{\times}\{\Sigma'\}$. 
Now, let $\eta \in k^{\times}$, $\zeta\in l^{\times}$ such 
that $1+\eta\neq 0$ and $\tau (\eta')=\zeta'$. 
Let $x\in X^{\cl}$
be a zero of $f-\eta$ and set 
$y\defeq \phi (x)$. We have $f\equiv \eta\ (\mod x)$ which implies that  $1+f\equiv 1+\eta\ (\mod x)$,
and $g\equiv \xi\ (\mod y)$ where $\zeta' =\xi'\in (\bar l^{\times})^{\Sigma}$ by
the preservation of the $\Sigma$-values of pseudo-functions (cf. Lemma 4.6 (ii)), i.e. 
there exists $\epsilon\in \bar l^{\times}\{\Sigma'\}$ such that $\xi=\epsilon\zeta$.
Then 
$$\tau ((1+\eta)')=\tau((1+f)'(x))=\rho((1+f)')(y)=(\alpha+\beta g)'(y)=(\alpha +\beta \xi)'
=(\alpha+\beta \epsilon\zeta)',$$ 
where the second equality
results from the preservation of the $\Sigma$-values of pseudo-functions (cf. Lemma 4.6 (ii)). 
As $\alpha, \beta \epsilon\in \bar l^{\times}\{\Sigma'\}$, the assertion follows. 
\qed
\enddemo

\proclaim
{Lemma 4.10} The isomorphisms 
$$\tau _{x,y}:(k(x)^{\times})^{\Sigma} \isom (k(y)^{\times})^{\Sigma},$$
in Lemma 2.3 (i) satisfy the following property: 
For $\eta\in k(x)^{\times}$ and $\zeta\in k(y)^{\times}$, 
if 
$$\text{$1+\eta\neq0$ and $\tau_{x,y} (\eta')=\zeta'$},$$  
then there exist $\alpha, \beta \in \bar l^{\times}\{\Sigma'\}$, such that 
$$\text{$\alpha+\beta\zeta \neq 0$ and $\tau_{x,y} ((1+\eta)')=(\alpha+\beta \zeta)'$}.$$ 
\endproclaim

\demo
{Proof} After passing to finite extensions of scalars, this follows directly from 
Lemma 4.9 and the commutativity of the diagram (4.4) in Lemma 4.5. 
\qed
\enddemo

Next, recall the definition of $H_K^{\times}$:
$$H_K^{\times}\defeq \{f\in K^{\times}\ \vert\ \bar f\in \overline H_K\}.$$
Then $H_K^{\times}$ is a finite index subgroup of  $K^{\times}$, and the (finite) quotient
$K^{\times}/ H_K^{\times}$ is killed by $m\defeq \sharp(\pi _1(X)^{\ab,\tor}\{\Sigma'\})
=\sharp(\pi _1(Y)^{\ab,\tor}\{\Sigma'\})$, hence is $\Sigma'$-primary. 

\proclaim
{Lemma 4.11} Let $f\in H_K^{\times}$, and assume that $1+f\neq 0$. Then $1+f\in H_K^{\times}$.
\endproclaim

\demo
{Proof} 
Write $\rho(f')=g'$ with $g\in H_L^{\times}$. 
First, we have 
$((1+f)')^m\in H_K'$. 
Thus, we may write $\rho(((1+f)')^m)=h'$ with $h\in H_L^{\times}$. 
Next,  let $x\in X^{\cl}$ such that 
$x$ is 
neither a pole of $f$ nor a zero of $1+f$ 
and set $y\defeq \phi (x)$. 
Then we have 
$$
\align
&h'(y)=
\rho(((1+f)')^m)(y)\\
=&\tau _{x,y}(((1+f)')^m(x))
=\tau _{x,y}((1+f(x))')^{m} \\
=&((\alpha_{y}+\beta _{y} g(y))')^m
\endalign
$$
for some $\alpha_{y}, \beta_{y}\in \bar l^{\times}\{\Sigma'\}$, 
by Lemma 4.6 (ii) and Lemma 4.10. Equivalently, for some 
$\alpha_{y}, \beta_{y}, \gamma_{y}\in \bar l^{\times}\{\Sigma'\}$, 
we have 
$$h(y)=\gamma_{y}(\alpha_{y}+\beta _{y} g(y))^m.$$
Thus, we have 
$$h(y)=a_{y}(1+c_{y} g(y))^m,$$
where $a_{y}\defeq\gamma_{y}\alpha_{y}^m \in \bar l^{\times}\{\Sigma'\}$ 
and $c_{y}\defeq \beta_{y}/\alpha_{y}\in \bar l^{\times}\{\Sigma'\}$. 
By Definition/Proposition 3.12 (i)(ii), 
this implies that
$h=a(1+cg)^m$ 
for some $a, c\in l^{\times}\{\Sigma'\}$. 
Accordingly, $\rho(((1+f)')^m)=h'=((1+cg)')^m$ in $H_Y'\subset (L^{\times})^{(\Sigma)}$, 
hence $\rho((1+f)')^m=((1+cg)')^m$ in $\Ker(\psi_L^{(\Sigma)})$, where 
$\psi_L^{(\Sigma)}: J_{L}^{(\Sigma)}\to G_Y^{(\Sigma),\ab}$ is naturally induced by 
Artin's reciprocity map in global class field theory. 

Now, 
since $\Ker (\psi_L^{(\Sigma)})\subset J_{L}^{(\Sigma)}
= \prod_{y\in Y^{\cl}}'(L_{y}^{\times})^{(\Sigma)}$ does not 
admits a nontrivial $\Sigma'$-primary torsion, we conclude 
$\rho((1+f)')=(1+cg)'$ in $\Ker(\psi_L^{(\Sigma)})$. As 
$(1+cg)'\in (L^{\times})^{(\Sigma)}$, we have $(1+f)'\in H_K'$ by definition, 
as desired. 
\qed
\enddemo

We set
$$H_K\defeq H_K^{\times}\cup \{0\}\subset K,$$
and
$$H_L\defeq H_L^{\times}\cup \{0\}\subset L.$$

\proclaim
{Lemma 4.12} 
{\rm (i)} 
The subset $H_K$ of $K$ is a subfield. 

\noindent
{\rm (ii)} We have 
$H_K=K$, $H_L=L$, $H_K'=(K^{\times})^{(\Sigma)}$, 
$H_L'=(L^{\times})^{(\Sigma)}$,
$\overline H_K=K^{\times}/k^{\times}$,
and 
$\overline H_L=L^{\times}/l^{\times}$.

\endproclaim

\demo
{Proof} (i) First note that $k\subset H_K$, and
$H_K$ is closed under multiplication by its definition. Also, $H_K$ is closed under
addition. Indeed, let us show $f+g\in H_K$ for any $f,g\in H_K$. This is clear 
if either one of $f,g,f+g$ is zero. So, assume that none of $f,g,f+g$ is zero. 
Then, as $\frac{g}{f}\in H_K^\times$ and $1+\frac {g}{f}=\frac{f+g}{f}\neq 0$, 
Lemma 4.11 implies that $1+\frac {g}{f}\in H_K^{\times}$, 
hence $f+g=f(1+\frac {g}{f})\in H_K^{\times}\subset H_K$. 
Thus, $H_K$ is a $k$-subfield of $K$.

\noindent
(ii) Since $H_K^{\times}$ is a finite index subgroup of $K^{\times}$, there exist 
finitely many $f_1,\dots, f_r\in K^{\times}$ such that 
$K=H_Kf_1\cup\cdots\cup H_Kf_r$ holds. In particular, 
$K=H_Kf_1+\cdots+H_Kf_r$, hence 
$H_K\subset K$ is a finite field extension. This implies that 
$H_K$ is an infinite field. Then $K$ cannot be covered by finitely many proper 
$H_K$-vector subspaces. 
Thus, the equality $K=H_Kf_1\cup\cdots\cup H_Kf_r$ implies that $K$ is $1$-dimensional 
over $H_K$, i.e., $K=H_K$, as desired. In particular, 
$H_K'=(K^{\times})^{(\Sigma)}$ and $\overline H_K=K^{\times}/k^{\times}$. 

As $((K^{\times})^{(\Sigma)}: H_K')=(\Ker(\psi_K^{(\Sigma)}): H_K')/\sharp(\pi_1(X)^{\ab,\tor}\{\Sigma'\})$ 
and
$((L^{\times})^{(\Sigma)}: H_L')=(\Ker(\psi_L^{(\Sigma)}): H_L')/\sharp(\pi_1(Y)^{\ab,\tor}\{\Sigma'\})$, 
we have  
$((K^{\times})^{(\Sigma)}: H_K')=((L^{\times})^{(\Sigma)}: H_L')$. Thus, 
$H_L'=L^{\times}/{\Sigma}$ also holds, from which $H_L=L$ and 
$\overline H_L=L^{\times}/l^{\times}$ follow. 
\qed
\enddemo

It follows from Lemma 4.12 above that $\bar \rho$ is an isomorphism:
$$\bar \rho: K^{\times}/k^{\times}\isom  L^{\times}/l^{\times}$$
which is naturally induced by $\sigma$.

Next, we will think of elements of  $K^{\times}/k^{\times}$ (resp. $L^{\times}/l^{\times}$) 
as points of the infinite-dimensional projective space over $k$ (resp. $l$)
associated to the vector space $K$ (resp. $L$) over $k$ (resp. $l$). 
In particular, points of this 
projective space correspond to one-dimensional 
$k$-linear (resp. $l$-linear) 
subspaces in $K$ (resp. $L$), and lines correspond to two-dimensional 
$k$-linear (resp. $l$-linear) 
subspaces of $K$ (resp. $L$).  

\proclaim
{Lemma 4.13} (Recovering the Additive Structure of 
Function Fields)
The natural isomorphism 
$\bar \rho :K^{\times}/k^{\times}\isom L^{\times}/l^{\times}$ which follows from 
Lemmas 4.5 and 4.12, viewed as a set-theoretic bijection between points of projective spaces, 
preserves colineations. 
Accordingly, $\bar \rho$ arises from a $\psi _0$-isomorphism 
$$\psi: (K,+) \isom  (L,+),$$
where  $\psi _0:k\isom l$ is a field isomorphism. Namely, 
$\psi$ is an isomorphism of abelian groups which is compatible with $\psi _0$ 
in the sense that $\psi(ax)=\psi_0(a)\psi(x)$ for $a\in k$ and $x\in K$. 
Further, $\psi_0$ is uniquely determined and $\psi$ is uniquely determined 
up to scalar multiplication. 
\endproclaim

\demo
{Proof} 
In order to show the first assertion that the map $\bar \rho$ preserves colineations, it suffices to show that
for a non-constant function $f\in K^{\times}\setminus k^{\times}$, 
if $\bar \rho (\bar f)=\bar g$, then 
$\bar \rho (\overline {1+f})=\overline  {\alpha+\beta g}$, where $\alpha,\beta \in l$.
By replacing $g\in L^{\times}$ if necessary, we may and shall assume that 
$\rho(f')=g'$ holds. 
Write $\rho ((1+f)')=h'$ with $h\in L^{\times}$. Let $x\in X^{\cl}$ with 
$f(x)\not\in \{\infty, 0, -1\}$, 
and set $y\defeq \phi (x)$.  
Then $\tau _{x,y} (f'(x))=g'(y)$ and $\tau _{x,y} ((1+f)'(x))=h'(y)$ by Lemma 4.6 (ii).
Let $\eta \defeq f(x)$ and $\zeta \defeq g(y)$. 
Then 
$\tau _{x,y}(\eta')=\zeta'$. 
But $\tau _{x,y}((1+\eta)')=(\alpha_{y}+\beta_{y} \zeta)'$ by Lemma 4.10,
where $\alpha_{y},\beta_{y}\in \bar l^{\times}\{\Sigma'\}$.
Thus, $h(y)=a_{y}+b_{y}g(y)$, where $a_{y},b_{y}\in \bar l^{\times}\{\Sigma'\}$, 
$\forall 'x\in X^{\cl}$.
But this implies that $h=a+bg$ for some $a,b\in l^{\times}\{\Sigma'\}$ by Proposition 3.11 
and Definition/Proposition 3.5, as required.

The second and the third assertions follow from the first assertion and 
the fundamental theorem of projective geometry (cf. [Artin]). 
\qed
\enddemo

\proclaim
{Lemma 4.14} (Recovering the Field Structure of Function Fields) 
If we normalize the isomorphism 
$$\psi: (K,+) \isom  (L,+)$$
in Lemma 4.13 by the condition $\psi(1)=1$, it becomes a field isomorphism such that 
the diagram
$$
\CD
K @>\psi>>  L\\
@AAA                    @AAA   \\
k @>\psi_0>> l 
\endCD
$$
commutes. 
\endproclaim

\demo
{Proof} (See also the end of the proof of Theorem 5.11 in [Pop2].) 
Take any $f\in K^{\times}$, then $\psi\circ\mu_f$ and $\mu_{\psi(f)}\circ\psi$ 
are $\psi_0$-isomorphisms $(K,+)\isom (L,+)$, where $\mu_g$ denotes the 
$g$-multiplication map. The isomorphisms 
$K^{\times}/k^{\times}\isom L^{\times}/l^{\times}$ 
they induce coincide with each other: 
$$\overline{\psi\circ\mu_f}=\bar\rho\circ\mu_{\bar f}=\mu_{\bar\rho(\bar f)}
\circ \bar\rho=\overline{\mu_{\psi(f)}\circ\psi},$$
where the second equality follows from the multiplicativity of $\bar\rho$. 
Further, we have 
$$\psi\circ\mu_f(1)=\psi(f)=\mu_{\psi(f)}(1)=\mu_{\psi(f)}\circ\psi(1).$$ 
Thus, the equality $\psi\circ\mu_f=\mu_{\psi(f)}\circ\psi$ 
follows from the uniqueness in the fundamental theorem of projective geometry, 
which shows the multiplicativity of $\psi$. 

For any $a\in k$, we have 
$$\psi(a)=\psi(a\cdot 1)=\psi_0(a)\psi(1)=\psi_0(a)\cdot 1=\psi_0(a),$$
which shows the commutativity of the diagram. 
\qed
\enddemo

Let $U$ be an open subgroup of $G_K^{(\Sigma)}$, and let $V\defeq \sigma (U)$.
Let $K'/K$ (resp. $L'/L$) be the subextension of $K\sptilde /K$ 
(resp. $L\sptilde /L$) corresponding to $U$ (resp. $V$). 
Then $\sigma$ induces, by restriction to $U$, 
an isomorphism 
$$\sigma:U\isom V.$$ 
Note that it is unclear in general if $\Sigma$ satisfies condition $(\epsilon_{X'})$, 
where $X'$ 
denotes the normalization of $X$  
in $K'$.  
However, this condition 
is only used to establish Lemma 4.9 by resorting to Lemma 4.8. Since the assertion of 
Lemma 4.9 for $\sigma:U(=G_{K'}^{(\Sigma)})\isom V(=G_{L'}^{(\Sigma)})$ 
is just the same as that for $\sigma: G_K^{(\Sigma)}\isom G_L^{(\Sigma)}$, we can deduce from $\sigma:U\isom V$,
by Lemma 4.14 (without the need to assume that $\Sigma$ satisfies condition $(\epsilon_{X'})$), 
a natural field isomorphism:
$$\psi': K'\isom L'.$$

\proclaim
{Lemma 4.15} 
{\rm (i)} The following diagram is commutative: 
$$
\CD
K' @>{\psi'}>>  L'\\
@AAA                    @AAA   \\
K @>{\psi}>> L
\endCD
$$
where the vertical arrows are the natural inclusions and $\psi$, $\psi'$ are 
the field isomorphisms induced by $\sigma$. 

\noindent
{\rm (ii)} If, moreover, 
$U$ is normal in $G_K^{(\Sigma)}$, then $V$ is normal in $G_L^{(\Sigma)}$ 
and the above diagram is Galois-equivariant with respect to the 
isomorphism $G_K^{(\Sigma)}/U\isom G_L^{(\Sigma)}/V$ induced by $\sigma$. 
\endproclaim

\demo
{Proof} (i) Let $k'$ (resp. $l'$) denote the constant field of $K'$ (resp. $L'$). 
Then, the commutativity of the diagram (4.7) in Lemma 4.7 implies that 
the diagram 
$$
\CD
(K')^{\times}/(k')^{\times} @>{\bar\rho'}>>  (L')^{\times}/(l')^{\times}\\
@AAA                    @AAA   \\
K^{\times}/k^{\times} @>{\bar\rho}>> L^{\times}/l^{\times}
\endCD
\tag{4.9}
$$
commutes. 

Now, write $i:K\to K'$ and $j:L\to L'$ for the natural inclusions. 
To prove $\psi'\circ i=j\circ \psi$, we shall first show that the image 
$\psi'(K)$ of the left-hand side map $\psi'\circ i$ 
and the image $j(\psi(K))=L$ of the right-hand side map $j\circ \psi$ 
coincide with each other. But by the commutativity of (4.9), we have 
at least: 
$\psi'(K^{\times})\cdot (l')^{\times}=L^{\times}\cdot(l')^{\times}$. 
Set $H\defeq\psi'(K)\cap L$, which is a subfield of $\psi'(K)$ 
and a subfield of $L$ at a time. 
Further, 
$H^{\times}=
\psi'(K^{\times})\cap L^{\times}$ is of finite index 
(dividing $\sharp((l')^{\times})$) both in $\psi'(K^{\times})$ and 
in $L^{\times}$. 
Thus, as in the proof of Lemma 4.12, we deduce $\psi'(K)=H=L$. 

Finally, the desired equality $\psi'\circ i=j\circ \psi$ follows from 
the uniqueness in the fundamental theorem of projective geometry, since 
the diagram (4.9) commutes and 
$\psi'\circ i(1)=1=j\circ \psi(1)$. 

\noindent
(ii) Assume that $U$ is normal in $G_{K}^{(\Sigma)}$, then 
$V=\sigma(U)$ is normal in $G_{L}^{(\Sigma)}$ and 
$\sigma$ induces an isomorphism 
$G_{K}^{(\Sigma)}/U\isom G_{L}^{(\Sigma)}/V$. 

Since the action of $G_{K}^{(\Sigma)}/U$ 
(resp. $G_{L}^{(\Sigma)}/V$) on $J_{K'}^{(\Sigma)}$ 
(resp. $J_{L'}^{(\Sigma)}$)
arises from the conjugation on the decomposition groups, the isomorphism 
$J_{K'}^{(\Sigma)}\isom J_{L'}^{(\Sigma)}$ is Galois-equivariant. 

Further, since the diagram 
$$
\matrix
(K')^{\times}/(k')^{\times}&\twoheadleftarrow&((K')^{\times})^{(\Sigma)}&
\hookrightarrow&J_{K'}^{(\Sigma)}\\
&&&&\\
\phantom{\wr}\downarrow\wr&&\phantom{\wr}\downarrow\wr&&\downarrow\wr\\
&&&&\\
(L')^{\times}/(l')^{\times}&\twoheadleftarrow&((L')^{\times})^{(\Sigma)}&
\hookrightarrow&J_{L'}^{(\Sigma)}
\endmatrix
$$
is commutative, the isomorphisms 
$\rho': ((K')^{\times})^{(\Sigma)}\isom ((L')^{\times})^{(\Sigma)}$ 
and 
$\bar\rho': (K')^{\times}/(k')^{\times}\isom (L')^{\times}/(l')^{\times}$ 
are Galois-equivariant. 

Finally, it follows from the uniqueness in the fundamental theorem of projective geometry that 
the isomorphism $\psi': K'\isom L'$ is Galois-equivariant, as desired. 
\qed
\enddemo

By considering various open subgroups of $G_K^{(\Sigma)}$ and $G_L^{(\Sigma)}$ as above, 
corresponding to each other via $\sigma$, 
and using Lemmas 4.14 and 4.15 (i), we obtain naturally a field isomorphism 
$$\tilde \psi : K\sptilde \isom L\sptilde .$$

\proclaim
{Lemma 4.16} The following diagram is commutative: 
$$
\CD
K\sptilde @>{\tilde \psi}>> L\sptilde \\
@AAA              @AAA \\
K @>{\psi}>> L \\
\endCD
$$
where the vertical arrows are the natural inclusions and $\psi$, $\tilde \psi$ are 
field isomorphisms induced by $\sigma$. Further, $\tilde\psi$ is 
Galois-equivariant with respect to $\sigma: G_K^{(\Sigma)}\isom G_L^{(\Sigma)}$. 
\endproclaim
\demo
{Proof} This follows directly from Lemma 4.15. 
\qed
\enddemo

Thus, from Lemma 4.16 we deduce a commutative diagram
$$
\CD
L\sptilde @>{{\tilde \psi}^{-1}}>> K\sptilde \\
@AAA              @AAA \\
L @>{\psi^{-1}}>> K \\
\endCD
$$
which is Galois-equivariant with respect to $\sigma: G_K^{(\Sigma)}\isom G_L^{(\Sigma)}$. 
This completes the proof of the existence part of Theorem 4.1. For the uniqueness part 
of Theorem 4.1, suppose that for $j=1,2$ there is a commutative diagram 
$$
\CD
L\sptilde @>{{\tilde \psi_j}^{-1}}>> K\sptilde \\
@AAA              @AAA \\
L @>{{\psi_j}^{-1}}>> K \\
\endCD
$$
as above, where the horizontal maps are isomorphisms, and which is Galois-equivariant with respect to $\sigma: G_{K}^{(\Sigma)}\isom G_{L}^{(\Sigma)}$. 
Set $\tilde\alpha\defeq\tilde\psi_2^{-1}\circ\tilde\psi_1\in\Aut(K\sptilde)$ and 
$\alpha\defeq\psi_2^{-1}\circ\psi_1\in\Aut(K)$. Then they fit into the following commutative diagram 
$$
\CD
K\sptilde @>{\tilde \alpha}>> K\sptilde \\
@AAA              @AAA \\
K @>{\alpha}>> K \\
\endCD
$$
which is Galois-equivariant with respect to $\id: G_{K}^{(\Sigma)}\isom G_{K}^{(\Sigma)}$. 
Namely, $\tilde\alpha$ commutes with $G_{K}^{(\Sigma)}
(=\Gal(K\sptilde/K))$ in $\Aut(K\sptilde)$, or, equivalently, 
the conjugation action of $\tilde\alpha$ on $G_{K}^{(\Sigma)}$ is trivial. Then, in particular, 
every finite Galois extension $K\subset K'\subset K\sptilde$ is preserved by $\tilde\alpha$. 
Further, considering the action of $\tilde\alpha$ on 
$J_{K'}^{(\Sigma)}\hookleftarrow ((K')^{\times})^{(\Sigma)}\twoheadrightarrow 
(K')^{\times}/(k')^{\times}$, we conclude that the action of $\tilde\alpha$ on 
$(K')^{\times}/(k')^{\times}$ is trivial. Now, it follows from the uniqueness in the 
fundamental theorem of projective geometry that the action of $\tilde\alpha$ on 
$(K')^{\times}$ is trivial. Since $K\subset {K'}\subset K\sptilde$ is an arbitrary finite 
Galois extension, we conclude that the action of $\tilde\alpha$ on $K\sptilde$ 
is trivial, i.e., $\tilde\alpha=1$, as desired. This completes the proof of Theorem 4.1. 
\qed

\bigskip
\bigskip
$$\text{References.}$$
\noindent
[Artin] Artin, E., Geometric algebra, Interscience Publishers, Inc., New York, 1957.

\noindent
[Bogomolov] Bogomolov, F.A., 
On two conjectures in birational algebraic geometry, 
in Algebraic geometry and analytic geometry (Tokyo, 1990), 26--52,
ICM-90 Satell. Conf. Proc., Springer, Tokyo, 1991. 




\noindent
[Grothendieck] Grothendieck, A., Brief an G. Faltings, 
with an English 
translation on pp. 285--293, London Math. Soc. Lecture Note Ser., 242, Geometric Galois 
actions, 1, 49--58, Cambridge Univ. Press, Cambridge, 1997.

\noindent
[Grunewald-Segal] Grunewald, F. J. and Segal, D., 
On congruence topologies in number fields, 
J. Reine Angew. Math., 311/312 (1979), 389--396. 

\noindent
[Harbater] Harbater, D., Fundamental groups and embedding problems in 
characteristic $p$, in Recent developments in the inverse Galois problem 
(Seattle, WA, 1993), 353--369, Contemp. Math., 186, Amer. 
Math. Soc., Providence, RI, 1995.



\noindent
[Mochizuki] Mochizuki, S., Absolute anabelian cuspidalizations of proper 
hyperbolic curves, 
J. Math. Kyoto Univ. 47 (2007), no. 3, 451--539. 

\noindent
[Neukirch] Neukirch, J., Kennzeichnung der endlich-algebraischen Zahlk\"orper 
durch die Galoisgruppe der maximal aufl\"osbaren Erweiterungen. 
J. Reine Angew. Math.  238  (1969), 135--147.

\noindent
[Neukirch-Schmidt-Wingberg] Neukirch, J., Schmidt, A., Winberg, K., Cohomology of number fields, first edition,
Springer, Grundlehren der mathematischen Wissenschaften Bd. 323, 2000.

\noindent
[Pop1] Pop, F., \'Etale Galois covers of affine smooth curves, 
The geometric case of a conjecture of Shafarevich, On Abhyankar's 
conjecture, Invent. Math.  120  (1995),  no. 3, 555--578.


\noindent
[Pop2] Pop, F., The birational anabelian conjecture -- revisited --, manuscript (2002), 
available at http://www.math.upenn.edu/$\sim$pop/Research/Papers.html. 

\noindent
[Raynaud] Raynaud, M., Sur le groupe fondamental d'une courbe compl\`ete 
en caract\'eristique $p>0$, 
in Arithmetic fundamental groups and noncommutative algebra 
(Berkeley, CA, 1999), 335--351, Proc. Sympos. Pure Math., 
70, Amer. Math. Soc., Providence, RI, 
2002.

\noindent
[Sa\"\i di-Tamagawa1] Sa\"\i di, M. and Tamagawa, A., 
A prime-to-$p$ version of Grothendieck's anabelian conjecture for hyperbolic curves 
over finite fields of characteristic $p>0$,  Publ. Res. Inst. Math. Sci. 45 (2009), no. 1, 135--186.

\noindent
[Sa\"\i di-Tamagawa2] Sa\"\i di, M. and Tamagawa, A., On the anabelian geometry of 
hyperbolic curves over finite fields, in Algebraic number theory and related topics 2007, 
67--89, RIMS K\^oky\^uroku Bessatsu, B12, Res. Inst. Math. Sci. (RIMS), Kyoto, 2009. 

\noindent
[Sa\"\i di-Tamagawa3] Sa\"\i di, M. and Tamagawa, A., On the Hom-form of Grothendieck's 
birational anabelian conjecture in characteristic $p>0$, Algebra \&\  Number Theory 
5 (2011), no. 2, 131--184. 

\noindent
[Serre1] 
Serre, J.-P., Cohomologie galoisienne, Seconde \'edition, 
Lecture Notes in Mathematics, 5, Springer-Verlag, Berlin-Heidelberg-New York, 1962/1963.

\noindent
[Serre2] Serre, J.-P., 
Corps locaux, Publications de l'Institut de Math\'ematique de l'Universit\'e de Nancago, VIII,  
Actualit\'es Sci. Indust., No. 1296, Hermann, Paris, 1962. 

\noindent
[Serre3] Serre. J.-P., Local class field theory, in Algebraic Number 
Theory (Proc. Instructional Conf., Brighton, 1965), 128--161, Thompson, 
Washington, D.C., 1967.

\noindent
[Tamagawa] Tamagawa, A., The Grothendieck conjecture for affine curves, 
Compositio Math. 109 (1997), 135--194. 


\noindent
[Uchida] Uchida, K., Isomorphisms of Galois groups of algebraic function fields.  
Ann. Math. (2)  106  (1977), no. 3, 589--598. 

\bigskip
\noindent
Mohamed Sa\"\i di

\noindent
College of Engineering, Mathematics and Physical Sciences

\noindent
University of Exeter

\noindent
Harrison Building

\noindent
North Park Road

\noindent
EXETER EX4 4QF 

\noindent
United Kingdom

\noindent
M.Saidi\@exeter.ac.uk

\bigskip
\noindent
Akio Tamagawa

\noindent
Research Institute for Mathematical Sciences

\noindent
Kyoto University

\noindent
KYOTO 606-8502

\noindent
Japan

\noindent
tamagawa\@kurims.kyoto-u.ac.jp

\enddocument